\documentclass[preprints,article,accept,moreauthors,pdftex]{Definitions/mdpi}

\firstpage{1}
\pubvolume{1}
\issuenum{1}
\articlenumber{0}
\pubyear{2021}
\copyrightyear{2021}
\externaleditor{Academic Editor: {Francesco Mainardi and Jaan Janno} 
} 
\datereceived{24 November 2021}
\dateaccepted{12 December 2021}
\datepublished{}
\hreflink{https://doi.org/10.3390/math9243255} 

\usepackage{graphicx}
\usepackage{hyperref}
\usepackage{lineno}

\usepackage{float}
\Title{The Integral Mittag-Leffler, Whittaker and Wright Functions}

\TitleCitation{The Integral Mittag-Leffler, Whittaker and Wright Functions}

\Author{Alexander Apelblat 
 $^{1}$ and Juan Luis Gonz\'{a}lez-Santander $^{2,}$*\orcidA{}}

\AuthorNames{Firstname Lastname, Firstname Lastname and Firstname Lastname}

\AuthorCitation{Apelblat, A.; Gonz\'{a}lez-Santander, J.L.}

\address{%
$^{1}$ \quad Department of Chemical Engineering, Ben Gurion University of the Negev, Beer Sheva 84105, Israel; apelblat@bgu.ac.il\\
$^{2}$ \quad Department of Mathematics, Universidad de Oviedo, 
33007 Oviedo, Spain}

\corres{Correspondence: gonzalezmarjuan@uniovi.es}

\abstract{Integral Mittag-Leffler, Whittaker and Wright functions with integrands
similar to those which already exist in mathematical literature are
introduced for the first time. For particular values of parameters, they can
be presented in closed-form. In most reported cases, these new integral
functions are expressed as generalized hypergeometric functions but also
in terms of elementary and special functions. The behavior of some of the new
integral functions is presented in graphical form. By using the MATHEMATICA
program to obtain infinite sums that define the Mittag-Leffler, Whittaker, and
Wright functions and also their corresponding integral functions, these
functions and many new Laplace transforms of them are also reported in
the Appendices for integral and fractional values of~parameters.
}


\begin{document}
\section{Introduction}

The appearance of special functions of mathematical physics was associated with
solutions of particular ordinary differential equations, while the integral
special functions arrived much later in mathematical literature after
properties of these functions were investigated. Integral special functions
were introduced as new special functions, which can be applied in many
circumstances, especially in operational calculus, where they are frequently
serving as direct and inverse integral transforms. The form of an integrand is
identical for all integral functions, but limits of integration are
different in order to assure the convergence of defined integrals.
There are two types of integral special functions: those with elementary functions in their integrands and those with special functions.
To the first group belong the exponential integral $-\mathrm{Ei}\left( -x\right) $, the sine
and cosine integrals, $\mathrm{si}\left( x\right) $, $\mathrm{Si}\left(
x\right) $, $\mathrm{ci}(x)$ and $\mathrm{Ci}\left( x\right) $, and the
corresponding integrals of hyperbolic trigonometric functions, $\mathrm{Shi}%
\left( x\right) $ and $\mathrm{Chi}\left( x\right) $. These functions are
defined in the following way \cite{Ederlyi,Ederlyi2,Abramowitz,Magnus,DLMF}%
\begin{equation}
\begin{array}{l}
\displaystyle%
\mathrm{E}_{1}\left( x\right) =-\mathrm{Ei}\left( -x\right)
=\int_{x}^{\infty }\frac{e^{-t}}{t}dt,\quad x>0, \\
\displaystyle%
\mathrm{Si}\left( x\right) =\int_{0}^{x}\frac{\sin t}{t}dt, \\
\displaystyle%
\mathrm{si}\left( x\right) =-\int_{x}^{\infty }\frac{\sin t}{t}dt=\mathrm{Si}%
\left( x\right) -\frac{\pi }{2}, \\
\displaystyle%
\mathrm{Ci}\left( x\right) =-\int_{x}^{\infty }\frac{\cos t}{t}dt=\gamma +%
\mathrm{\ln }x-\int_{0}^{x}\frac{1-\cos t}{t}dt=-\mathrm{ci}(x), \\
\displaystyle%
\mathrm{Shi}\left( x\right) =\int_{0}^{x}\frac{\sinh t}{t}dt, \\
\displaystyle%
\mathrm{Chi}\left( x\right) =\gamma +\mathrm{\ln }x-\int_{0}^{x}\frac{%
1-\cosh t}{t}dt,%
\end{array}
\label{Trigonometric_integral_functions}
\end{equation}%
where $\gamma $ is the Euler--Mascheroni constant. As can be observed in (\ref%
{Trigonometric_integral_functions}), the integral special functions have
integrands in the form, $f\left( t\right) /t$, and the intervals of
integrations are $0<t<x$ or $x<t<\infty $. Few direct and inverse integral
transforms are presented below to illustrate their applications, for example,
in the Laplace transformation \cite{Roberts,Oberhettinger,ApelblatBook},%
\begin{equation}
F\left( s\right) :=\mathcal{L}\left[ f\left( t\right) \right]
:=\int_{0}^{\infty }e^{-st}f\left( t\right) dt,
\label{Laplace_transform_def}
\end{equation}%
we have%
\begin{equation}
\begin{array}{l}
\displaystyle%
\mathcal{L}\left[ \frac{1}{\sqrt{t}}\mathrm{Ei}\left( -t\right) \right] =-2%
\sqrt{\frac{\pi}{s}}\,\mathrm{\ln }\left( \sqrt{s}+\sqrt{s+1}\right) ,\quad \mathrm{Re}%
\,s>0, \\
\displaystyle%
\mathcal{L}\left[ \mathrm{Si}\left( t\right) \right] =\frac{\cot ^{-1}s}{s}%
,\quad \mathrm{Re}\,s>0, \\
\displaystyle%
\mathcal{L}\left[ \mathrm{si}\left( t\right) \right] =\frac{\tan ^{-1}s}{s},
\\
\displaystyle%
\mathcal{L}\left[ \mathrm{Ci}\left( t\right) \right] =-\frac{\mathrm{\ln }%
\left( 1+s^{2}\right) }{2s}, \\
\displaystyle%
\mathcal{L}^{-1}\left[ \frac{\mathrm{\ln }\left( s+b\right) }{s+a}\right]
=e^{-at}\left[ \mathrm{\ln }\left( b-a\right) -\mathrm{Ei}\left( \left(
a-b\right) t\right) \right] ,\quad \mathrm{Re}\,\left( s-a\right) >0, \\
\displaystyle%
\mathcal{L}^{-1}\left[ \frac{\mathrm{\ln }s}{s^{2}+1}\right] =\cos t\,%
\mathrm{Si}\left( t\right) -\sin t\,\mathrm{Ci}\left( t\right) ,\quad
\mathrm{Re}\,s>0, \\
\displaystyle%
\mathcal{L}^{-1}\left[ \frac{s \, \mathrm{\ln }s}{s^{2}+1}\right] =-\sin t\,%
\mathrm{Si}\left( t\right) -\cos t\,\mathrm{Ci}\left( t\right) .%
\end{array}%
\label{Laplace_transforms_examples}
\end{equation}

Integrands in the second group of integral special functions include special functions,
the most well-known and applied of which are the integral Bessel functions (see, \linebreak e.g., \cite%
{Abramowitz,Oberhettinger,vanderPol,vanderPolBook,Humbert,ApelblatIntegral,ApelblatBessel})%
\begin{equation}
\begin{array}{l}
\displaystyle%
\mathrm{Ji}_{\nu }\left( x\right) =-\int_{x}^{\infty }\frac{J_{\nu }\left(
t\right) }{t}dt, \\
\displaystyle%
\mathrm{Yi}_{\nu }\left( x\right) =-\int_{x}^{\infty }\frac{Y_{\nu }\left(
t\right) }{t}dt, \\
\displaystyle%
\mathrm{Ii}_{\nu }\left( x\right) =-\int_{0}^{x}\frac{I_{\nu }\left(
t\right) }{t}dt, \\
\displaystyle%
\mathrm{Ki}_{\nu }\left( x\right) =-\int_{x}^{\infty }\frac{K_{\nu }\left(
t\right) }{t}dt.%
\end{array}
\label{Integral_Bessel_functions}
\end{equation}

Already in 1929, van der Pol \cite{vanderPol} showed that it is possible to
express the differentiation with respect to the order of the Bessel function
of the first kind as a convolution integral, which includes the integral
Bessel function of the zero-order:%
\begin{equation}
\frac{\partial J_{\nu }\left( t\right) }{\partial \nu }=\frac{1}{2}%
\int_{0}^{t}\mathrm{Ji}_{0}\left( t - x\right) \left[ J_{\nu -1}\left( x\right)
-J_{\nu +1}\left( x\right) \right] dx.  \label{DJ_nu/d_nu}
\end{equation}

The integral Bessel functions of the zero-order are inverse transforms of
the following Laplace transforms \cite{Oberhettinger}%
\begin{equation}
\begin{array}{l}
\displaystyle%
\mathcal{L}^{-1}\left[ \frac{\sinh ^{-1}s}{s}\right] =\mathrm{Ji}_{0}\left(
t\right) , \\
\displaystyle%
\mathcal{L}^{-1}\left[ \frac{\left( \sinh ^{-1}s\right) ^{2}}{s}\right] =%
\mathrm{Yi}_{0}\left( t\right) , \\
\displaystyle%
\mathcal{L}^{-1}\left[ \mathrm{\ln }\left( s+\sqrt{s^{2}+1}\right) -\frac{%
\pi i}{2}\right] =\mathrm{Ii}_{0}\left( t\right) , \\
\displaystyle%
\mathcal{L}^{-1}\left[ \frac{\left( \cosh ^{-1}s\right) ^{2}}{2s}+\frac{\pi
^{2}}{8s}\right] =\mathrm{Ki}_{0}\left( t\right) .%
\end{array}
\label{Inverse_Laplace=Integral_Bessel}
\end{equation}

In analogy to the integral Bessel functions and with the possibility of extension to other special functions,
this work introduces three new integral functions. Furthermore, these integral functions guide us toward the establishment of integrals and series.
Section \ref{Section_Mittag-Leffler} explores the integral Mittag-Leffler functions. Sections \ref{Section_Whittaker} and \ref{Section_Wright} discuss the integral Whittaker and Wright functions, respectively.
Section \ref{Section_Conclusions} contains concluding remarks.

In order to preserve the applied form of notation, the following two
integral functions are introduced:%
\begin{equation}
\mathrm{Fi}\left( x\right) =\int_{0}^{x}\frac{f\left( t\right) -f\left(
0\right) }{t}dt,  \label{Fi_notation}
\end{equation}%
and%
\begin{equation}
\mathrm{fi}\left( x\right) =\int_{x}^{\infty }\frac{f\left( t\right) }{t}dt.
\label{fi_notation}
\end{equation}

To ensure convergence of integrals in (\ref{Fi_notation}) or in (\ref%
{fi_notation}), which depends on the behavior of $f\left( t\right) /t$ integrands
at the origin and at infinity, the forms of integral functions $\mathrm{Fi}%
\left( x\right) $ or $\mathrm{fi}\left( x\right) $ are chosen. Since the
explicit expressions for $f\left( t\right) $ functions are sometimes given
in the form of $f\left( t^{\alpha }\right) $ where $\alpha =\pm \frac{1}{2}%
,\pm 1,2,3,\ldots $ the corresponding change of integration variables for
these equations is desired.

In the case of Mittag-Leffler, Whittaker and Wright functions, for some
values of parameters, by using the MATHEMATICA program, it was possible to
obtain these integral functions in a closed-form. Derived integral functions
are tabulated and also in some cases graphically presented (see \cite%
{Abramowitz}).

\section{The Integral Mittag-Leffler Functions} \label{Section_Mittag-Leffler}

The classical one-parameter and the two-parameter Mittag-Leffler functions
are defined by \cite{GorenfloBook}:%
\begin{equation}
\begin{array}{l}
\displaystyle%
\mathrm{E}_{\alpha }\left( x\right) =\sum_{k=0}^{\infty }\frac{x^{k}}{\Gamma
\left( \alpha k+1\right) },\quad \mathrm{Re}\,\alpha >0, \\
\displaystyle%
\mathrm{E}_{\alpha ,\beta }\left( x\right) =\sum_{k=0}^{\infty }\frac{x^{k}}{%
\Gamma \left( \alpha k+\beta \right) },\quad \mathrm{Re}\,\alpha >0,\
\mathrm{Re}\,\beta >0.%
\end{array}
\label{Mittag_Leffler_def}
\end{equation}

In this investigation, they are only considered for positive real values of
the argument, i.e., $x>0$. In the particular case of positive rational $\alpha $
with $\alpha =p/q$ and $p$ and $q$ positive coprimes, Mittag-Leffler
functions are given as a finite sum of generalized hypergeometric functions
(see (\ref{Mittag-Leffler_p/q_resultado}) in Appendix \ref{Appendix:
Mittag-Leffler}).

The Laplace transforms of the Mittag-Leffler functions are derived directly
from (\ref{Laplace_transform_def})\ and (\ref{Mittag_Leffler_def}), and we
have:%
\begin{equation}
\begin{array}{l}
\displaystyle%
\mathcal{L}\left[ \mathrm{E}_{\alpha }\left( t\right) \right]
=\int_{0}^{\infty }e^{-st}\left[ \sum_{k=0}^{\infty }\frac{t^{k}}{\Gamma
\left( \alpha k+1\right) }\right] dt=\sum_{k=0}^{\infty }\frac{k!}{\Gamma
\left( \alpha k+1\right) }\left( \frac{1}{s}\right) ^{k+1}, \\
\displaystyle%
\mathcal{L}\left[ \mathrm{E}_{\alpha ,\beta }\left( t\right) \right]
=\int_{0}^{\infty }e^{-st}\left[ \sum_{k=0}^{\infty }\frac{t^{k}}{\Gamma
\left( \alpha k+\beta \right) }\right] dt=\sum_{k=0}^{\infty }\frac{k!}{%
\Gamma \left( \alpha k+\beta \right) }\left( \frac{1}{s}\right) ^{k+1}, \\
\mathrm{Re}\,s>1.%
\end{array}
\label{LT_Mittag-Leffler}
\end{equation}

For particular values of parameters $\alpha $ and $\beta $, the explicit
form of the Mittag-Leffler functions can be obtained by applying the MATHEMATICA
program to sums of infinite series in (\ref{Mittag_Leffler_def}), and these
results are presented in Appendix \ref{Appendix: Mittag-Leffler}. Using
Equation (\ref{LT_Mittag-Leffler}), many new Laplace transforms of the
Mittag-Leffler functions were evaluated, and they are also reported in
Appendix \ref{Appendix: Mittag-Leffler}. Similarly as in the case when $%
\alpha $ is positive rational, the Laplace transforms of the Mittag-Leffler
functions can be expressed by the finite sum of products of generalized
hypergeometric functions (see (\ref{Laplace_Mittag-Leffler_p/q_resultado})
in Appendix \ref{Appendix: Mittag-Leffler}).

The integral Mittag-Leffler functions are introduced by considering their
exponential behavior as a function of real, positive variable $x$ (see
Appendix \ref{Appendix: Mittag-Leffler}).%
\begin{equation}
\begin{array}{l}
\displaystyle%
\mathrm{Ei}_{\alpha }\left( x\right) =\int_{0}^{x}\frac{\mathrm{E}_{\alpha
}\left( t\right) -1}{t}dt, \\
\displaystyle%
\mathrm{Ei}_{\alpha ,\beta }\left( x\right) =\int_{0}^{x}\frac{\mathrm{E}%
_{\alpha ,\beta }\left( t\right) -1/\Gamma \left( \beta \right) }{t}dt.%
\end{array}
\label{Integral_Mittag-Leffler_def}
\end{equation}

Formally, by introducing (\ref{Mittag_Leffler_def}) into (\ref%
{Integral_Mittag-Leffler_def}) we have%
\begin{equation}
\begin{array}{l}
\displaystyle%
\mathrm{Ei}_{\alpha }\left( x\right) =\sum_{k=1}^{\infty }\frac{x^{k}}{k\
\Gamma \left( \alpha k+1\right) }, \\
\displaystyle%
\mathrm{Ei}_{\alpha ,\beta }\left( x\right) =\sum_{k=1}^{\infty }\frac{x^{k}%
}{k\ \Gamma \left( \alpha k+\beta \right) }.%
\end{array}
\label{Integral_Mittag-Leffler_series}
\end{equation}

\begin{figure}[H]
\includegraphics[width=10cm]{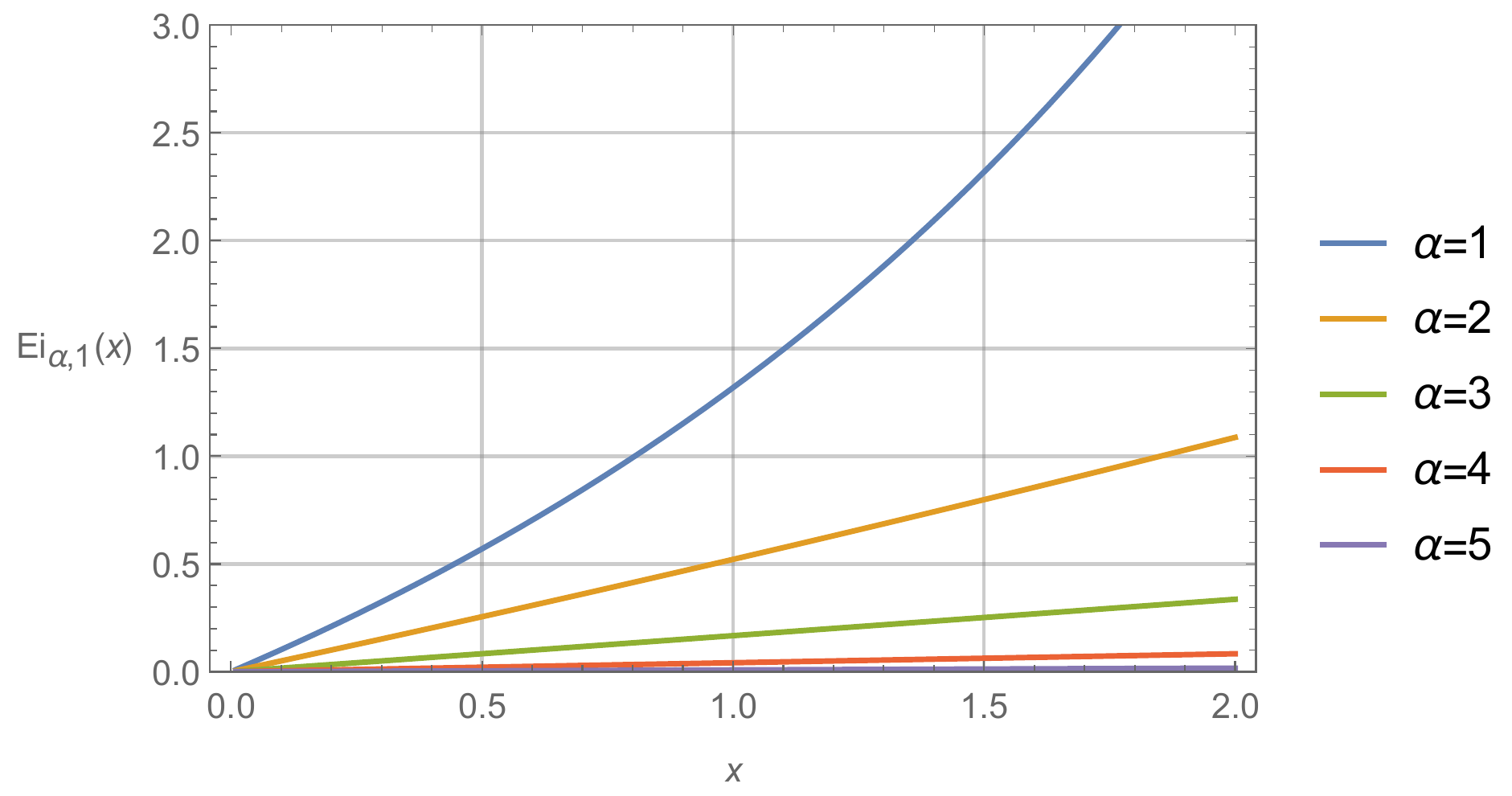}
\caption{The integral one-parameter Mittag-Leffler function $\mathrm{Ei}_{%
\protect\alpha ,1}\left( x\right) $ as a function of variable $x$ and
parameters $\protect\alpha $.}
\label{Figure: Mittag-Leffler a}
\end{figure}

For several values of parameters $\alpha $ and $\beta $, it is possible to
derive the integral Mittag-Leffler functions in a closed-form by applying
the MATHEMATICA program to the sums of infinite series in (\ref%
{Integral_Mittag-Leffler_series}). These functions are presented in Tables %
\ref{Table1a} and \ref{Table1b}. As it is observable, most of these integral functions are expressed as generalized hypergeometric series.
Typical behavior of one-parameter and two-parameter integral Mittag-Leffler
functions is illustrated in Figures \ref{Figure: Mittag-Leffler a} and \ref%
{Figure: Mittag-Leffler b}.

Evidently, also direct integration, by using (\ref%
{Integral_Mittag-Leffler_def}), leads to the integral Mittag-Leffler
functions. For example, for $\mathrm{E}_{2}\left( x\right) =\cosh \sqrt{x}$,
according to (\ref{Trigonometric_integral_functions}), we have%
\begin{equation}
\mathrm{Ei}_{2}\left( x\right) =\int_{0}^{x}\frac{\cosh \sqrt{t}-1}{t}%
dt=-2\gamma -\mathrm{\ln }x+2\,\mathrm{Chi}\sqrt{x},  \label{Ei_2(x)}
\end{equation}%
and as expected, this result is identical to that derived from (\ref%
{Integral_Mittag-Leffler_series}) (see Tables \ref{Table1a} and \ref{Table1b}%
).

Applying the formulas (\ref{Sum_split})\ and (\ref{Multiplication_Gamma})
given in Appendix \ref{Appendix: Mittag-Leffler}, the integral
Mittag-Leffler function for positive rational values of parameter $\alpha $
with $\alpha =p/q$ and $p,q$ positive coprimes is
\begin{eqnarray}
&&\mathrm{Ei}_{p/q,\beta }\left( x\right)
\label{Integral_Mittag-Leffler_p/q} \\
&=&\sum_{k=1}^{q}\frac{x^{k}}{k\,\Gamma \left( k/q+\beta \right) }%
\,_{2}F_{q}\left( \left.
\begin{array}{c}
1,k/q \\
b_{0},\ldots ,b_{p-1},k/q+1%
\end{array}%
\right\vert \frac{x^{q}}{p^{p}}\right) .  \nonumber
\end{eqnarray}%
where%
\[
b_{j}=\frac{k}{q}+\frac{\beta +j}{p}.
\]

\vspace{-12pt}
\begin{specialtable}[H] \centering%
\caption{The integral Mittag-Leffler functions derived for some values of
parameters $\alpha$ and $\beta$ by using
(\ref{Integral_Mittag-Leffler_series}).}%
\begin{tabular*}{\hsize}{c@{\extracolsep{\fill}}cc}
\toprule
\boldmath$\alpha $ & \boldmath$\beta $ & \boldmath$\mathrm{Ei}_{\alpha ,\beta }\left( x\right) $ \\
\midrule
$\frac{1}{3}$ & $\frac{1}{5}$ & $\frac{x}{\Gamma \left( \frac{8}{15}\right) }%
\,_{2}F_{2}\left( \left.
\begin{array}{c}
1,\frac{1}{3} \\
\frac{8}{15},\frac{4}{3}%
\end{array}%
\right\vert x^{3}\right) +\frac{x^{2}}{2\,\Gamma \left( \frac{13}{15}\right)
}\,_{2}F_{2}\left( \left.
\begin{array}{c}
1,\frac{2}{3} \\
\frac{13}{15},\frac{5}{3}%
\end{array}%
\right\vert x^{3}\right) +\frac{5\,x^{3}}{3\,\Gamma \left( \frac{1}{5}%
\right) }\,_{2}F_{2}\left( \left.
\begin{array}{c}
1,1 \\
\frac{6}{5},2%
\end{array}%
\right\vert x^{3}\right) $ \\ \midrule
$\frac{1}{3}$ & $\frac{1}{4}$ & $\frac{x}{\Gamma \left( \frac{7}{12}\right) }%
\,_{2}F_{2}\left( \left.
\begin{array}{c}
1,\frac{1}{3} \\
\frac{7}{12},\frac{4}{3}%
\end{array}%
\right\vert x^{3}\right) +\frac{x^{2}}{2\,\Gamma \left( \frac{11}{12}\right)
}\,_{2}F_{2}\left( \left.
\begin{array}{c}
1,\frac{2}{3} \\
\frac{11}{12},\frac{5}{3}%
\end{array}%
\right\vert x^{3}\right) +\frac{4x^{3}}{3\,\Gamma \left( \frac{1}{4}\right) }%
\,_{2}F_{2}\left( \left.
\begin{array}{c}
1,1 \\
\frac{5}{4},2%
\end{array}%
\right\vert x^{3}\right) $ \\ \midrule
$\frac{1}{3}$ & $\frac{1}{2}$ & $\frac{2\,x^{3}}{3\sqrt{\pi }}%
\,_{2}F_{2}\left( \left.
\begin{array}{c}
1,1 \\
\frac{3}{2},2%
\end{array}%
\right\vert x^{3}\right) +\frac{x}{\Gamma \left( \frac{5}{6}\right) }%
\,_{2}F_{2}\left( \left.
\begin{array}{c}
1,\frac{1}{3} \\
\frac{5}{6},\frac{4}{3}%
\end{array}%
\right\vert x^{3}\right) +\frac{3\,x^{2}}{\Gamma \left( \frac{1}{6}\right) }%
\,_{2}F_{2}\left( \left.
\begin{array}{c}
1,\frac{2}{3} \\
\frac{7}{6},\frac{5}{3}%
\end{array}%
\right\vert x^{3}\right) $ \\ \midrule

$\frac{1}{3}$ & $\frac{3}{2}$ & $\frac{x}{\Gamma \left( \frac{11}{6}\right) }%
\,_{2}F_{2}\left( \left.
\begin{array}{c}
1,\frac{1}{3} \\
\frac{4}{3},\frac{11}{6}%
\end{array}%
\right\vert x^{3}\right) +\frac{18\,x^{2}}{7\,\Gamma \left( \frac{1}{6}%
\right) }\,_{2}F_{2}\left( \left.
\begin{array}{c}
1,\frac{2}{3} \\
\frac{5}{3},\frac{13}{6}%
\end{array}%
\right\vert x^{3}\right) +\frac{4\,x^{3}}{9\sqrt{\pi }}\,_{2}F_{2}\left(
\left.
\begin{array}{c}
1,1 \\
2,\frac{5}{2}%
\end{array}%
\right\vert x^{3}\right) $ \\ \midrule

$\frac{1}{2}$ & $\frac{1}{2}$ & $\frac{x^{2}}{\sqrt{\pi }}\,_{2}F_{2}\left(
\left.
\begin{array}{c}
1,1 \\
\frac{3}{2},2%
\end{array}%
\right\vert x^{2}\right) +e^{x^{2}}\,F\left( x\right) ,\quad F\left(
x\right) =e^{-x^{2}}\int_{0}^{x}e^{t^{2}}dt$ \\ \midrule
$\frac{1}{2}$ & $1$ & $-\frac{\gamma }{2}-\mathrm{\ln }x+\frac{\mathrm{Ei}%
\left( x^{2}\right) }{2}+\frac{2\,x}{\sqrt{\pi }}\,_{2}F_{2}\left( \left.
\begin{array}{c}
\frac{1}{2},1 \\
\frac{3}{2},\frac{3}{2}%
\end{array}%
\right\vert x^{2}\right) $ \\ \midrule
$\frac{1}{2}$ & $2$ & $\frac{1}{2}\left( 1-\gamma +\mathrm{Ei}\left(
x^{2}\right) +\frac{1-e^{x^{2}}}{x^{2}}\right) -\mathrm{\ln }x+\frac{4\,x}{3%
\sqrt{\pi }}\,_{2}F_{2}\left( \left.
\begin{array}{c}
\frac{1}{2},1 \\
\frac{3}{2},\frac{3}{2}%
\end{array}%
\right\vert x^{2}\right) $ \\

\bottomrule
\end{tabular*}
\end{specialtable}

\begin{specialtable}[H]\ContinuedFloat
\small
\caption{{\em Cont.}}
\begin{tabular*}{\hsize}{c@{\extracolsep{\fill}}cc}
\toprule
\boldmath$\alpha $ & \boldmath$\beta $ & \boldmath$\mathrm{Ei}_{\alpha ,\beta }\left( x\right) $ \\
\midrule

$\frac{1}{2}$ & $3$ & $\frac{2+4x^{2}+\left( 3-2\gamma \right)
x^{4}-2e^{x^{2}}\left( 1+x^{2}\right) +2x^{4}\left( \mathrm{Ei}\left(
x^{2}\right) -2\mathrm{\ln }x\right) }{8x^{4}}+\frac{8\,x}{15\sqrt{\pi }}%
\,_{2}F_{2}\left( \left.
\begin{array}{c}
\frac{1}{2},1 \\
\frac{3}{2},\frac{7}{2}%
\end{array}%
\right\vert x^{2}\right) $ \\ \midrule

$\frac{1}{2}$ & $4$ & $\frac{12+18x^{2}\left( 1+x^{2}\right) +\left(
11-6\gamma \right) x^{6}-6e^{x^{2}}\left( 2+x^{2}+x^{4}\right) +3x^{6}\left(
2\mathrm{Ei}\left( x^{2}\right) -4\mathrm{\ln }x\right) }{72x^{6}}+\frac{%
16\,x}{105\sqrt{\pi }}\,_{2}F_{2}\left( \left.
\begin{array}{c}
\frac{1}{2},1 \\
\frac{3}{2},\frac{9}{2}%
\end{array}%
\right\vert x^{2}\right) $ \\ \midrule

$\frac{1}{2}$ & $\beta $ & $\frac{x^{2}}{2\,\Gamma \left( \beta +1\right) }%
\,_{2}F_{2}\left( \left.
\begin{array}{c}
1,1 \\
2,\beta +1%
\end{array}%
\right\vert x^{2}\right) +\frac{\,x}{\Gamma \left( \beta +\frac{1}{2}\right)
}\,_{2}F_{2}\left( \left.
\begin{array}{c}
\frac{1}{2},1 \\
\frac{3}{2},\beta +\frac{1}{2}%
\end{array}%
\right\vert x^{2}\right) $ \\ \midrule
$1$ & $\frac{1}{4}$ & $\frac{x}{\Gamma \left( \frac{5}{4}\right) }%
\,_{2}F_{2}\left( \left.
\begin{array}{c}
1,1 \\
\frac{5}{4},2%
\end{array}%
\right\vert x\right) $ \\ \midrule
$1$ & $\frac{1}{3}$ & $\frac{x}{\Gamma \left( \frac{4}{3}\right) }%
\,_{2}F_{2}\left( \left.
\begin{array}{c}
1,1 \\
\frac{4}{3},2%
\end{array}%
\right\vert x\right) $ \\ \midrule
$1$ & $\frac{1}{2}$ & $\frac{2x}{\sqrt{\pi }}\,_{2}F_{2}\left( \left.
\begin{array}{c}
1,1 \\
\frac{3}{2},2%
\end{array}%
\right\vert x\right) $ \\ \midrule
$1$ & $1$ & $-\gamma -\mathrm{\ln }x+\mathrm{Chi}\left( x\right) +\mathrm{Shi%
}\left( x\right) $ \\ \midrule
$1$ & $\frac{3}{2}$ & $\frac{4x}{3\sqrt{\pi }}\,_{2}F_{2}\left( \left.
\begin{array}{c}
1,1 \\
\frac{5}{2},2%
\end{array}%
\right\vert x\right) $ \\ \midrule
$1$ & $\beta $ & $\frac{x}{\Gamma \left( \beta +1\right) }\,_{2}F_{2}\left(
\left.
\begin{array}{c}
1,1 \\
2,1+\beta%
\end{array}%
\right\vert x\right) $ \\ \bottomrule
\end{tabular*}%
\label{Table1a}%
\end{specialtable}%
\vspace{-6pt}
\vspace{-6pt}
\begin{figure}[H]
\includegraphics[width=10cm]{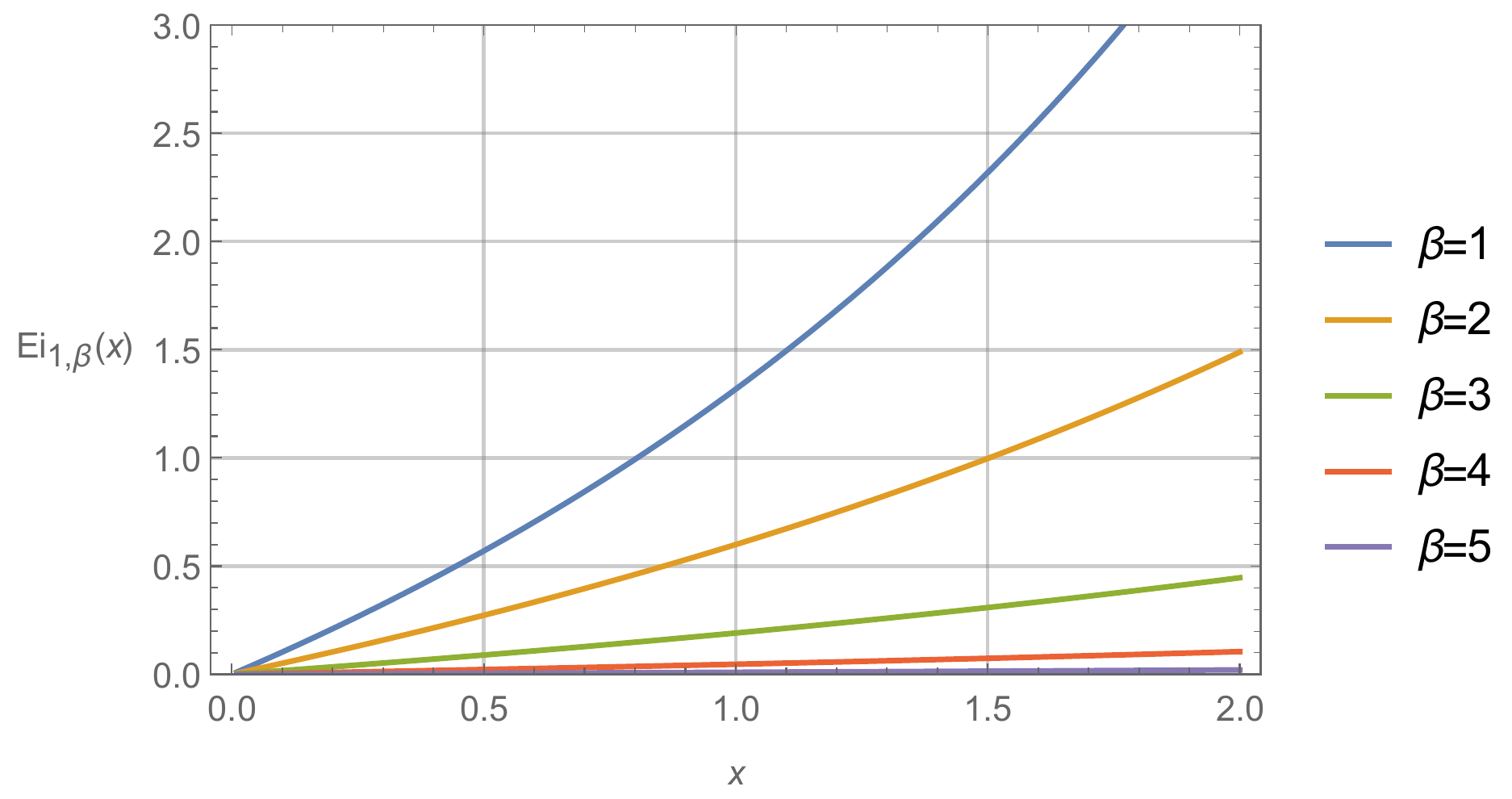}
\caption{The integral two-parameter Mittag-Leffler function $\mathrm{Ei}_{1,%
\protect\beta }\left( x\right) $ as a function of variable $x$ and
parameters $\protect\beta $.}
\label{Figure: Mittag-Leffler b}
\end{figure}

\begin{specialtable}[H] \centering%
\caption{The integral Mittag-Leffler functions derived for some values of
parameters $\alpha$ and $\beta$ by using
(\ref{Integral_Mittag-Leffler_series}).}%
\begin{tabular*}{\hsize}{c@{\extracolsep{\fill}}cc}
\toprule
\boldmath$\alpha $ & \boldmath$\beta $ & \boldmath$\mathrm{Ei}_{\alpha ,\beta }\left( x\right) $ \\
\midrule
$\frac{3}{2}$ & $\frac{1}{2}$ & $\frac{4x^{2}}{15\sqrt{\pi }}%
\,_{2}F_{4}\left( \left.
\begin{array}{c}
1,1 \\
\frac{7}{6},\frac{3}{2},\frac{11}{6},2%
\end{array}%
\right\vert \frac{x^{2}}{27}\right) +\frac{x}{48}\,_{1}F_{3}\left( \left.
\begin{array}{c}
\frac{1}{2} \\
\frac{2}{3},\frac{4}{3},\frac{3}{2}%
\end{array}%
\right\vert \frac{x^{2}}{27}\right) $ \\ \midrule
$\frac{3}{2}$ & $1$ & $\frac{4x}{3\sqrt{\pi }}\,_{2}F_{4}\left( \left.
\begin{array}{c}
\frac{1}{2},1 \\
\frac{5}{6},\frac{7}{6},\frac{3}{2},\frac{3}{2}%
\end{array}%
\right\vert \frac{x^{2}}{27}\right) +\frac{x^{2}}{12}\,_{2}F_{4}\left(
\left.
\begin{array}{c}
1,1 \\
\frac{4}{3},\frac{5}{3},2,2%
\end{array}%
\right\vert \frac{x^{2}}{27}\right) $ \\ \midrule
$\frac{3}{2}$ & $\frac{3}{2}$ & $\frac{x}{2}\,_{1}F_{3}\left( \left.
\begin{array}{c}
\frac{1}{2} \\
\frac{4}{3},\frac{3}{2},\frac{5}{3}%
\end{array}%
\right\vert \frac{x^{2}}{27}\right) +\frac{8x^{2}}{105\sqrt{\pi }}%
\,_{2}F_{4}\left( \left.
\begin{array}{c}
1,1 \\
\frac{3}{2},\frac{11}{6},2,\frac{13}{6}%
\end{array}%
\right\vert \frac{x^{2}}{27}\right) $ \\ \midrule
$\frac{3}{2}$ & $2$ & $\frac{8x}{15\sqrt{\pi }}\,_{2}F_{4}\left( \left.
\begin{array}{c}
\frac{1}{2},1 \\
\frac{7}{6},\frac{3}{2},\frac{3}{2},\frac{11}{6}%
\end{array}%
\right\vert \frac{x^{2}}{27}\right) +\frac{x^{2}}{48}\,_{2}F_{4}\left(
\left.
\begin{array}{c}
1,1 \\
\frac{5}{3},2,2,\frac{7}{3}%
\end{array}%
\right\vert \frac{x^{2}}{27}\right) $ \\ \midrule
$2$ & $\frac{1}{4}$ & $\frac{16x}{5\,\Gamma \left( \frac{1}{4}\right) }%
\,_{2}F_{3}\left( \left.
\begin{array}{c}
1,1 \\
\frac{9}{8},\frac{11}{8},2%
\end{array}%
\right\vert \frac{x}{4}\right) $ \\ \midrule
$2$ & $\frac{1}{3}$ & $\frac{9x}{4\,\Gamma \left( \frac{1}{3}\right) }%
\,_{2}F_{3}\left( \left.
\begin{array}{c}
1,1 \\
\frac{7}{6},\frac{5}{3},2%
\end{array}%
\right\vert \frac{x}{4}\right) $ \\ \midrule
$2$ & $\frac{1}{2}$ & $\frac{9x}{3\sqrt{\pi }}\,_{2}F_{3}\left( \left.
\begin{array}{c}
1,1 \\
\frac{5}{4},\frac{7}{4},2%
\end{array}%
\right\vert \frac{x}{4}\right) $ \\ \midrule
$2$ & $1$ & $-2\gamma -\mathrm{\ln }x+2\mathrm{Chi}\left( \sqrt{x}\right) $
\\ \midrule
$2$ & $2$ & $2-2\gamma -\mathrm{\ln }x-\frac{2\sinh \sqrt{x}}{\sqrt{x}}+2%
\mathrm{Chi}\left( \sqrt{x}\right) $ \\ \midrule
$2$ & $3$ & $\frac{x}{24}\,_{2}F_{3}\left( \left.
\begin{array}{c}
1,1 \\
2,\frac{5}{2},3%
\end{array}%
\right\vert \frac{x}{4}\right) $ \\ \midrule
$2$ & $4$ & $\frac{x}{120}\,_{2}F_{3}\left( \left.
\begin{array}{c}
1,1 \\
2,\frac{7}{2},3%
\end{array}%
\right\vert \frac{x}{4}\right) $ \\ \midrule
$2$ & $\beta $ & $\frac{x}{\Gamma \left( \beta +2\right) }\,_{2}F_{3}\left(
\left.
\begin{array}{c}
1,1 \\
2,\frac{\beta }{2}+1,\frac{\beta +3}{2}%
\end{array}%
\right\vert \frac{x}{4}\right) $ \\ \midrule
$3$ & $1$ & $\frac{x}{6}\,_{2}F_{4}\left( \left.
\begin{array}{c}
1,1 \\
\frac{4}{3},\frac{5}{3},2,2%
\end{array}%
\right\vert \frac{x}{27}\right) $ \\ \midrule
$3$ & $\beta $ & $\frac{x}{\Gamma \left( \beta +3\right) }\,_{2}F_{4}\left(
\left.
\begin{array}{c}
1,1 \\
2,\frac{\beta }{3}+1,\frac{\beta +4}{3},\frac{\beta +5}{3}%
\end{array}%
\right\vert \frac{x}{27}\right) $ \\ \midrule
$4$ & $1$ & $\frac{x}{24}\,_{2}F_{5}\left( \left.
\begin{array}{c}
1,1 \\
\frac{5}{4},\frac{3}{2},\frac{7}{4},2,2%
\end{array}%
\right\vert \frac{x}{256}\right) $ \\ \midrule
$4$ & $\beta $ & $\frac{x}{\Gamma \left( \beta +4\right) }\,_{2}F_{5}\left(
\left.
\begin{array}{c}
1,1 \\
2,\frac{\beta }{4}+1,\frac{\beta +5}{4},\frac{3}{2}+\frac{\beta }{4},\frac{%
\beta +7}{4}%
\end{array}%
\right\vert \frac{x}{256}\right) $ \\ \midrule
$5$ & $1$ & $\frac{x}{120}\,_{2}F_{6}\left( \left.
\begin{array}{c}
1,1 \\
\frac{6}{5},\frac{7}{5},\frac{8}{5},\frac{9}{5},2,2%
\end{array}%
\right\vert \frac{x}{3125}\right) $ \\ \midrule
$5$ & $\beta $ & $\frac{x}{\Gamma \left( \beta +5\right) }\,_{2}F_{6}\left(
\left.
\begin{array}{c}
1,1 \\
2,\frac{\beta }{5}+1,\frac{\beta +6}{5},\frac{\beta +7}{5},\frac{\beta +8}{5}%
,\frac{\beta +9}{5}%
\end{array}%
\right\vert \frac{x}{3125}\right) $ \\ \bottomrule
\end{tabular*}%
\label{Table1b}%
\end{specialtable}%

In addition, using the sums in (\ref{Integral_Mittag-Leffler_series}), it is
possible to derive the Laplace transforms of the integral Mittag-Leffler
functions:%
\begin{equation}
\begin{array}{l}
\displaystyle%
\mathcal{L}\left[ \mathrm{Ei}_{\alpha }\left( t\right) \right]
=\int_{0}^{\infty }e^{-st}\left[ \sum_{k=0}^{\infty }\frac{t^{k+1}}{\left(
k+1\right) \Gamma \left( \alpha \left( k+1\right) +1\right) }\right] dt \\
\displaystyle%
=\sum_{k=0}^{\infty }\frac{\left( k+1\right) !}{\left( k+1\right) \Gamma
\left( \alpha \left( k+1\right) +1\right) }\left( \frac{1}{s}\right)
^{k+2},\quad \mathrm{Re}\,s>1. \\
\displaystyle%
\mathcal{L}\left[ \mathrm{Ei}_{\alpha ,\beta }\left( t\right) \right]
=\int_{0}^{\infty }e^{-st}\left[ \sum_{k=0}^{\infty }\frac{t^{k+1}}{\left(
k+1\right) \Gamma \left( \alpha \left( k+\beta \right) +1\right) }\right] dt
\\
\displaystyle%
=\sum_{k=0}^{\infty }\frac{\left( k+1\right) !}{\left( k+1\right) \Gamma
\left( \alpha \left( k+1\right) +\beta \right) }\left( \frac{1}{s}\right)
^{k+2},\quad \mathrm{Re}\,s>1.%
\end{array}
\label{LT_Integral_Mittag-Leffler_series}
\end{equation}

The evaluated Laplace transforms of the integral Mittag-Leffler functions
are presented in Tables \ref{Table2} and \ref{Table2a}.

The Laplace transforms of the integral Mittag-Leffler functions with
positive rational parameter $\alpha $ with $\alpha =p/q$ and $p,q$ positive
coprimes can be evaluated from:%
\begin{eqnarray}
&&\mathcal{L}\left[ \mathrm{Ei}_{p/q,\beta }\left( t\right) \right]
\label{Laplace_Integral_Mittag-Leffler_p/q} \\
&=&\frac{1}{s^{2}}\sum_{k=0}^{q-1}\frac{k!\,s^{-k}}{\Gamma \left( \frac{p}{q}%
\left( k+1\right) +\beta \right) }\,_{q+1}F_{p}\left( \left.
\begin{array}{c}
1,a_{0},\ldots ,a_{q-1} \\
b_{0},\ldots ,b_{p-1}%
\end{array}%
\right\vert \frac{\left( q/s\right) ^{q}}{p^{p}}\right) .  \nonumber
\end{eqnarray}%
where%
\begin{eqnarray*}
a_{j} &=&\frac{k+1+j}{q}, \\
b_{j} &=&\frac{k}{q}+\frac{\beta +j}{p}.
\end{eqnarray*}

Furthermore, the following relation is satisfied:%
\begin{equation}
\mathcal{L}\left[ \mathrm{Ei}_{p/q,\beta }\left( t\right) \right] =\frac{1}{%
p^{p/q}s}\mathcal{L}\left[ \mathrm{E}_{p/q,\beta }\left( t\right) \right] .
\label{LT_Mittag-Leffler_relation_p/q}
\end{equation}

\vspace{-6pt}

\begin{specialtable}[H] \centering%
\caption{The Laplace transforms of the integral Mittag-Leffler functions $\mathrm{Ei}_{\alpha , \beta}$  derived for some values of
parameters $\alpha$ and $\beta$ by using
(\ref{LT_Integral_Mittag-Leffler_series}).}%
\begin{tabular*}{\hsize}{c@{\extracolsep{\fill}}cc}
\toprule
\boldmath$\alpha $ & \boldmath$\beta $ & \boldmath$\mathcal{L}\left[ \mathrm{Ei}_{\alpha ,\beta }\left(
t\right) \right] $ \\ \midrule
$1$ & $\frac{1}{5}$ & $\frac{5}{\,\Gamma \left( \frac{1}{5}\right) s^{2}}%
\,_{2}F_{1}\left( \left.
\begin{array}{c}
1,1 \\
\frac{6}{5}%
\end{array}%
\right\vert \frac{1}{s}\right) $ \\ \midrule
$1$ & $\frac{1}{4}$ & $\frac{4}{\Gamma \left( \frac{1}{4}\right) s^{2}}%
\,_{2}F_{1}\left( \left.
\begin{array}{c}
1,1 \\
\frac{5}{4}%
\end{array}%
\right\vert \frac{1}{s}\right) $ \\ \midrule
$1$ & $\frac{1}{3}$ & $\frac{3}{\Gamma \left( \frac{1}{6}\right) s^{2}}%
\,_{2}F_{1}\left( \left.
\begin{array}{c}
1,1 \\
\frac{4}{3}%
\end{array}%
\right\vert \frac{1}{s}\right) $ \\ \midrule
$1$ & $\frac{1}{2}$ & $\frac{2\csc ^{-1}\left( \sqrt{s}\right) }{\sqrt{\pi }s%
\sqrt{s-1}}$ \\ \midrule
$1$ & $1$ & $-\frac{1}{s}\mathrm{\ln }\left( 1-\frac{1}{s}\right) $ \\ \midrule
$1$ & $\frac{3}{2}$ & $\frac{4}{\sqrt{\pi }s}\left[ 1-\sqrt{s-1}\csc
^{-1}\left( \sqrt{s}\right) \right] $ \\ \midrule
$1$ & $\beta $ & $\frac{1}{s^{2}\Gamma \left( \beta +1\right) }%
\,_{2}F_{1}\left( \left.
\begin{array}{c}
1,1 \\
\beta +1%
\end{array}%
\right\vert \frac{1}{s}\right) $ \\

\bottomrule
\end{tabular*}
\end{specialtable}

\begin{specialtable}[H]\ContinuedFloat
\small
\caption{{\em Cont.}}
\begin{tabular*}{\hsize}{c@{\extracolsep{\fill}}cc}
\toprule
\boldmath$\alpha $ & \boldmath$\beta $ & \boldmath$\mathcal{L}\left[ \mathrm{Ei}_{\alpha ,\beta }\left(
t\right) \right] $ \\ \midrule

$\frac{3}{2}$ & $\frac{1}{2}$ & $\frac{8}{15\sqrt{\pi }s^{3}}%
\,_{2}F_{2}\left( \left.
\begin{array}{c}
1,1 \\
\frac{7}{6},\frac{11}{6}%
\end{array}%
\right\vert \frac{4}{27s^{2}}\right) +\frac{1}{s^{2}}\,_{2}F_{2}\left(
\left.
\begin{array}{c}
\frac{1}{2},1 \\
\frac{2}{3},\frac{4}{3}%
\end{array}%
\right\vert \frac{4}{27s^{2}}\right) $ \\ \midrule
$\frac{3}{2}$ & $1$ & $\frac{4}{3\sqrt{\pi }s^{2}}\,_{3}F_{3}\left( \left.
\begin{array}{c}
\frac{1}{2},1,1 \\
\frac{5}{6},\frac{7}{6},\frac{3}{2}%
\end{array}%
\right\vert \frac{4}{27s^{2}}\right) +\frac{1}{6s^{3}}\,_{3}F_{3}\left(
\left.
\begin{array}{c}
1,1,\frac{3}{2} \\
\frac{4}{3},\frac{5}{6},2%
\end{array}%
\right\vert \frac{4}{27s^{2}}\right) $ \\ \midrule
$\frac{3}{2}$ & $\frac{3}{2}$ & $\frac{1}{2s^{2}}\,_{2}F_{2}\left( \left.
\begin{array}{c}
\frac{1}{2},1 \\
\frac{4}{3},\frac{5}{3}%
\end{array}%
\right\vert \frac{4}{27s^{2}}\right) +\frac{16}{105\sqrt{\pi }s^{3}}%
\,_{2}F_{2}\left( \left.
\begin{array}{c}
1,1 \\
\frac{11}{6},\frac{13}{6}%
\end{array}%
\right\vert \frac{4}{27s^{2}}\right) $ \\ \midrule
$\frac{3}{2}$ & $2$ & $\frac{8}{15\sqrt{\pi }s^{2}}\,_{3}F_{3}\left( \left.
\begin{array}{c}
\frac{1}{2},1,1 \\
\frac{7}{6},\frac{3}{2},\frac{11}{6}%
\end{array}%
\right\vert \frac{4}{27s^{2}}\right) +\frac{1}{24s^{3}}\,_{3}F_{3}\left(
\left.
\begin{array}{c}
1,1,\frac{3}{2} \\
\frac{5}{3},2,\frac{7}{3}%
\end{array}%
\right\vert \frac{4}{27s^{2}}\right) $ \\ \midrule
$2$ & $\frac{1}{5}$ & $\frac{25}{6\,\Gamma \left( \frac{1}{5}\right) s^{2}}%
\,_{2}F_{2}\left( \left.
\begin{array}{c}
1,1 \\
\frac{11}{10},\frac{8}{5}%
\end{array}%
\right\vert \frac{1}{4s}\right) $ \\ \midrule
$2$ & $\frac{1}{4}$ & $\frac{16}{5\,\Gamma \left( \frac{1}{4}\right) s^{2}}%
\,_{2}F_{2}\left( \left.
\begin{array}{c}
1,1 \\
\frac{9}{8},\frac{13}{8}%
\end{array}%
\right\vert \frac{1}{4s}\right) $ \\ \midrule
$2$ & $\frac{1}{3}$ & $\frac{9}{4\,\Gamma \left( \frac{1}{3}\right) s^{2}}%
\,_{2}F_{2}\left( \left.
\begin{array}{c}
1,1 \\
\frac{7}{6},\frac{5}{3}%
\end{array}%
\right\vert \frac{1}{4s}\right) $ \\ \midrule
$2$ & $\frac{1}{2}$ & $\frac{4}{3\sqrt{\pi }s^{2}}\,_{2}F_{2}\left( \left.
\begin{array}{c}
1,1 \\
\frac{5}{4},\frac{7}{4}%
\end{array}%
\right\vert \frac{1}{4s}\right) $ \\ \midrule
$2$ & $1$ & $\frac{1}{2s^{2}}\,_{2}F_{2}\left( \left.
\begin{array}{c}
1,1 \\
\frac{3}{2},2%
\end{array}%
\right\vert \frac{1}{4s}\right) $ \\ \midrule
$2$ & $2$ & $\frac{1}{6s^{2}}\,_{2}F_{2}\left( \left.
\begin{array}{c}
1,1 \\
\frac{5}{2},2%
\end{array}%
\right\vert \frac{1}{4s}\right) $ \\ \midrule
$2$ & $3$ & $\frac{1}{24s^{2}}\,_{2}F_{2}\left( \left.
\begin{array}{c}
1,1 \\
\frac{5}{2},3%
\end{array}%
\right\vert \frac{1}{4s}\right) $ \\ \midrule
$2$ & $4$ & $\frac{1}{120s^{2}}\,_{2}F_{2}\left( \left.
\begin{array}{c}
1,1 \\
\frac{7}{2},3%
\end{array}%
\right\vert \frac{1}{4s}\right) $ \\ \midrule
$2$ & $\beta $ & $\frac{1}{\Gamma \left( \beta +2\right) s^{2}}%
\,_{2}F_{2}\left( \left.
\begin{array}{c}
1,1 \\
1+\frac{\beta }{2},\frac{\beta +3}{2}%
\end{array}%
\right\vert \frac{1}{4s}\right) $ \\ \bottomrule
\end{tabular*}%
\label{Table2}%
\end{specialtable}%

\vspace{-12pt}

\begin{specialtable}[H] \centering%
\caption{The Laplace transforms of the integral Mittag-Leffler functions $\mathrm{Ei}_{\alpha , \beta}$  derived for some values of
parameters $\alpha$ and $\beta$ by using
(\ref{LT_Integral_Mittag-Leffler_series}).}%
\begin{tabular*}{\hsize}{c@{\extracolsep{\fill}}cc}
\toprule
\boldmath$\alpha $ & \boldmath$\beta $ & \boldmath$\mathcal{L}\left[ \mathrm{Ei}_{\alpha ,\beta }\left(
t\right) \right] $ \\ \midrule
$3$ & $\frac{1}{5}$ & $\frac{125}{66\,\Gamma \left( \frac{1}{5}\right) s^{2}}%
\,_{2}F_{3}\left( \left.
\begin{array}{c}
1,1 \\
\frac{16}{15},\frac{7}{5},\frac{26}{15}%
\end{array}%
\right\vert \frac{1}{27s}\right) $ \\ \midrule
$3$ & $\frac{1}{4}$ & $\frac{64}{25\,\Gamma \left( \frac{1}{4}\right) s^{2}}%
\,_{2}F_{3}\left( \left.
\begin{array}{c}
1,1 \\
\frac{13}{12},\frac{17}{12},\frac{7}{4}%
\end{array}%
\right\vert \frac{1}{27s}\right) $ \\ \midrule
$3$ & $\frac{1}{3}$ & $\frac{27}{28\,\Gamma \left( \frac{1}{3}\right) s^{2}}%
\,_{2}F_{3}\left( \left.
\begin{array}{c}
1,1 \\
\frac{10}{9},\frac{13}{9},\frac{10}{9}%
\end{array}%
\right\vert \frac{1}{27s}\right) $ \\ \midrule
$3$ & $\frac{1}{2}$ & $\frac{27}{15\sqrt{\pi }s^{2}}\,_{2}F_{3}\left( \left.
\begin{array}{c}
1,1 \\
\frac{7}{6},\frac{3}{2},\frac{11}{6}%
\end{array}%
\right\vert \frac{1}{27s}\right) $ \\  \midrule

$3$ & $1$ & $\frac{1}{6\,s^{2}}\,_{2}F_{3}\left( \left.
\begin{array}{c}
1,1 \\
\frac{4}{3},\frac{5}{3},2%
\end{array}%
\right\vert \frac{1}{27s}\right) $ \\ \midrule
$3$ & $3$ & $\frac{1}{120\,s^{2}}\,_{2}F_{3}\left( \left.
\begin{array}{c}
1,1 \\
2,\frac{7}{3},\frac{8}{3}%
\end{array}%
\right\vert \frac{1}{27s}\right) $ \\

\bottomrule
\end{tabular*}
\end{specialtable}

\begin{specialtable}[H]\ContinuedFloat
\small
\caption{{\em Cont.}}
\begin{tabular*}{\hsize}{c@{\extracolsep{\fill}}cc}
\toprule
\boldmath$\alpha $ & \boldmath$\beta $ & \boldmath$\mathcal{L}\left[ \mathrm{Ei}_{\alpha ,\beta }\left(
t\right) \right] $ \\ \midrule

$3$ & $\beta $ & $\frac{1}{\Gamma \left( \beta +3\right) s^{2}}%
\,_{2}F_{3}\left( \left.
\begin{array}{c}
1,1 \\
1+\frac{\beta }{3},\frac{\beta +4}{3},\frac{\beta +5}{3}%
\end{array}%
\right\vert \frac{1}{27s}\right) $ \\ \midrule
$4$ & $1$ & $\frac{1}{24\,s^{2}}\,_{2}F_{4}\left( \left.
\begin{array}{c}
1,1 \\
\frac{5}{4},\frac{3}{2},\frac{7}{4},2%
\end{array}%
\right\vert \frac{1}{256\,s}\right) $ \\ \midrule
$4$ & $4$ & $\frac{1}{5040\,s^{2}}\,_{2}F_{4}\left( \left.
\begin{array}{c}
1,1 \\
2,\frac{9}{4},\frac{5}{2},\frac{11}{4}%
\end{array}%
\right\vert \frac{1}{256\,s}\right) $ \\ \midrule
$4$ & $\beta $ & $\frac{1}{\Gamma \left( \beta +4\right) s^{2}}%
\,_{2}F_{4}\left( \left.
\begin{array}{c}
1,1 \\
1+\frac{\beta }{4},\frac{\beta +5}{3},\frac{3}{2}+\frac{\beta }{4},\frac{%
\beta +7}{4}%
\end{array}%
\right\vert \frac{1}{256\,s}\right) $ \\ \midrule
$5$ & $1$ & $\frac{1}{120\,s^{2}}\,_{2}F_{5}\left( \left.
\begin{array}{c}
1,1 \\
\frac{6}{5},\frac{7}{5},\frac{8}{5},\frac{9}{5}%
\end{array}%
\right\vert \frac{1}{3125\,s}\right) $ \\ \midrule
$5$ & $\beta $ & $\frac{1}{\Gamma \left( \beta +5\right) s^{2}}%
\,_{2}F_{5}\left( \left.
\begin{array}{c}
1,1 \\
1+\frac{\beta }{5},\frac{\beta +6}{5},\frac{\beta +7}{5},\frac{\beta +8}{5},%
\frac{\beta +9}{5}%
\end{array}%
\right\vert \frac{1}{3125\,s}\right) $ \\ \bottomrule
\end{tabular*}%
\label{Table2a}%
\end{specialtable}%


\section{The Integral Whittaker Functions} \label{Section_Whittaker}

In 1903, Whittaker \cite{Whittaker} showed that it is possible to express some special functions such as Bessel
functions, parabolic cylinder functions, error functions, incomplete gamma
functions, and logarithm and cosine integrals in terms of a new function suggested by him, i.e., the Whittaker function. Two Whittaker functions are
applied today, and they are defined by using the Kummer confluent
hypergeometric function \cite{Abramowitz,Magnus}:%
\begin{equation}
\begin{array}{l}
\displaystyle%
\mathrm{M}_{\kappa ,\mu }\left( x\right) =x^{\mu
-1/2}e^{-x/2}\,_{1}F_{1}\left( \left.
\begin{array}{c}
\mu -\kappa +\frac{1}{2} \\
1+2\mu%
\end{array}%
\right\vert x\right) , \\
\displaystyle%
\mathrm{W}_{\kappa ,\mu }\left( x\right) =\frac{\Gamma \left( -2\mu \right)
}{\Gamma \left( \frac{1}{2}-k-\mu \right) }\mathrm{M}_{\kappa ,\mu }\left(
x\right) +\frac{\Gamma \left( 2\mu \right) }{\Gamma \left( \frac{1}{2}-k+\mu
\right) }\mathrm{M}_{\kappa ,-\mu }\left( x\right) .%
\end{array}
\label{Whittaker_def}
\end{equation}

This permits us to introduce four integral Whittaker functions:%
\begin{equation}
\begin{array}{l}
\displaystyle%
\mathrm{Mi}_{\kappa ,\mu }\left( x\right) =\int_{0}^{x}\frac{\mathrm{M}%
_{\kappa ,\mu }\left( t\right) }{t}dt, \\
\displaystyle%
\mathrm{mi}_{\kappa ,\mu }\left( x\right) =\int_{x}^{\infty }\frac{\mathrm{M}%
_{\kappa ,\mu }\left( t\right) }{t}dt.%
\end{array}
\label{Integral_Whittaker_M_def}
\end{equation}%
and%
\begin{equation}
\begin{array}{l}
\displaystyle%
\mathrm{Wi}_{\kappa ,\mu }\left( x\right) =\int_{0}^{x}\frac{\mathrm{W}%
_{\kappa ,\mu }\left( t\right) }{t}dt, \\
\displaystyle%
\mathrm{wi}_{\kappa ,\mu }\left( x\right) =\int_{x}^{\infty }\frac{\mathrm{W}%
_{\kappa ,\mu }\left( t\right) }{t}dt.%
\end{array}
\label{Integral_Whittaker_W_def}
\end{equation}

\begin{figure}[H]
\includegraphics[width=9.1cm]{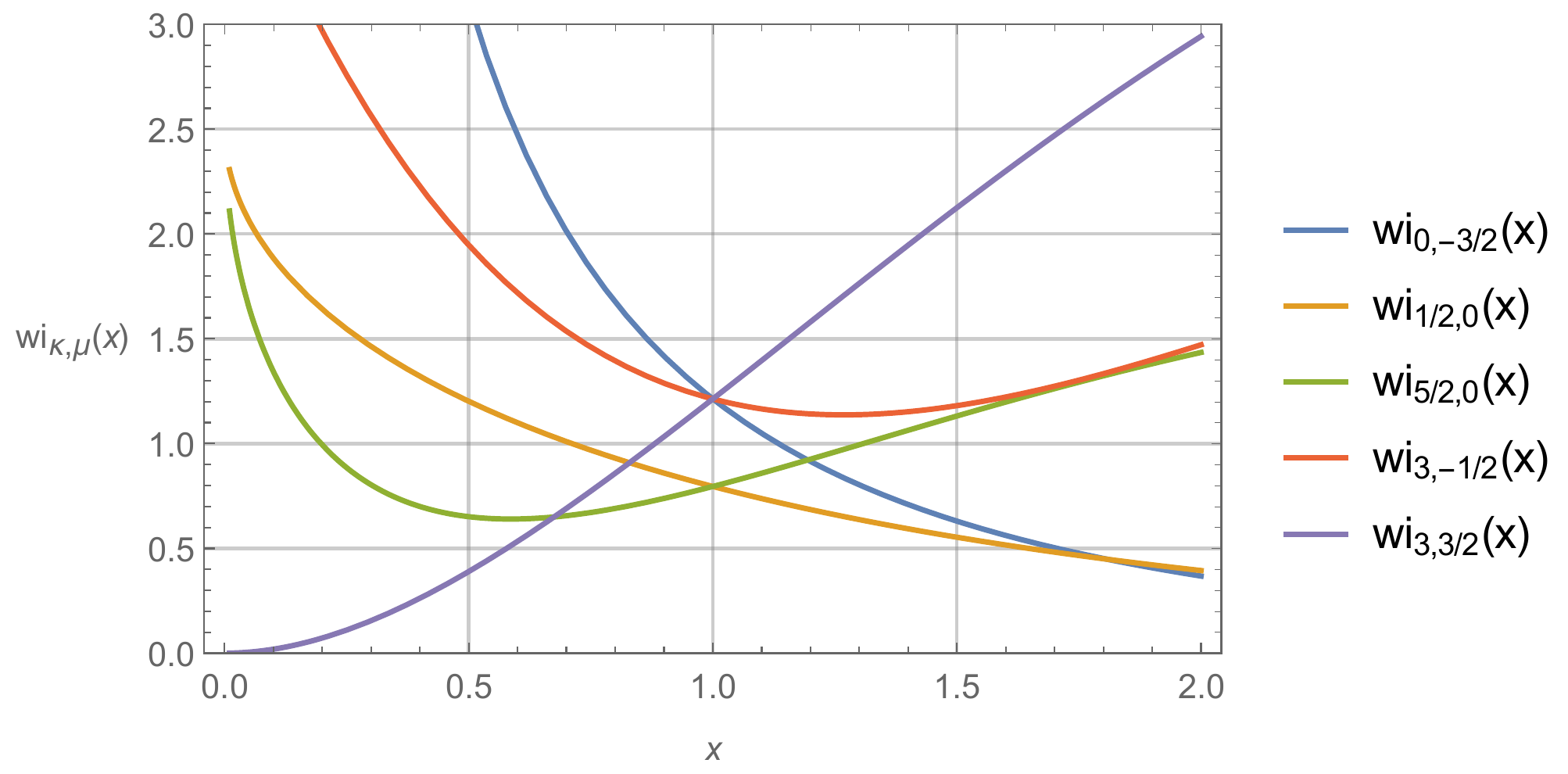}
\caption{The integral Whittaker functions $\mathrm{wi}_{\protect\kappa ,%
\protect\mu }\left( x\right) $ as a function of variable $x$ at fixed values
of parameters $\protect\kappa $ and $\protect\mu $.}
\label{Figure: wi_kappa_mu}
\end{figure}

The integral Whittaker functions with particular values of parameters $%
\kappa $ and $\mu $ can be expressed in terms of elementary and special
functions. These cases, derived using the MATHEMATICA program, are presented in
Tables \ref{Table3}--\ref{Table6}. Several integral Whittaker functions $\mathrm{Mi}_{\kappa ,\mu
}\left( x\right) $, $\mathrm{mi}_{\kappa ,\mu }\left( x\right) $, $\mathrm{Wi%
}_{\kappa ,\mu }\left( x\right) $ and $\mathrm{wi}_{\kappa ,\mu }\left(
x\right) $ as a function of variable $x$ at fixed values of parameters $%
\kappa $ and $\mu $ are plotted in Figures \ref{Figure: Mi_kappa,mu}--\ref{Figure: wi_kappa_mu}. Similarly, a long list of the Whittaker functions $\mathrm{M}_{\kappa ,\mu
}\left( x\right) $ and $\mathrm{W}_{\kappa ,\mu }\left( x\right) $ with
integer and fractional parameters was prepared (see Appendix \ref{Appendix: Whittaker}).
In some cases, it was possible to obtain for them their Laplace transforms, and
they are also reported in Appendix \ref{Appendix: Whittaker}.

\begin{figure}[H]
\includegraphics[width=9.1cm]{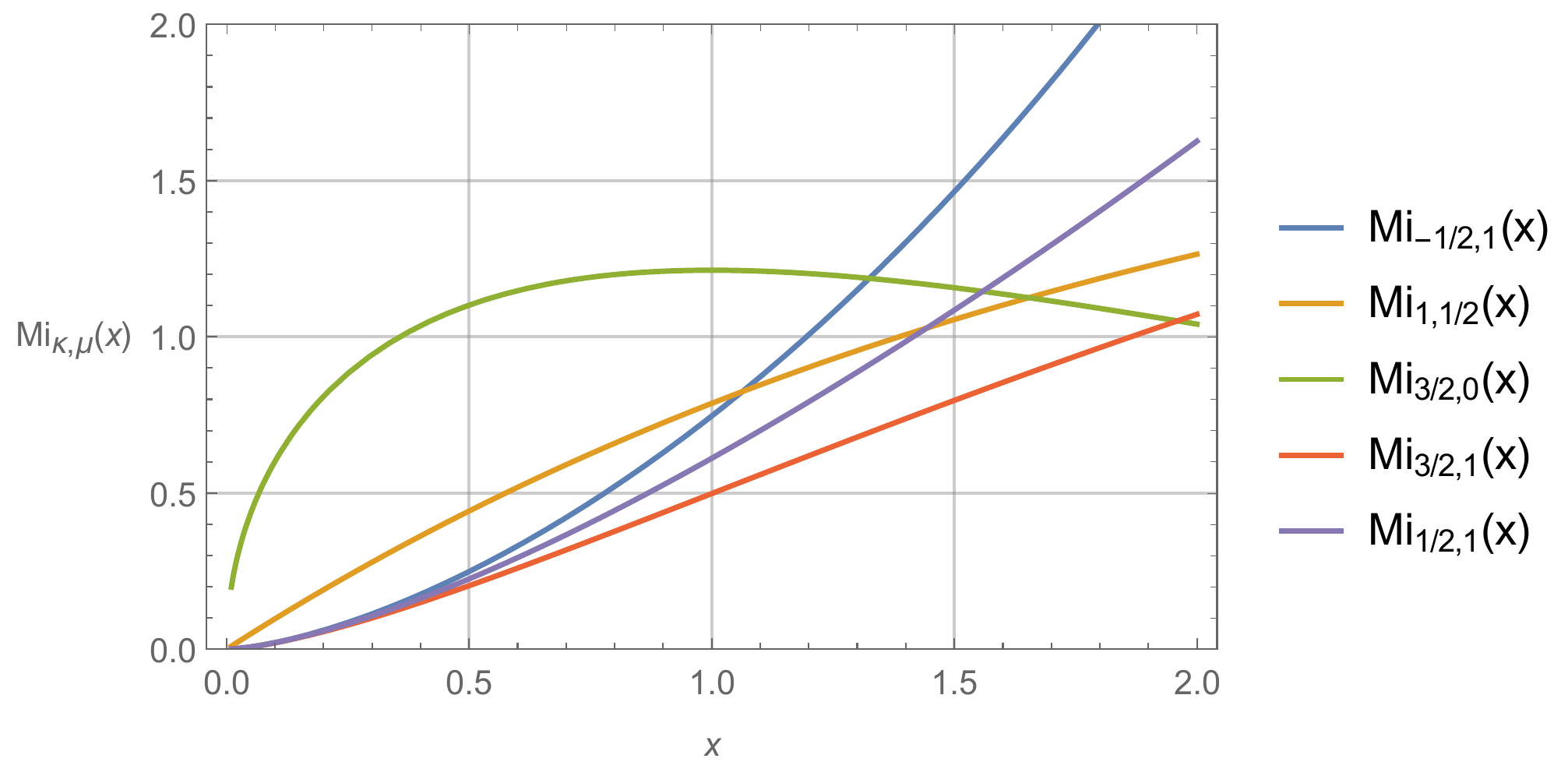}
\caption{The integral Whittaker functions $\mathrm{Mi}_{\protect\kappa ,%
\protect\mu }\left( x\right) $ as a function of variable $x$ at fixed values
of parameters $\protect\kappa $ and $\protect\mu $.}
\label{Figure: Mi_kappa,mu}
\end{figure}

\vspace{-5pt}

\begin{figure}[H]
\includegraphics[width=9.1cm]{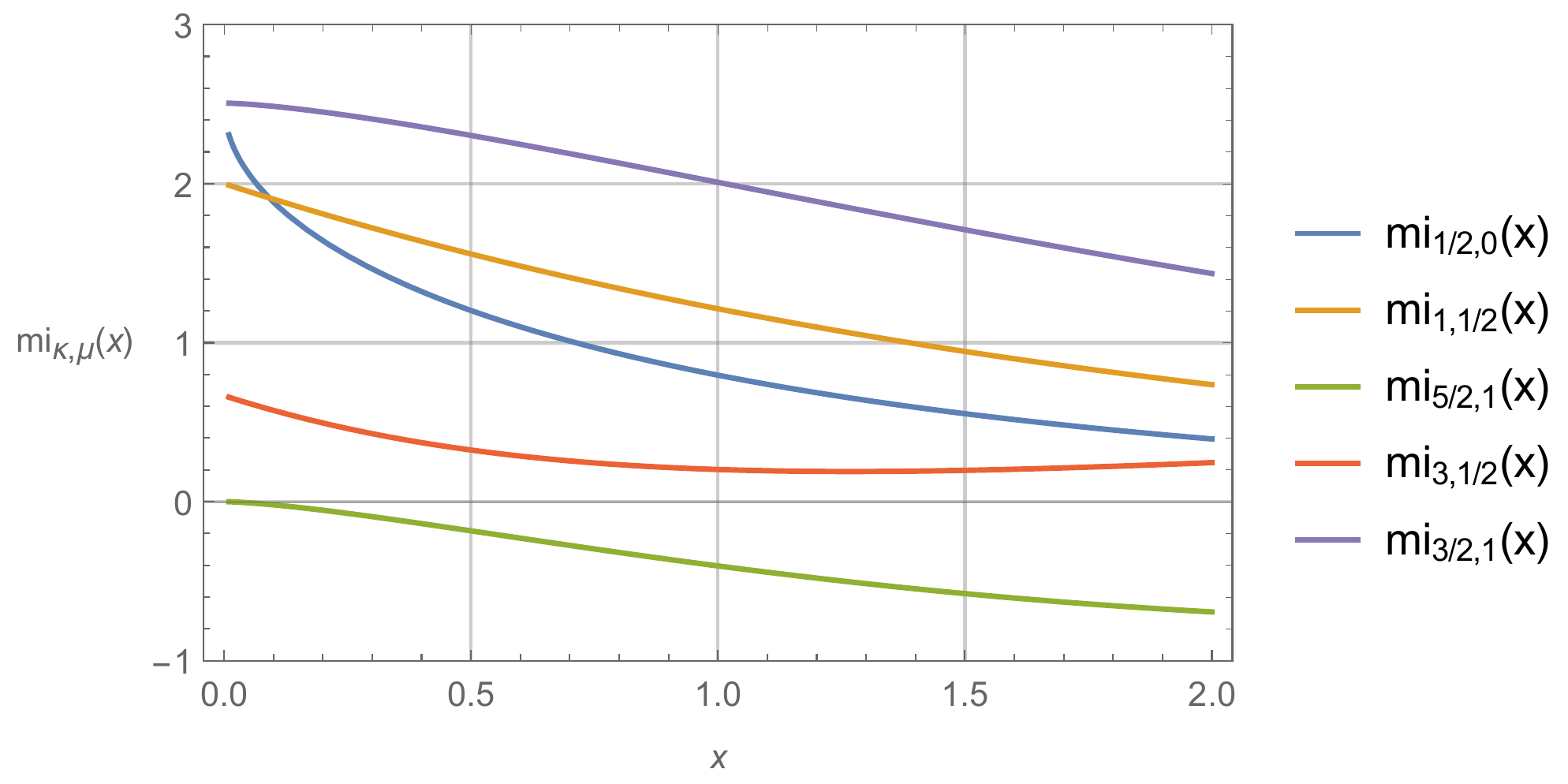}
\caption{The integral Whittaker functions $\mathrm{mi}_{\protect\kappa ,%
\protect\mu }\left( x\right) $ as a function of variable $x$ at fixed values
of parameters $\protect\kappa $ and $\protect\mu $.}
\label{Figure: mmi_kappa_mu}
\end{figure}
\vspace{-5pt}

\begin{figure}[H]
\includegraphics[width=9.1cm]{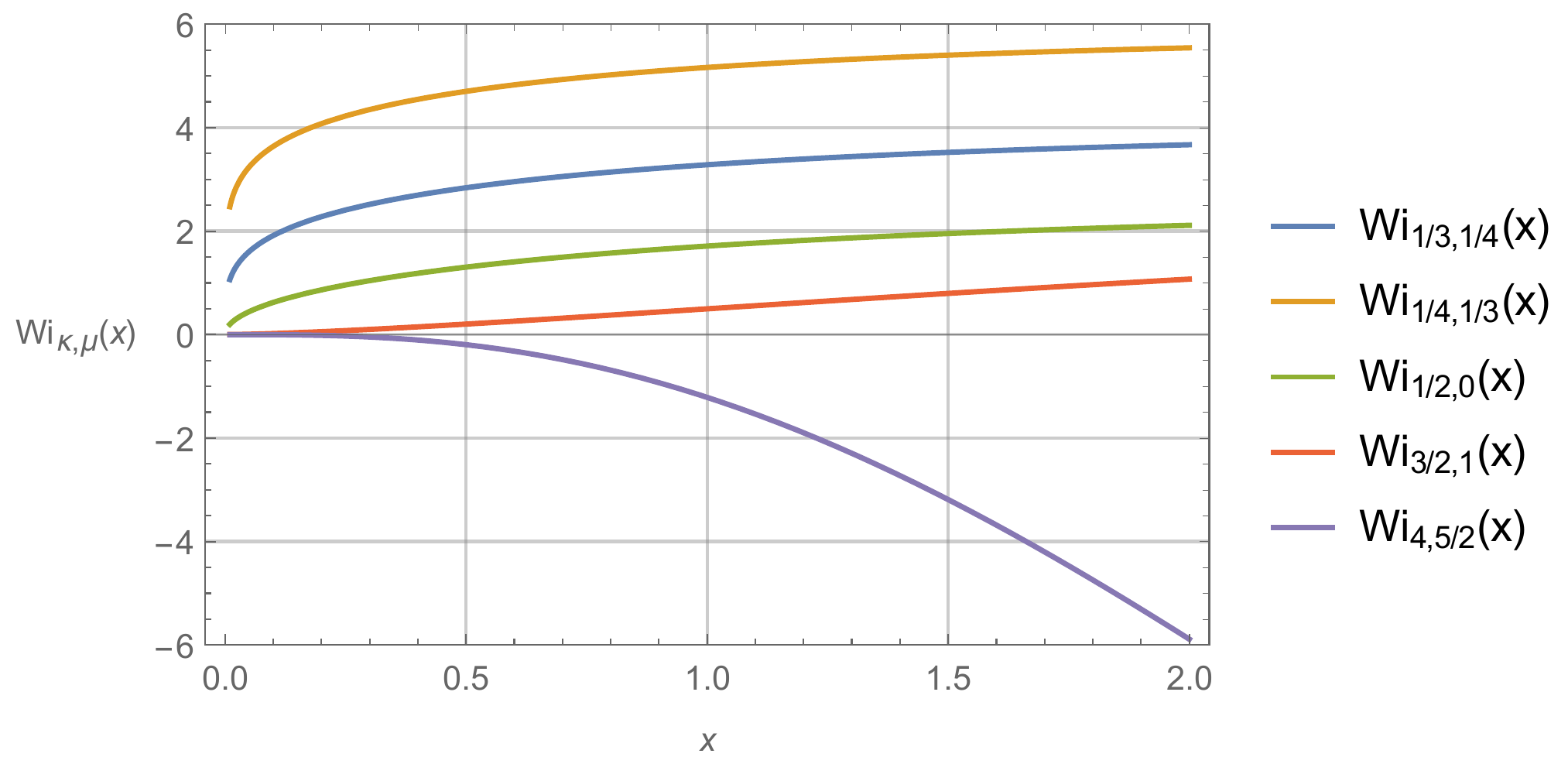}
\caption{The integral Whittaker functions $\mathrm{Wi}_{\protect\kappa ,%
\protect\mu }\left( x\right) $ as a function of variable $x$ at fixed values
of parameters $\protect\kappa $ and $\protect\mu $.}
\label{Figure: Wi_kappa_mu}
\end{figure}

\vspace{-5pt}

\begin{center}
\begin{specialtable}[H] \centering%
\caption{The integral Whittaker functions $\mathrm{Mi}_{\kappa, \mu}$ derived for some values of
parameters $\kappa$ and $\mu$ by using (\ref{Whittaker_def}) and (\ref{Integral_Whittaker_M_def}).}%
\begin{tabular*}{\hsize}{c@{\extracolsep{\fill}}cc}
\toprule
\boldmath$\kappa $ & \boldmath$\mu $ & \boldmath$\mathrm{Mi}_{\kappa ,\mu }\left( x\right) $ \\
\midrule
$-\frac{5}{2}$ & $0$ & $e^{x/2}\left[ \sqrt{x}\left( x+1\right) +2F\left(
\sqrt{\frac{x}{2}}\right) \right] ,\quad F\left( x\right)
=e^{-x^{2}}\int_{0}^{x}e^{t^{2}}dt$ \\ \midrule
$-\frac{3}{2}$ & $0$ & $2\sqrt{x}e^{x/2}$ \\ \midrule
$-\frac{3}{2}$ & $\frac{1}{2}$ & $\frac{1}{6}\left[ \left( 8x+3\pi x\,%
\mathbf{L}_{0}\left( \frac{x}{2}\right) \right) I_{1}\left( \frac{x}{2}%
\right) +\left( -\pi x\,\mathbf{L}_{1}\left( \frac{x}{2}\right) +6x+4\right)
I_{0}\left( \frac{x}{2}\right) -4\right] $ \\ \midrule
$-\frac{3}{2}$ & $1$ & $2e^{x/2}\left[ \sqrt{x}-\sqrt{2}F\left( \sqrt{\frac{x%
}{2}}\right) \right] $ \\ \midrule
$-\frac{3}{2}$ & $\frac{3}{2}$ & $\frac{2}{5}\left[ \left( -8+8x+3\pi x\,%
\mathbf{L}_{0}\left( \frac{x}{2}\right) \right) I_{1}\left( \frac{x}{2}%
\right) +\left( -3\pi x\,\mathbf{L}_{1}\left( \frac{x}{2}\right)
+2x-12\right) I_{0}\left( \frac{x}{2}\right) +12\right] $ \\ \midrule
$-\frac{3}{2}$ & $2$ & $8\left[ \sqrt{2\pi }\,\mathrm{erf}\left( \sqrt{\frac{%
x}{2}}\right) +\frac{e^{-x/2}}{x^{-3/2}}\left( 2\left( x-1\right)
+e^{x}\left( x\left( x-4\right) +2\right) \right) \right] $ \\ \midrule
$-\frac{3}{2}$ & $\frac{5}{2}$ & {\small$\frac{16}{7}\left[ 20+\left( 5\pi x\,%
\mathbf{L}_{1}\left( \frac{x}{2}\right) +18x-44\right) I_{0}\left( \frac{x}{2%
}\right) +\frac{1}{x}\left( -5\pi x^{2}\,\mathbf{L}_{0}\left( \frac{x}{2}%
\right) +8x\left( x-9\right) +96\right) I_{1}\left( \frac{x}{2}\right) %
\right] $ }\\ \midrule
$-\frac{3}{2}$ & $3$ & $30\left[ x^{-5/2}e^{-x/2}\left( 48+2e^{x}\left(
x^{3}-12x^{2}+24x-24\right) \right) +3\sqrt{2\pi }\,\mathrm{erfi}\left(
\sqrt{\frac{x}{2}}\right) \right] $ \\ \midrule
$-\frac{1}{2}$ & $0$ & $\sqrt{2\pi }\,\mathrm{erfi}\left( \sqrt{\frac{x}{2}}%
\right) $ \\ \midrule
$-\frac{1}{2}$ & $\frac{1}{2}$ & $\frac{1}{2}\left[ -\pi x\,\mathbf{L}%
_{0}\left( \frac{x}{2}\right) I_{1}\left( \frac{x}{2}\right) +\left( \pi x\,%
\mathbf{L}_{1}\left( \frac{x}{2}\right) +2x+4\right) I_{0}\left( \frac{x}{2}%
\right) -4\right] $ \\ \midrule
$-\frac{1}{2}$ & $1$ & $8x^{-1/2}\sinh \left( \frac{x}{2}\right) -2\sqrt{%
2\pi }\,\mathrm{erf}\left( \sqrt{\frac{x}{2}}\right) $ \\ \midrule
$-\frac{1}{2}$ & $\frac{3}{2}$ & $2\left( \pi x\,\mathbf{L}_{0}\left( \frac{x%
}{2}\right) +8\right) I_{1}\left( \frac{x}{2}\right) -2\left( \pi x\,\mathbf{%
L}_{1}\left( \frac{x}{2}\right) +2x-4\right) I_{0}\left( \frac{x}{2}\right)
-8$ \\ \midrule
$-\frac{1}{2}$ & $2$ & $48x^{-3/2}e^{-x/2}\left[ \left( x-1\right) e^{x}+1%
\right] -12\sqrt{2\pi }\,\mathrm{erfi}\left( \sqrt{\frac{x}{2}}\right) $ \\
\midrule
$0$ & $\frac{1}{8}$ & $\frac{8}{5}x^{5/8}\,_{1}F_{2}\left( \left.
\begin{array}{c}
\frac{5}{16} \\
\frac{9}{8},\frac{21}{16}%
\end{array}%
\right\vert \frac{x^{2}}{16}\right) $ \\ \midrule
$0$ & $\frac{1}{7}$ & $\frac{14}{9}x^{9/14}\,_{1}F_{2}\left( \left.
\begin{array}{c}
\frac{9}{28} \\
\frac{8}{7},\frac{37}{28}%
\end{array}%
\right\vert \frac{x^{2}}{16}\right) $ \\ \midrule
$0$ & $\frac{1}{6}$ & $\frac{3}{2}x^{2/3}\,_{1}F_{2}\left( \left.
\begin{array}{c}
\frac{1}{3} \\
\frac{7}{6},\frac{4}{3}%
\end{array}%
\right\vert \frac{x^{2}}{16}\right) $ \\ \midrule
$0$ & $\frac{1}{5}$ & $\frac{10}{7}x^{7/10}\,_{1}F_{2}\left( \left.
\begin{array}{c}
\frac{7}{20} \\
\frac{6}{5},\frac{27}{20}%
\end{array}%
\right\vert \frac{x^{2}}{16}\right) $ \\ \midrule
$0$ & $\frac{1}{4}$ & $\frac{4}{3}x^{3/4}\,_{1}F_{2}\left( \left.
\begin{array}{c}
\frac{3}{8} \\
\frac{5}{4},\frac{11}{8}%
\end{array}%
\right\vert \frac{x^{2}}{16}\right) $ \\ \midrule
$0$ & $\frac{1}{3}$ & $\frac{6}{5}x^{5/6}\,_{1}F_{2}\left( \left.
\begin{array}{c}
\frac{5}{12} \\
\frac{4}{3},\frac{17}{12}%
\end{array}%
\right\vert \frac{x^{2}}{16}\right) $ \\ \midrule
$0$ & $\frac{1}{2}$ & $2\,\mathrm{Shi}\left( \frac{x}{2}\right) $ \\ \midrule
$0$ & $1$ & $\frac{2}{3}x^{3/2}\,_{1}F_{2}\left( \left.
\begin{array}{c}
\frac{3}{4} \\
\frac{7}{4},2%
\end{array}%
\right\vert \frac{x^{2}}{16}\right) $ \\ \midrule
$0$ & $\frac{3}{2}$ & $\frac{24}{x}\sinh \left( \frac{x}{2}\right) -12$ \\
\midrule
$0$ & $2$ & $\frac{2}{5}x^{5/2}\,_{1}F_{2}\left( \left.
\begin{array}{c}
\frac{5}{4} \\
\frac{9}{4},3%
\end{array}%
\right\vert \frac{x^{2}}{16}\right) $ \\ \midrule
$0$ & $\frac{5}{2}$ & $60\left[ x^{-2}\left( 6x\cosh \left( \frac{x}{2}%
\right) -2\sinh \left( \frac{x}{2}\right) \right) -\mathrm{Shi}\left( \frac{x%
}{2}\right) \right] $ \\ \bottomrule
\end{tabular*}%
\label{Table3}%
\end{specialtable}%
\end{center}

\begin{specialtable}[H] \centering%
\caption{The integral Whittaker functions $\mathrm{Mi}_{\kappa, \mu}$ derived for some values of
parameters $\kappa$ and $\mu$ by using (\ref{Whittaker_def}) and (\ref{Integral_Whittaker_M_def}).}%
\begin{tabular*}{\hsize}{c@{\extracolsep{\fill}}cc}
\toprule
\boldmath$\kappa $ & \boldmath$\mu $ & \boldmath$\mathrm{Mi}_{\kappa ,\mu }\left( x\right) $ \\
\midrule
$\frac{1}{2}$ & $0$ & $\sqrt{2\pi }\,\mathrm{erf}\left( \sqrt{\frac{x}{2}}%
\right) $ \\ \midrule
$\frac{1}{2}$ & $\frac{1}{2}$ & $\frac{1}{2}\left[ -\pi x\,\mathbf{L}%
_{0}\left( \frac{x}{2}\right) I_{1}\left( \frac{x}{2}\right) +\left( \pi x\,%
\mathbf{L}_{1}\left( \frac{x}{2}\right) +2x-4\right) I_{0}\left( \frac{x}{2}%
\right) +4\right] $ \\ \midrule
$\frac{1}{2}$ & $1$ & $2\sqrt{2\pi }\,\mathrm{erfi}\left( \sqrt{\frac{x}{2}}%
\right) -8x^{-1/2}\sinh \left( \frac{x}{2}\right) $ \\ \midrule
$\frac{1}{2}$ & $\frac{3}{2}$ & $-2\left( \pi x\,\mathbf{L}_{0}\left( \frac{x%
}{2}\right) +8\right) I_{1}\left( \frac{x}{2}\right) +2\left( \pi x\,\mathbf{%
L}_{1}\left( \frac{x}{2}\right) +2x+4\right) I_{0}\left( \frac{x}{2}\right)
+8$ \\ \midrule
$\frac{1}{2}$ & $\frac{5}{2}$ & $%
\begin{array}{l}
16\left[ \left( \pi x\,\mathbf{L}_{0}\left( \frac{x}{2}\right) +8\right)
I_{1}\left( \frac{x}{2}\right) -\left( \pi x\,\mathbf{L}_{1}\left( \frac{x}{2%
}\right) +2x+4\right) I_{0}\left( \frac{x}{2}\right) +4\right] \\
+8x^{2}\,_{2}F_{3}\left( \left.
\begin{array}{c}
1,1 \\
2,2,2%
\end{array}%
\right\vert \frac{x^{2}}{16}\right) -4x^{2}\,_{2}F_{3}\left( \left.
\begin{array}{c}
1,1 \\
2,2,3%
\end{array}%
\right\vert \frac{x^{2}}{16}\right)%
\end{array}%
$ \\ \midrule
$1$ & $0$ & $\frac{\sqrt{x}}{30}\left[ 60\,_{1}F_{2}\left( \left.
\begin{array}{c}
\frac{1}{4} \\
1,\frac{5}{4}%
\end{array}%
\right\vert \frac{x^{2}}{16}\right) +3x^{2}\,_{1}F_{2}\left( \left.
\begin{array}{c}
\frac{5}{4} \\
2,\frac{9}{4}%
\end{array}%
\right\vert \frac{x^{2}}{16}\right) -20x\,_{1}F_{2}\left( \left.
\begin{array}{c}
\frac{3}{4} \\
1,\frac{7}{4}%
\end{array}%
\right\vert \frac{x^{2}}{16}\right) \right] $ \\ \midrule
$1$ & $\frac{1}{2}$ & $2\left( 1-e^{-x/2}\right) $ \\ \midrule
$1$ & $1$ & $-\frac{2x^{3/2}}{45}\left[ -20\,_{1}F_{2}\left( \left.
\begin{array}{c}
\frac{3}{4} \\
1,\frac{7}{4}%
\end{array}%
\right\vert \frac{x^{2}}{16}\right) +5x^{2}\,_{1}F_{2}\left( \left.
\begin{array}{c}
\frac{3}{4} \\
\frac{7}{4},\frac{3}{4}%
\end{array}%
\right\vert \frac{x^{2}}{16}\right) +3x\,_{1}F_{2}\left( \left.
\begin{array}{c}
\frac{5}{4} \\
2,\frac{9}{4}%
\end{array}%
\right\vert \frac{x^{2}}{16}\right) \right] $ \\ \midrule
$\frac{3}{2}$ & $0$ & $2\sqrt{x}e^{-x/2}$ \\ \midrule
$\frac{3}{2}$ & $\frac{1}{2}$ & $\frac{1}{6}\left[ x\left( \pi \,\mathbf{L}%
_{0}\left( \frac{x}{2}\right) -8\right) I_{1}\left( \frac{x}{2}\right)
+\left( \pi x\,\mathbf{L}_{1}\left( \frac{x}{2}\right) -6x+4\right)
I_{0}\left( \frac{x}{2}\right) +4\right] $ \\ \midrule
$\frac{3}{2}$ & $1$ & $\sqrt{2\pi }\,\mathrm{erf}\left( \sqrt{\frac{x}{2}}%
\right) -2\sqrt{x}e^{-x/2}$ \\ \midrule
$\frac{3}{2}$ & $\frac{3}{2}$ & $\frac{1}{5}\left[ 2\left( -3\pi \,\mathbf{L}%
_{0}\left( \frac{x}{2}\right) +8x+8\right) I_{1}\left( \frac{x}{2}\right)
+\left( 6\pi x\,\mathbf{L}_{1}\left( \frac{x}{2}\right) -4x+24\right)
I_{0}\left( \frac{x}{2}\right) +24\right] $ \\ \midrule
$2$ & $\frac{1}{2}$ & $xe^{-x/2}$ \\ \midrule
$2$ & $\frac{3}{2}$ & $4-2\left( 2+x\right) e^{-x/2}$ \\ \midrule
$2$ & $\frac{5}{2}$ & $-5x^{-2}\left( 6\left( 2+x\right) e^{x}-2x\left(
3+x\right) ^{2}-12\right) +30\,\mathrm{Shi}\left( \frac{x}{2}\right) $ \\
\bottomrule
\end{tabular*}%
\label{Table3a}%
\end{specialtable}%
\vspace{-12pt}
\begin{specialtable}[H] \centering%
\caption{The integral Whittaker functions $\mathrm{mi}_{\kappa, \mu}$ derived for some values of
parameters $\kappa$ and $\mu$ by using (\ref{Whittaker_def}) and (\ref{Integral_Whittaker_M_def}).}%
\begin{tabular*}{\hsize}{c@{\extracolsep{\fill}}cc}
\toprule
\boldmath$\kappa $ & \boldmath$\mu $ & \boldmath$\mathrm{mi}_{\kappa ,\mu }\left( x\right) $ \\
\midrule
$\frac{1}{2}$ & $0$ & $\sqrt{2\pi }\,\mathrm{erfc}\left( \sqrt{\frac{x}{2}}%
\right) $ \\ \midrule
$1$ & $\frac{1}{2}$ & $2e^{-x/2}$ \\ \midrule
$\frac{3}{2}$ & $1$ & $2\sqrt{x}e^{-x/2}+\sqrt{2\pi }\,\mathrm{erfc}\left(
\sqrt{\frac{x}{2}}\right) $ \\ \midrule
$2$ & $\frac{1}{2}$ & $-e^{-x/2}$ \\ \midrule
$2$ & $\frac{3}{2}$ & $2\left( 2+x\right) e^{-x/2}$ \\ \midrule
$\frac{5}{2}$ & $1$ & $-\frac{2}{3}x^{3/2}e^{-x/2}$ \\ \midrule
$\frac{5}{2}$ & $2$ & $2\sqrt{x}\left( 3+x\right) e^{-x/2}+3\sqrt{2\pi }\,%
\mathrm{erfc}\left( \sqrt{\frac{x}{2}}\right) $ \\ \midrule
$3$ & $\frac{1}{2}$ & $\frac{1}{3}\left[ 2+\left( x-2\right) x\right]
e^{-x/2}$ \\ \midrule
$3$ & $\frac{3}{2}$ & $-\frac{1}{2}x^{2}e^{-x/2}$ \\ \midrule
$4$ & $\frac{1}{2}$ & $-\frac{x}{12}\left[ 12+\left( x-6\right) x\right]
e^{-x/2}$ \\ \midrule
$4$ & $\frac{3}{2}$ & $\frac{1}{10}\left[ 8+\left( x-2\right) ^{2}x\right]
e^{-x/2}$ \\ \bottomrule
\end{tabular*}%
\label{Table4}%
\end{specialtable}%


\begin{specialtable}[H] \centering%
\caption{The integral Whittaker functions $\mathrm{Wi}_{\kappa, \mu}$ derived for some values of
parameters $\kappa$ and $\mu$ by using (\ref{Whittaker_def}) and (\ref{Integral_Whittaker_W_def}).}%
\begin{tabular*}{\hsize}{c@{\extracolsep{\fill}}cc}
\toprule
\boldmath$\kappa $ & \boldmath$\mu $ & \boldmath$\mathrm{Wi}_{\kappa ,\mu }\left( x\right) $ \\
\midrule
$-\frac{3}{2}$ & $0$ & $-2\sqrt{x}e^{-x/2}\mathrm{Ei}\left( -x\right) +\sqrt{%
2\pi }\,\mathrm{erfc}\left( \sqrt{\frac{x}{2}}\right) $ \\ \midrule
$\frac{1}{4}$ & $\frac{1}{4}$ & $2^{1/4}\gamma \left( \frac{1}{4},\frac{x}{2}%
\right) $ \\ \midrule
$\frac{1}{2}$ & $0$ & $\sqrt{2\pi }\,\mathrm{erf}\left( \sqrt{\frac{x}{2}}%
\right) $ \\ \midrule
$1$ & $-\frac{1}{2}$ & $2\left( 1-e^{-x/2}\right) $ \\ \midrule
$1$ & $\frac{1}{2}$ & $2\left( 1-e^{-x/2}\right) $ \\ \midrule
$\frac{3}{2}$ & $0$ & $-2\sqrt{x}e^{-x/2}$ \\ \midrule
$\frac{3}{2}$ & $1$ & $\sqrt{2\pi }\,\mathrm{erf}\left( \sqrt{\frac{x}{2}}%
\right) -2\sqrt{x}e^{-x/2}$ \\ \midrule
$2$ & $-\frac{1}{2}$ & $-2xe^{-x/2}$ \\ \midrule
$2$ & $\frac{1}{2}$ & $-2xe^{-x/2}$ \\ \midrule
$2$ & $\frac{3}{2}$ & $4-2\left( 2+x\right) e^{-x/2}$ \\ \midrule
$\frac{5}{2}$ & $0$ & $\sqrt{2\pi }\,\mathrm{erf}\left( \sqrt{\frac{x}{2}}%
\right) -2\sqrt{x}\left( x-1\right) e^{-x/2}$ \\ \midrule
$3$ & $\frac{3}{2}$ & $-2x^{2}e^{-x/2}$ \\ \midrule
$3$ & $\frac{5}{2}$ & $16-2\left[ 8+\left( 4+x\right) x\right] e^{-x/2}$ \\
\midrule
$4$ & $-\frac{1}{2}$ & $-2x\left[ 12+\left( x-6\right) x\right] e^{-x/2}$ \\
\midrule
$4$ & $\frac{1}{2}$ & $-2x\left[ 12+\left( x-6\right) x\right] e^{-x/2}$ \\
\midrule
$4$ & $\frac{3}{2}$ & $16-2\left[ 8+x\left( x-2\right) ^{2}\right] e^{-x/2}$
\\ \midrule
$4$ & $\frac{5}{2}$ & $-2x^{3}e^{-x/2}$ \\ \bottomrule
\end{tabular*}%
\label{Table5}%
\end{specialtable}%

\vspace{-12pt}

\begin{specialtable}[H] \centering%
\caption{The integral Whittaker functions $\mathrm{wi}_{\kappa, \mu}$ derived for some values of
parameters $\kappa$ and $\mu$ by using (\ref{Whittaker_def}) and (\ref{Integral_Whittaker_W_def}).}%
\begin{tabular*}{\hsize}{c@{\extracolsep{\fill}}cc}
\toprule
$\kappa $ & $\mu $ & $\mathrm{wi}_{\kappa ,\mu }\left( x\right) $ \\
\midrule
$-\frac{1}{2}$ & $1$ & $2x^{-1/2}e^{-x/2}-\sqrt{2\pi }\,\mathrm{erfc}\left(
\sqrt{\frac{x}{2}}\right) $ \\ \midrule
$-\frac{1}{2}$ & $2$ & $2x^{-3/2}e^{-x/2}$ \\ \midrule
$-\frac{1}{2}$ & $3$ & $\frac{1}{3}\left[ -2\left( x-6\right) \left(
x+2\right) x^{-5/2}e^{-x/2}+\sqrt{2\pi }\,\mathrm{erfc}\left( \sqrt{\frac{x}{%
2}}\right) \right] $ \\ \midrule
$0$ & $-\frac{5}{2}$ & $\frac{1}{2}\mathrm{Ei}\left( -\frac{x}{2}\right)
+3x^{-2}\left( x+2\right) e^{-x/2}$ \\ \midrule
$0$ & $-\frac{3}{2}$ & $2x^{-1}e^{-x/2}$ \\ \midrule
$0$ & $-\frac{1}{2}$ & $-\mathrm{Ei}\left( -\frac{x}{2}\right) $ \\ \midrule
$\frac{1}{2}$ & $0$ & $\sqrt{2\pi }\,\mathrm{erfc}\left( \sqrt{\frac{x}{2}}%
\right) $ \\ \midrule
$1$ & $-\frac{1}{2}$ & $2e^{-x/2}$ \\ \midrule
$1$ & $\frac{1}{2}$ & $2e^{-x/2}$ \\ \midrule
$1$ & $\frac{5}{2}$ & $2x^{-2}\left[ 6+x\left( 6+x\right) \right] e^{-x/2}$
\\ \midrule
$\frac{3}{2}$ & $0$ & $2\sqrt{x}e^{-x/2}$ \\

\bottomrule
\end{tabular*}
\end{specialtable}

\vspace{-12pt}
\begin{specialtable}[H]\ContinuedFloat
\small
\caption{{\em Cont.}}
\begin{tabular*}{\hsize}{c@{\extracolsep{\fill}}cc}
\toprule
$\kappa $ & $\mu $ & $\mathrm{wi}_{\kappa ,\mu }\left( x\right) $ \\
\midrule

$\frac{3}{2}$ & $1$ & $\sqrt{2\pi }\,\mathrm{erfc}\left( \sqrt{\frac{x}{2}}%
\right) +2\sqrt{x}e^{-x/2}$ \\ \midrule
$2$ & $\frac{1}{2}$ & $2xe^{-x/2}$ \\ \midrule
$2$ & $\frac{3}{2}$ & $2\left( 2+x\right) e^{-x/2}$ \\ \midrule
$\frac{5}{2}$ & $0$ & $\sqrt{2\pi }\,\mathrm{erfc}\left( \sqrt{\frac{x}{2}}%
\right) +2\sqrt{x}\left( x-1\right) e^{-x/2}$ \\ \midrule
$\frac{5}{2}$ & $1$ & $2x^{3/2}e^{-x/2}$ \\ \midrule
$3$ & $-\frac{1}{2}$ & $2\left[ 2+x\left( x-2\right) \right] e^{-x/2}$ \\
\midrule
$3$ & $\frac{1}{2}$ & $2\left[ 2+x\left( x-2\right) \right] e^{-x/2}$ \\
\midrule
$3$ & $\frac{3}{2}$ & $2x^{2}e^{-x/2}$ \\ \midrule
$3$ & $\frac{5}{2}$ & $2\left[ 8+x\left( 4+x\right) \right] e^{-x/2}$ \\
\midrule
$4$ & $-\frac{1}{2}$ & $2x\left[ 12+\left( x-6\right) x\right] e^{-x/2}$ \\
\midrule
$4$ & $\frac{1}{2}$ & $2x\left[ 12+\left( x-6\right) x\right] e^{-x/2}$ \\
\midrule
$4$ & $\frac{3}{2}$ & $2\left[ 8+\left( x-2\right) ^{2}x\right] e^{-x/2}$ \\
\bottomrule
\end{tabular*}%
\label{Table6}%
\end{specialtable}%

There is a number of recurrence relations between the Whittaker functions,
for example~\cite{Abramowitz,Magnus}%
\begin{equation}
\begin{array}{c}
2\mu \left[ \mathrm{M}_{\kappa -1/2,\mu -1/2}\left( t\right) -\mathrm{M}%
_{\kappa +1/2,\mu -1/2}\left( t\right) \right] =t^{1/2}\mathrm{M}_{\kappa
,\mu }\left( t\right) , \\
\left( \kappa +\mu \right) \mathrm{W}_{\kappa -1/2,\mu }\left( t\right) +%
\mathrm{W}_{\kappa +1/2,\mu }\left( t\right) =t^{1/2}\mathrm{W}_{\kappa ,\mu
+1/2}\left( t\right) ,%
\end{array}
\label{Recurrence_Whittaker}
\end{equation}%
and this leads to integrals that are expressed in terms of the integral
Whittaker functions%
\begin{equation}
\begin{array}{c}
\displaystyle%
\int_{0}^{x}\frac{\mathrm{M}_{\kappa ,\mu }\left( t\right) }{t^{1/2}}dt=2\mu %
\left[ \mathrm{Mi}_{\kappa -1/2,\mu -1/2}\left( t\right) -\mathrm{Mi}%
_{\kappa +1/2,\mu -1/2}\left( t\right) \right] , \\
\displaystyle%
\int_{0}^{x}\frac{\mathrm{W}_{\kappa ,\mu +1/2}\left( t\right) }{t^{1/2}}%
dt=\left( \kappa +\mu \right) \mathrm{Wi}_{\kappa -1/2,\mu }\left( t\right) +%
\mathrm{Wi}_{\kappa +1/2,\mu }\left( t\right) .%
\end{array}
\label{Recurrence_Whittaker_Integrals}
\end{equation}

Using the following representation of the Whittaker functions \cite{DLMF}%
\begin{equation}
\mathrm{M}_{\kappa ,\mu }\left( t\right) =t^{\mu +1/2}\sum_{n=0}^{\infty
}\,_{2}F_{1}\left( \left.
\begin{array}{c}
-n,\mu -\kappa +\frac{1}{2} \\
1+2\mu%
\end{array}%
\right\vert 2\right) \frac{\left( -t/2\right) ^{n}}{n!},
\label{Alternative_Whittaker_representation}
\end{equation}%
it is possible to obtain the integral Whittaker functions in terms of a
rapidly convergent alternating series as follows:%
\begin{equation}
\mathrm{Mi}_{\kappa ,\mu }\left( x\right) =x^{\mu +1/2}\sum_{n=0}^{\infty
}\,_{2}F_{1}\left( \left.
\begin{array}{c}
-n,\mu -\kappa +\frac{1}{2} \\
1+2\mu%
\end{array}%
\right\vert 2\right) \frac{\left( -x/2\right) ^{n}}{n!\left( \frac{1}{2}+\mu
+n\right) }.  \label{Mi_alternating}
\end{equation}

There is a number of particular cases where the integral Whittaker functions
can be written in a closed-form, for example, from \cite{DLMF}%
\begin{equation}
\mathrm{M}_{\kappa ,\kappa -1/2}\left( x\right) =x^{\kappa }e^{-x/2},
\label{Whittaker_reduction_1}
\end{equation}%
we have%
\begin{equation}
\mathrm{Mi}_{\kappa ,\kappa -1/2}\left( x\right) =2^{\kappa }\,\gamma \left(
\kappa ,\frac{x}{2}\right) ,  \label{Whittaker_Integral_reduction_1}
\end{equation}%
but \cite{DLMF}%
\begin{equation}
\mathrm{M}_{\kappa ,\kappa -1/2}\left( x\right) =\mathrm{W}_{\kappa ,\kappa
-1/2}\left( x\right) =\mathrm{W}_{\kappa ,-\kappa +1/2}\left( x\right)
=e^{-x/2}x^{\kappa },  \label{Whittaker_reduction_2}
\end{equation}%
and therefore%
\begin{equation}
\mathrm{Mi}_{\kappa ,\kappa -1/2}\left( x\right) =\mathrm{Wi}_{\kappa
,\kappa -1/2}\left( x\right) =\mathrm{Wi}_{\kappa ,-\kappa +1/2}\left(
x\right) =2^{\kappa }\,\gamma \left( \kappa ,\frac{x}{2}\right) .
\label{Whittaker_Integral_reduction_2}
\end{equation}

Furthermore, from \cite{DLMF}%
\begin{equation}
\mathrm{M}_{0,\mu }\left( x\right) =2^{2\mu +1/2}\Gamma \left( \mu +1\right)
\sqrt{\frac{t}{2}}\ I_{\mu }\left( \frac{x}{2}\right) ,
\label{Whittaker_reduction_0_mu}
\end{equation}%
follows that%
\begin{equation}
\mathrm{Mi}_{0,\mu }\left( x\right) =\frac{x^{\mu +1/2}}{\mu +1/2}%
\,_{1}F_{2}\left( \left.
\begin{array}{c}
\frac{2\mu +1}{4} \\
\mu +1,\frac{2\mu +5}{4}%
\end{array}%
\right\vert \frac{x^{2}}{16}\right) .
\label{Whittaker_Integral_reduction_0_mu}
\end{equation}

Similary from%
\begin{equation}
\mathrm{W}_{0,\mu }\left( x\right) =\sqrt{\frac{x}{\pi }}K_{\mu }\left(
\frac{x}{2}\right) ,  \label{Whittaker_2_reduction_0_mu}
\end{equation}%
we have%
\begin{equation}
\mathrm{Wi}_{0,\mu }\left( x\right) =\frac{\sqrt{\pi }}{2\sin \pi \mu }\left[
\frac{4^{\mu }\,\mathrm{Mi}_{0,-\mu }\left( x\right) }{\Gamma \left( 1-\mu
\right) }-\frac{4^{-\mu }\,\mathrm{Mi}_{0,\mu }\left( x\right) }{\Gamma
\left( 1+\mu \right) }\right] ,  \label{Whittaker_2_Integral_reduction_0_mu}
\end{equation}%
and in a general case%
\begin{equation}
\mathrm{Wi}_{\kappa ,\mu }\left( x\right) =\frac{\Gamma \left( -2\mu \right)
\mathrm{Mi}_{\kappa ,\mu }\left( x\right) }{\Gamma \left( \frac{1}{2}-\kappa
-\mu \right) }+\frac{\Gamma \left( 2\mu \right) \mathrm{Mi}_{\kappa ,-\mu
}\left( x\right) }{\Gamma \left( \frac{1}{2}-\kappa +\mu \right) }.
\label{Whittaker_2_Integral_reduction_kappa_mu}
\end{equation}

For $\kappa =\pm 1/2$, it is possible to obtain%
\begin{equation}
\mathrm{Wi}_{\pm \frac{1}{2},\mu }\left( x\right) =\mathrm{F}_{\mu }^{\pm
}\left( x\right) +\mathrm{F}_{-\mu }^{\pm }\left( x\right) ,
\label{Whittaker_2_Integral_reduction_1/2_mu}
\end{equation}%
where we have set
\begin{eqnarray}
\mathrm{F}_{\mu }^{\pm }\left( x\right) &=&\frac{2x^{1/2+\mu }\Gamma \left(
-2\mu \right) }{\left( 1+2\mu \right) \Gamma \left( \frac{1}{2}\mp \frac{1}{2%
}-\mu \right) }  \label{F_mu(x)_def} \\
&&\left[ _{1}F_{2}\left( \left.
\begin{array}{c}
\frac{1}{4}+\frac{\mu }{2} \\
\frac{1}{2}+\mu ,\frac{3}{4}+\frac{\mu }{2}%
\end{array}%
\right\vert \frac{x^{2}}{16}\right) \mp \frac{x/2}{3+2\mu }\,_{1}F_{2}\left(
\left.
\begin{array}{c}
\frac{3}{4}+\frac{\mu }{2} \\
\frac{3}{2}+\mu ,\frac{7}{4}+\frac{\mu }{2}%
\end{array}%
\right\vert \frac{x^{2}}{16}\right) \right] .  \nonumber
\end{eqnarray}

Since \cite{Prudnikov}%
\begin{eqnarray}
&&_{2}F_{1}\left( \left.
\begin{array}{c}
-n,\lambda \\
2\lambda +1%
\end{array}%
\right\vert 2\right)  \label{Prudnikov_+1/2} \\
&=&\frac{\Gamma \left( \lambda +\frac{1}{2}\right) }{\sqrt{\pi }}\left[
\left( \frac{1+\left( -1\right) ^{n}}{2}\right) \frac{\Gamma \left( \frac{n+1%
}{2}\right) }{\Gamma \left( \lambda +\frac{n+1}{2}\right) }+\left( \frac{%
1-\left( -1\right) ^{n}}{2}\right) \frac{\Gamma \left( \frac{n}{2}+1\right)
}{\Gamma \left( \lambda +\frac{n}{2}+1\right) }\right] ,  \nonumber
\end{eqnarray}%
and%
\begin{eqnarray}
&&_{2}F_{1}\left( \left.
\begin{array}{c}
-n,\lambda \\
2\lambda -1%
\end{array}%
\right\vert 2\right)  \label{Prudnikov_-1/2} \\
&=&\frac{\Gamma \left( \lambda -\frac{1}{2}\right) }{\sqrt{\pi }}\left[
\left( \frac{1+\left( -1\right) ^{n}}{2}\right) \frac{\Gamma \left( \frac{n+1%
}{2}\right) }{\Gamma \left( \lambda +\frac{n-1}{2}\right) }-\left( \frac{%
1-\left( -1\right) ^{n}}{2}\right) \frac{\Gamma \left( \frac{n}{2}+1\right)
}{\Gamma \left( \lambda +\frac{n}{2}\right) }\right] ,  \nonumber
\end{eqnarray}%
by introducing $\lambda =\mu $ and $\lambda =\mu +1$, after some steps, it
leads to%
\begin{eqnarray}
&&\mathrm{Mi}_{\pm \frac{1}{2},\mu }\left( x\right)  \label{Mi_1/2_mu} \\
&=&\frac{x^{\mu +1/2}}{\mu +1/2}\,\left[ _{1}F_{2}\left( \left.
\begin{array}{c}
\frac{\mu }{2}+\frac{1}{4} \\
\mu +\frac{1}{2},\frac{\mu }{2}+\frac{3}{4}%
\end{array}%
\right\vert \frac{x^{2}}{16}\right) \mp \frac{x/2}{2\mu +3}\,_{1}F_{2}\left(
\left.
\begin{array}{c}
\frac{\mu }{2}+\frac{3}{4} \\
\mu +\frac{3}{2},\frac{\mu }{2}+\frac{7}{4}%
\end{array}%
\right\vert \frac{x^{2}}{16}\right) \right] .  \nonumber
\end{eqnarray}%


\section{The Integral Wright Functions} \label{Section_Wright}

In 1933 \cite{Wright} and in 1940 \cite{Wright2}, Wright introduced new
special functions that were considered as a kind of generalization of the
Bessel functions. However, today they play a significant independent role in
mathematics and in solutions of physical problems by modeling space
diffusion, stochastic processes, probability distributions and other diverse
natural phenomena \cite{Gorenflo,MainardiBook}. The Wright functions are
defined by the following series%
\begin{equation}
\mathrm{W}_{\alpha ,\beta }\left( x\right) =\sum_{k=0}^{\infty }\frac{x^{k}}{%
k!\,\Gamma \left( \alpha k+\beta \right) }.  \label{Wright_def}
\end{equation}

If the parameter $\alpha $ is a positive real number, they are called the
Wright functions of the first kind, and when $-1<\alpha <0$, the Wright
functions of the second kind.

Furthermore, consider the following functions:\
\begin{equation}
\begin{array}{l}
\mathrm{F}_{\alpha }\left( x\right) =\mathrm{W}_{-\alpha ,0}\left( -x\right)
,\quad 0<\alpha <1, \\
\mathrm{M}_{\alpha }\left( x\right) =\mathrm{W}_{-\alpha ,1-\alpha }\left(
-x\right) ,\quad 0<\alpha <1, \\
\mathrm{F}_{\alpha }\left( x\right) =\alpha \,x\,M_{\alpha }\left( x\right) .%
\end{array}
\label{Mainardi_def}
\end{equation}

These functions with negative arguments $x$ and with particular values of
parameters are frequently named as the Mainardi functions and are denoted as
$\mathrm{F}_{\alpha }\left( x\right) $ and \linebreak $\mathrm{M}_{\alpha }\left(
x\right) $ \cite{Gorenflo,MainardiBook}.

Their explicit form is%
\begin{equation}
\begin{array}{l}
\displaystyle%
\mathrm{F}_{\alpha }\left( x\right) =\sum_{k=1}^{\infty }\frac{\left(
-x\right) ^{k}}{k!\Gamma \left( -\alpha k\right) } \\
\displaystyle%
=-\frac{1}{\pi }\sum_{k=1}^{\infty }\frac{\left( -x\right) ^{k}}{k!}\Gamma
\left( \alpha k+1\right) \sin \left( \pi \alpha k\right) , \\
\displaystyle%
\mathrm{M}_{\alpha }\left( x\right) =\sum_{k=0}^{\infty }\frac{\left(
-x\right) ^{k}}{k!\Gamma \left( -\alpha \left( k+1\right) +1\right) } \\
\displaystyle%
=\frac{1}{\pi }\sum_{k=0}^{\infty }\frac{\left( -x\right) ^{k}}{k!}\Gamma
\left( \alpha \left( k+1\right) \right) \sin \left( \pi \alpha \left(
k+1\right) \right) .%
\end{array}
\label{Mainardi_series}
\end{equation}

For positive rational $\alpha =p/q$, where $p,q$ are positive coprimes, we
have obtained reduction formulas for $\mathrm{F}_{p/q}\left( x\right) $ and $%
\mathrm{M}_{p/q}\left( x\right) $ in Appendix \ref{Appendix: Wright}. Furthermore,
by applying the MATHEMATICA\ program to sums of infinite series in (\ref%
{Wright_def}), it is possible to obtain the Wright functions of the first
and second kinds for particular values of parameters $\alpha $ and $\beta $
in an explicit form (Appendix \ref{Appendix: Wright}). The Laplace
transforms of these functions are expressed in terms of the Mittag-Leffler
functions, so they are omitted here \cite%
{Gorenflo,MainardiBook,ApelblatDifferentiation}.

The two-parameter $\mathrm{E}_{\alpha ,\beta }\left( t\right) $
Mittag-Leffler functions defined in (\ref{Mittag_Leffler_def}) differ only
by the absence of factorials from the Wright functions and, 
 therefore, the
form of series in (\ref{Wright_def})\ leads to the integral Wright function,
which is similar to that introduced in (\ref{Integral_Mittag-Leffler_def})\
and (\ref{Integral_Mittag-Leffler_series}).
\begin{equation}
\mathrm{Wi}_{\alpha ,\beta }\left( x\right) =\int_{0}^{x}\frac{\mathrm{W}%
_{\alpha ,\beta }\left( t\right) -1/\Gamma \left( \beta \right) }{t}dt.
\label{Phi_i_alpha_beta_def}
\end{equation}

Unfortunately, the notation is the same as the integral Whittaker functions.
In an explicit form from (\ref{Wright_def}), we have
\begin{equation}
\mathrm{Wi}_{\alpha ,\beta }\left( x\right) =\sum_{k=1}^{\infty }\frac{x^{k}%
}{k\ k!\,\Gamma \left( \alpha k+\beta \right) }.  \label{Phi_i_series}
\end{equation}

For $p$ and $q$ positive coprimes, applying (\ref{Sum_split})\ and (\ref%
{Multiplication_Gamma}), the corresponding expression to (\ref%
{Integral_Mittag-Leffler_p/q}) is
\begin{eqnarray}
&&\mathrm{Wi}_{p/q,\beta }\left( x\right)  \label{Phi_i_(x)_p/q} \\
&=&\sum_{k=1}^{q-1}\frac{x^{k}}{k\,k!\Gamma \left( \frac{p}{q}k+\beta
\right) }\,_{2}F_{p+q}\left( \left.
\begin{array}{c}
1,k/q \\
b_{0},\ldots ,b_{p-1},c_{0},\ldots ,c_{q-1}%
\end{array}%
\right\vert \frac{x^{q}}{p^{p}q^{q}}\right) ,  \nonumber \\
&&b_{j}=\frac{k}{q}+\frac{\beta +j}{p},\ \quad c_{j}=\frac{k+1+j}{q}.
\nonumber
\end{eqnarray}

In the case of the Mainardi functions, we have%
\begin{equation}
\begin{array}{l}
\displaystyle%
\mathrm{Fi}_{p/q}\left( x\right) =-\frac{1}{\pi }\sum_{k=1}^{q}\frac{\left(
-x\right) ^{k}}{k\ k!}\Gamma \left( \frac{p}{q}k+1\right) \sin \left( \pi
\frac{p}{q}k\right) S_{k}\left( x\right) , \\
\displaystyle%
\mathrm{Mi}_{p/q}\left( x\right) =\frac{1}{\pi }\sum_{k=1}^{q}\frac{\left(
-x\right) ^{k}}{k\ k!}\sin \left( \pi \frac{p}{q}\left( k+1\right) \right)
\Gamma \left( \frac{p}{q}\left( k+1\right) \right) S_{k}\left( x\right) ,%
\end{array}
\label{Integral_Mainardi_p/q}
\end{equation}%
where%
\begin{eqnarray}
&&S_{k}\left( x\right) =\,_{p+2}F_{q+1}\left( \left.
\begin{array}{c}
1,\frac{k}{q},a_{0},\ldots ,a_{p-1} \\
\frac{k}{q}+1,b_{0},\ldots ,b_{q-1}%
\end{array}%
\right\vert \frac{\left( -1\right) ^{p+q}x^{q}p^{p}}{q^{q}}\right) ,
\label{S_k(x)_def} \\
&&a_{j}=\frac{k}{q}+\frac{j+1}{p},\quad b_{j}=\frac{k+1+j}{q}.  \nonumber
\end{eqnarray}

In Tables \ref{Table7} and \ref{Table7a}, the integral Wright functions
derived with the help of MATHEMATICA\ program for some values of parameters $%
\alpha $ and $\beta $ are derived. There are many other expressions for
these functions, which are available using this program, but being long and
complex, they were omitted. The integral Mainardi fuctions $\mathrm{Fi}%
_{\alpha }\left( x\right) $ and $\mathrm{Mi}_{\alpha }\left( x\right) $ for $%
0<\alpha <1$, are presented in Tables \ref{Table8} and \ref{Table9}. As can
be expected, most of these integral functions are expressed in terms of
generalized hypergeometric functions.

\begin{specialtable}[H] \centering%
\caption{The integral Wright functions $\mathrm{Wi}_{\alpha , \beta}$ derived for some values of
parameters $\alpha$ and $\beta$ by using
(\ref{Phi_i_series}).}%
\begin{footnotesize}
\begin{tabular*}{\hsize}{c@{\extracolsep{\fill}}cc}
\toprule
\boldmath$\alpha $ & \boldmath$\beta $ & \boldmath$\mathrm{Wi}_{\alpha ,\beta }\left( x\right) $ \\
\midrule
$-1$ & $\frac{1}{2}$ & $\frac{1}{\sqrt{\pi }}\left[ \mathrm{\ln }4-2\mathrm{%
\ln }\left( \sqrt{1+x}+1\right) \right] $ \\ \midrule
$-1$ & $\frac{3}{2}$ & $\frac{x}{\sqrt{\pi }}\,_{3}F_{2}\left( \left.
\begin{array}{c}
\frac{1}{2},1,1 \\
2,2%
\end{array}%
\right\vert -x\right) $ \\ \midrule
$-1$ & $\beta $ & $\frac{x}{\Gamma \left( \beta -1\right) }\,_{3}F_{2}\left(
\left.
\begin{array}{c}
1,1,2-\beta \\
2,2%
\end{array}%
\right\vert -x\right) $ \\ \midrule
$0$ & $-\frac{4}{3}$ & $\frac{-\gamma -\mathrm{\ln }x+\mathrm{Chi}\left(
x\right) +\mathrm{Shi}\left( x\right) }{\Gamma \left( -4/3\right) }$ \\
\midrule
$0$ & $\beta $ & $\frac{-\gamma -\mathrm{\ln }x+\mathrm{Chi}\left( x\right) +%
\mathrm{Shi}\left( x\right) }{\Gamma \left( \beta \right) }$ \\ \midrule
$\frac{1}{2}$ & $0$ & $-\frac{1}{2}+\frac{1}{2}\,_{0}F_{2}\left( \left.
\begin{array}{c}
- \\
\frac{1}{2},1%
\end{array}%
\right\vert \frac{x^{2}}{4}\right) +\frac{x}{\sqrt{\pi }}\,_{0}F_{2}\left(
\left.
\begin{array}{c}
- \\
\frac{1}{2},\frac{3}{2}%
\end{array}%
\right\vert \frac{x^{2}}{4}\right) $ \\

\bottomrule
\end{tabular*}
\end{footnotesize}
\end{specialtable}

\begin{specialtable}[H]\ContinuedFloat
\small
\caption{{\em Cont.}}
\begin{footnotesize}
\begin{tabular*}{\hsize}{c@{\extracolsep{\fill}}cc}
\toprule
\boldmath$\alpha $ & \boldmath$\beta $ & \boldmath$\mathrm{Wi}_{\alpha ,\beta }\left( x\right) $ \\
\midrule

$\frac{1}{2}$ & $\frac{1}{2}$ & $\frac{x^{2}}{2\sqrt{\pi }}\,_{2}F_{4}\left(
\left.
\begin{array}{c}
1,1 \\
\frac{3}{2},\frac{3}{2},2,2%
\end{array}%
\right\vert \frac{x^{2}}{4}\right) +x\,_{1}F_{3}\left( \left.
\begin{array}{c}
\frac{1}{2} \\
1,\frac{3}{2},\frac{3}{2}%
\end{array}%
\right\vert \frac{x^{2}}{4}\right) $ \\ \midrule
$\frac{1}{2}$ & $1$ & $\frac{x}{4}\left[ \frac{8}{\sqrt{\pi }}%
\,_{1}F_{3}\left( \left.
\begin{array}{c}
\frac{1}{2} \\
\frac{3}{2},\frac{3}{2},\frac{3}{2}%
\end{array}%
\right\vert \frac{x^{2}}{4}\right) +x\,_{2}F_{4}\left( \left.
\begin{array}{c}
1,1 \\
\frac{3}{2},\frac{3}{2},2,2%
\end{array}%
\right\vert \frac{x^{2}}{4}\right) \right] $ \\ \midrule
$\frac{1}{2}$ & $2$ & $\frac{x^{2}}{8}\,_{2}F_{4}\left( \left.
\begin{array}{c}
1,1 \\
\frac{3}{2},2,2,3%
\end{array}%
\right\vert \frac{x^{2}}{4}\right) +\frac{4x}{3\sqrt{\pi }}\,_{1}F_{3}\left(
\left.
\begin{array}{c}
\frac{1}{2} \\
\frac{3}{2},\frac{3}{2},\frac{5}{2}%
\end{array}%
\right\vert \frac{x^{2}}{4}\right) $ \\ \midrule
$\frac{1}{2}$ & $\beta $ & $\frac{x^{2}}{4\,\Gamma \left( 1+\beta \right) }%
\,_{2}F_{4}\left( \left.
\begin{array}{c}
1,1 \\
\frac{3}{2},2,2,\beta +1%
\end{array}%
\right\vert \frac{x^{2}}{4}\right) +\frac{x}{\Gamma \left( \frac{1}{2}+\beta
\right) }\,_{1}F_{3}\left( \left.
\begin{array}{c}
\frac{1}{2} \\
\frac{3}{2},\frac{3}{2},\beta +\frac{1}{2}%
\end{array}%
\right\vert \frac{x^{2}}{4}\right) $ \\ \midrule
$1$ & $-\frac{3}{2}$ & $-\frac{x}{2\sqrt{\pi }}\,_{2}F_{3}\left( \left.
\begin{array}{c}
1,1 \\
-\frac{1}{2},2,2%
\end{array}%
\right\vert x\right) $ \\ \midrule
$3$ & $\frac{3}{2}$ & $\frac{x}{\sqrt{\pi }}\,_{2}F_{3}\left( \left.
\begin{array}{c}
1,1 \\
\frac{1}{2},2,2%
\end{array}%
\right\vert x\right) $ \\ \midrule
$1$ & $0$ & $-1+I_{0}\left( 2\sqrt{x}\right) $ \\ \midrule
$1$ & $\frac{1}{4}$ & $\frac{x}{\Gamma \left( 5/4\right) }\,_{2}F_{3}\left(
\left.
\begin{array}{c}
1,1 \\
\frac{5}{2},2,2%
\end{array}%
\right\vert x\right) $ \\ \midrule
$1$ & $\frac{1}{2}$ & $-\frac{2\gamma +\mathrm{\ln }4+\mathrm{\ln }x-2\,%
\mathrm{Chi}\left( 2\sqrt{x}\right) }{\sqrt{\pi }}$ \\ \midrule
$1$ & $1$ & $x\,_{2}F_{3}\left( \left.
\begin{array}{c}
1,1 \\
2,2,2%
\end{array}%
\right\vert x\right) $ \\ \midrule
$1$ & $\frac{3}{2}$ & $-\frac{1}{\sqrt{\pi x}}\left[ 2\sinh \left( 2\sqrt{x}%
\right) -2\sqrt{x}\left( 2\gamma -2+\mathrm{\ln }4-2\,\mathrm{Chi}\left( 2%
\sqrt{x}\right) \right) \right] $ \\ \midrule
$1$ & $\beta $ & $x\,_{2}F_{3}\left( \left.
\begin{array}{c}
1,1 \\
2,2,\beta +1%
\end{array}%
\right\vert x\right) $ \\ \bottomrule
\end{tabular*}%
\end{footnotesize}
\label{Table7}%
\end{specialtable}%
\vspace{-12pt}
\begin{specialtable}[H] \centering%
\caption{The integral Wright functions $\mathrm{Wi}_{\alpha , \beta}$ derived for some values of
parameters $\alpha$ and $\beta$ by using
(\ref{Phi_i_series}).}%
\begin{footnotesize}
\begin{tabular*}{\hsize}{c@{\extracolsep{\fill}}cc}
\toprule
\boldmath$\alpha $ & \boldmath$\beta $ & \boldmath$\mathrm{Wi}_{\alpha ,\beta }\left( x\right) $ \\
\midrule
$\frac{3}{2}$ & $\frac{1}{2}$ & $\frac{2x^{2}}{15\sqrt{\pi }}%
\,_{2}F_{6}\left( \left.
\begin{array}{c}
1,1 \\
\frac{7}{2},\frac{3}{2},\frac{3}{2},\frac{11}{6},2,2%
\end{array}%
\right\vert \frac{x^{2}}{108}\right) +x\,_{1}F_{5}\left( \left.
\begin{array}{c}
\frac{1}{2} \\
\frac{2}{3},1,\frac{4}{3},\frac{3}{2},\frac{3}{2}%
\end{array}%
\right\vert \frac{x^{2}}{108}\right) $ \\ \midrule
$2$ & $\frac{1}{4}$ & $\frac{16\,x}{5\,\Gamma \left( 1/4\right) }%
\,_{2}F_{4}\left( \left.
\begin{array}{c}
1,1 \\
\frac{9}{8},\frac{13}{8},2,2%
\end{array}%
\right\vert \frac{x}{4}\right) $ \\ \midrule
$2$ & $\frac{1}{3}$ & $\frac{9\,x}{4\,\Gamma \left( 1/3\right) }%
\,_{2}F_{4}\left( \left.
\begin{array}{c}
1,1 \\
\frac{7}{6},\frac{5}{3},2,2%
\end{array}%
\right\vert \frac{x}{4}\right) $ \\ \midrule
$2$ & $\frac{1}{2}$ & $\frac{4\,x}{3\,\sqrt{\pi }}\,_{2}F_{4}\left( \left.
\begin{array}{c}
1,1 \\
\frac{5}{4},\frac{7}{4},2,2%
\end{array}%
\right\vert \frac{x}{4}\right) $ \\ \midrule
$2$ & $1$ & $\frac{\,x}{2}\,_{2}F_{4}\left( \left.
\begin{array}{c}
1,1 \\
\frac{3}{2},2,2,2%
\end{array}%
\right\vert \frac{x}{4}\right) $ \\ \midrule
$2$ & $2$ & $\frac{\,x}{6}\,_{2}F_{4}\left( \left.
\begin{array}{c}
1,1 \\
\frac{5}{2},2,2,2%
\end{array}%
\right\vert \frac{x}{4}\right) $ \\ \midrule
$2$ & $\beta $ & $\frac{\,x}{\Gamma \left( \beta +2\right) }%
\,_{2}F_{4}\left( \left.
\begin{array}{c}
1,1 \\
2,2,\frac{\beta }{2}+1,\frac{\beta +3}{2}%
\end{array}%
\right\vert \frac{x}{4}\right) $ \\ \midrule
$3$ & $1$ & $\frac{\,x}{6}\,_{2}F_{4}\left( \left.
\begin{array}{c}
1,1 \\
\frac{4}{3},\frac{5}{3},2,2%
\end{array}%
\right\vert \frac{x}{27}\right) $ \\ \midrule
$3$ & $\beta $ & $\frac{\,x}{\Gamma \left( \beta +3\right) }%
\,_{2}F_{5}\left( \left.
\begin{array}{c}
1,1 \\
2,2,\frac{\beta }{3}+1,\frac{\beta +4}{3},\frac{\beta +5}{3}%
\end{array}%
\right\vert \frac{x}{27}\right) $ \\ \midrule
$4$ & $\beta $ & $\frac{\,x}{\Gamma \left( \beta +4\right) }%
\,_{2}F_{6}\left( \left.
\begin{array}{c}
1,1 \\
2,2,\frac{\beta }{4}+1,\frac{\beta +5}{4},\frac{\beta +6}{4},\frac{\beta +7}{%
4}%
\end{array}%
\right\vert \frac{x}{256}\right) $ \\ \midrule
$5$ & $\beta $ & $\frac{\,x}{\Gamma \left( \beta +5\right) }%
\,_{2}F_{7}\left( \left.
\begin{array}{c}
1,1 \\
2,2,\frac{\beta }{4}+1,\frac{\beta +6}{5},\frac{\beta +7}{5},\frac{\beta +8}{%
5},\frac{\beta +9}{5}%
\end{array}%
\right\vert \frac{x}{3125}\right) $ \\ \bottomrule
\end{tabular*}%
\end{footnotesize}
\label{Table7a}%
\end{specialtable}%

\begin{specialtable}[H] \centering%
\caption{The integral Mainardi function $\mathrm{Fi}_{\alpha}$ derived for some values of
parameter $\alpha$ by using (\ref{Integral_Mainardi_p/q}).}%
\begin{tabular*}{\hsize}{c@{\extracolsep{\fill}}c}
\toprule
\boldmath$\alpha $ & \boldmath$\mathrm{Fi}_{\alpha }\left( x\right) $ \\ \midrule
$\frac{3}{4}$ & $%
\begin{array}{l}
-x\left[ \frac{1}{\Gamma \left( -\frac{3}{4}\right) }\,_{3}F_{3}\left(
\left.
\begin{array}{c}
\frac{1}{4},\frac{7}{12},\frac{11}{12} \\
\frac{1}{2},\frac{3}{4},\frac{5}{4}%
\end{array}%
\right\vert -\frac{27x^{4}}{256}\right) +\frac{x}{144}\right. \\
\left. \left( \frac{8x}{\Gamma \left( -\frac{9}{4}\right) }%
\,\,_{3}F_{3}\left( \left.
\begin{array}{c}
\frac{3}{4},\frac{13}{12},\frac{17}{12} \\
\frac{5}{4},\frac{3}{2},\frac{7}{4}%
\end{array}%
\right\vert -\frac{27x^{4}}{256}\right) -\frac{27}{\sqrt{\pi }}%
\,_{3}F_{4}\left( \left.
\begin{array}{c}
\frac{1}{2},\frac{5}{6},\frac{7}{6} \\
\frac{3}{4},\frac{5}{4},\frac{3}{2}%
\end{array}%
\right\vert -\frac{27x^{4}}{256}\right) \right) \right]%
\end{array}%
$ \\ \midrule
$\frac{2}{3}$ & $\frac{x}{4}\left[ \frac{x}{\Gamma \left( -\frac{4}{3}%
\right) }\,\,_{2}F_{2}\left( \left.
\begin{array}{c}
\frac{2}{3},\frac{7}{6} \\
\frac{4}{3},\frac{5}{3}%
\end{array}%
\right\vert -\frac{4x^{3}}{27}\right) -\frac{4}{\Gamma \left( -\frac{2}{3}%
\right) }\,_{2}F_{2}\left( \left.
\begin{array}{c}
\frac{1}{3},\frac{5}{6} \\
\frac{2}{3},\frac{4}{3}%
\end{array}%
\right\vert -\frac{4x^{3}}{27}\right) \right] $ \\ \midrule
$\frac{1}{2}$ & $\frac{1}{2}\,\mathrm{erf}\left( \frac{x}{2}\right) $ \\
\midrule
$\frac{1}{3}$ & $\frac{x}{4}\left[ \frac{x}{\Gamma \left( -\frac{2}{3}%
\right) }\,\,_{1}F_{2}\left( \left.
\begin{array}{c}
\frac{2}{3} \\
\frac{4}{3},\frac{5}{3}%
\end{array}%
\right\vert \frac{x^{3}}{27}\right) -\frac{4}{\Gamma \left( -\frac{1}{3}%
\right) }\,_{1}F_{2}\left( \left.
\begin{array}{c}
\frac{1}{3} \\
\frac{2}{3},\frac{4}{3}%
\end{array}%
\right\vert \frac{x^{3}}{27}\right) \right] $ \\ \midrule
$\frac{1}{4}$ & $%
\begin{array}{l}
-x\left[ \frac{1}{\Gamma \left( -\frac{1}{4}\right) }\,_{1}F_{3}\left(
\left.
\begin{array}{c}
\frac{1}{4} \\
\frac{1}{2},\frac{3}{4},\frac{5}{4}%
\end{array}%
\right\vert -\frac{x^{4}}{256}\right) +\frac{x}{72}\right. \\
\left. \left( \frac{9}{\sqrt{\pi }}\,\,_{1}F_{3}\left( \left.
\begin{array}{c}
\frac{1}{2} \\
\frac{3}{4},\frac{5}{4},\frac{3}{2}%
\end{array}%
\right\vert -\frac{x^{4}}{256}\right) +\frac{4x}{\Gamma \left( -\frac{3}{4}%
\right) }\,_{3}F_{4}\left( \left.
\begin{array}{c}
\frac{3}{4} \\
\frac{5}{4},\frac{3}{2},\frac{7}{4}%
\end{array}%
\right\vert -\frac{x^{4}}{256}\right) \right) \right]%
\end{array}%
$ \\ \bottomrule
\end{tabular*}%
\label{Table8}%
\end{specialtable}%

\vspace{-12pt}
\begin{specialtable}[H] \centering%
\caption{The integral Mainardi function $\mathrm{Mi}_{\alpha}$ derived for some values of
parameter $\alpha$ by using (\ref{Integral_Mainardi_p/q}).}%
\begin{tabular*}{\hsize}{c@{\extracolsep{\fill}}c}
\toprule
\boldmath$\alpha $ & \boldmath$\mathrm{Mi}_{\alpha }\left( x\right) $ \\ \midrule
$\frac{3}{4}$ & $%
\begin{array}{l}
\frac{x}{96}\left[ \frac{48}{\sqrt{\pi }}\,_{3}F_{3}\left( \left.
\begin{array}{c}
\frac{1}{4},\frac{5}{6},\frac{7}{6} \\
\frac{3}{4},\frac{5}{4},\frac{5}{4}%
\end{array}%
\right\vert -\frac{27x^{4}}{256}\right) \right. \\
\left. +\frac{24x}{\Gamma \left( -\frac{5}{4}\right) }\,\,_{3}F_{3}\left(
\left.
\begin{array}{c}
\frac{1}{2},\frac{13}{12},\frac{17}{12} \\
\frac{5}{4},\frac{3}{2},\frac{3}{2}%
\end{array}%
\right\vert -\frac{27x^{4}}{256}\right) +\frac{x^{3}}{\Gamma \left( -\frac{11%
}{4}\right) }\,_{4}F_{4}\left( \left.
\begin{array}{c}
1,1,\frac{19}{12},\frac{23}{12} \\
\frac{3}{2},\frac{7}{4},2,2%
\end{array}%
\right\vert -\frac{27x^{4}}{256}\right) \right]%
\end{array}%
$ \\ \midrule
$\frac{2}{3}$ & $-\frac{x^{3}}{18\,\Gamma \left( -\frac{5}{3}\right) }%
\,\,_{3}F_{3}\left( \left.
\begin{array}{c}
1,1,\frac{11}{6} \\
\frac{5}{3},2,2%
\end{array}%
\right\vert -\frac{4x^{3}}{27}\right) -\frac{x}{\Gamma \left( -\frac{1}{3}%
\right) }\,_{2}F_{2}\left( \left.
\begin{array}{c}
\frac{1}{3},\frac{7}{6} \\
\frac{4}{3},\frac{4}{3}%
\end{array}%
\right\vert -\frac{4x^{3}}{27}\right) $ \\ \midrule
$\frac{1}{2}$ & $\frac{1}{2\sqrt{\pi }}\left[ \mathrm{Chi}\left( \frac{x^{2}%
}{4}\right) -\mathrm{Shi}\left( \frac{x^{2}}{4}\right) -\mathrm{\ln }\left(
\frac{x^{2}}{4}\right) -\gamma \right] $ \\ \midrule
$\frac{1}{3}$ & $-\frac{x^{3}}{18\,\Gamma \left( -\frac{1}{3}\right) }%
\,\,_{2}F_{3}\left( \left.
\begin{array}{c}
1,1 \\
\frac{5}{3},2,2%
\end{array}%
\right\vert \frac{x^{3}}{27}\right) -\frac{x}{\Gamma \left( \frac{1}{3}%
\right) }\,_{1}F_{2}\left( \left.
\begin{array}{c}
\frac{1}{3} \\
\frac{4}{3},\frac{4}{3}%
\end{array}%
\right\vert \frac{x^{3}}{27}\right) $ \\ \midrule
$\frac{1}{4}$ & $%
\begin{array}{l}
-\frac{x}{\sqrt{\pi }}\,_{1}F_{3}\left( \left.
\begin{array}{c}
\frac{1}{4} \\
\frac{3}{4},\frac{5}{4},\frac{5}{4}%
\end{array}%
\right\vert -\frac{x^{4}}{256}\right) \\
\left. +\frac{x^{2}}{4\,\Gamma \left( \frac{1}{4}\right) }\,_{1}F_{3}\left(
\left.
\begin{array}{c}
\frac{1}{2} \\
\frac{5}{4},\frac{3}{2},\frac{3}{2}%
\end{array}%
\right\vert -\frac{x^{4}}{256}\right) +\frac{x^{4}}{96\,\Gamma \left( -\frac{%
1}{4}\right) }\,_{2}F_{4}\left( \left.
\begin{array}{c}
1,1 \\
\frac{3}{2},\frac{7}{4},2,2%
\end{array}%
\right\vert -\frac{x^{4}}{256}\right) \right]%
\end{array}%
$ \\ \bottomrule
\end{tabular*}%
\label{Table9}%
\end{specialtable}%


\section{Conclusions} \label{Section_Conclusions}

For the first time, three new special functions are presented in this investigation:
the integral Mittag-Leffler functions, the integral
Whittaker functions, and the integral Wright functions. These functions are
defined in the mathematical literature in the same manner as other elementary and special integral functions.
It is feasible to generate these functions in an explicit form for certain parameters values using the MATHEMATICA application.
These integral functions are often represented in terms of
generalized hypergeometric functions. The behavior of some of them is
shown graphically. In the Appendices, a large number of
Mittag-Leffler, Whittaker, and Wright functions with integral and fractional
parameters, as well as their Laplace transforms, are presented in tabular form.

It may be observed that, generally, it is highly possible to make general integral functions such as (\ref{Integral_Whittaker_M_def}) and (\ref{Integral_Whittaker_W_def}) by using generalized hypergeometric $_{p}F_{p}(t)$, because they converge in the whole complex t-plane, or, for every real number $t$.

\vspace{6pt}
\authorcontributions{Conceptualization, A.A. and J.L.G-S.; Methodology, A.A. and J.L.G-S.; Resources, A.A.; Writing---original draft, A.A. and J.L.G-S.; Writing---review \& editing, A.A. and J.L.G-S. All authors have read and agreed to the published version of the manuscript.}

\funding{This research received no external funding.}
\institutionalreview{Not applicable.}

\informedconsent{Not applicable.}

\dataavailability{Not applicable.}


\acknowledgments{We are grateful to Armando Consiglio from Institut f\"{u}r Theoretische
Physik und Astrophysik of University of W\"{u}rzburg, Germany, who performed
numerical evaluations of the integral functions presented in Figures \ref%
{Figure: Mittag-Leffler a} and \ref{Figure: Mittag-Leffler b}, and to
Professor Francesco Mainardi from the Department of Physics and Astronomy,
University of Bologna, Bologna, Italy, for his kind encouragement and
interest in our work.}

\conflictsofinterest{The authors declare no conflict of interest.} 

\appendixtitles{yes} 
\appendixstart
\appendix

\section{Representations of the One- and Two-Parameter Mittag-Leffler
Functions and Their Laplace Transforms \label{Appendix: Mittag-Leffler}}

The Mittag-Leffler functions are defined by the sums of infinite series
presented in (\ref{Mittag_Leffler_def}) and their Laplace transforms in (\ref%
{LT_Mittag-Leffler}). For positive variable $x$ and some values of
parameters $\alpha $ and $\beta $, these sums can be expressed in terms of
elementary and special functions, especially in terms of generalized
hypergeometric functions. They were derived by using the MATHEMATICA program and
presented in Tables \ref{TableA} and \ref{TableAa} for the Mittag-Leffler
functions, as well as Tables \ref{TableB} and \ref{TableBa} for the Laplace
transforms. These results, given in terms of infinite series, are mostly new,
and only they are only partly  known in the mathematical literature. Knowing that
any infinite sum can be split as\
\begin{equation}
\sum_{k=0}^{\infty }a\left( k\right) =\sum_{j=0}^{q-1}\sum_{k=0}^{\infty
}a\left( qk+j\right) ,  \label{Sum_split}
\end{equation}%
and applying the multiplication formula of the gamma function (\cite[]{DLMF} (Eqn.
5.5.6)), for $nt\neq 0,-1,-2,\ldots $%
\begin{equation}
\Gamma \left( nt\right) =\left( 2\pi \right) ^{\left( 1-n\right)
/2}n^{nt-1/2}\prod\limits_{j=0}^{n-1}\Gamma \left( t+\frac{j}{n}\right) ,
\label{Multiplication_Gamma}
\end{equation}%
it is possible to express from (\ref{Mittag_Leffler_def}) the Mittag-Leffler function in the case of
positive rational $\alpha =p/q$ with $p$ and $q$ positive coprimes,
\begin{equation}
\mathrm{E}_{p/q,\beta }\left( x\right) =\sum_{k=0}^{q-1}\frac{x^{k}}{\Gamma
\left( \frac{p}{q}k+\beta \right) }\,_{1}F_{p}\left( \left.
\begin{array}{c}
1 \\
b_{0},\ldots ,b_{p-1}%
\end{array}%
\right\vert \frac{x^{q}}{p^{p}}\right) ,
\label{Mittag-Leffler_p/q_resultado}
\end{equation}%
where%
\[
b_{j}=\frac{k}{q}+\frac{\beta +j}{p}.
\]

The corresponding Laplace transforms are%
\begin{eqnarray}
&&\mathcal{L}\left[ \mathrm{E}_{p/q,\beta }\left( t\right) \right]
\label{Laplace_Mittag-Leffler_p/q_resultado} \\
&=&\sum_{k=0}^{q-1}\frac{s^{-k-1}}{\Gamma \left( \frac{p}{q}k+\beta \right) }%
\,_{q+1}F_{p}\left( \left.
\begin{array}{c}
1,a_{0},\ldots ,a_{q-1} \\
b_{0},\ldots ,b_{p-1}%
\end{array}%
\right\vert \frac{\left( q/s\right) ^{q}}{p^{p}}\right) ,  \nonumber
\end{eqnarray}%
where%
\begin{eqnarray*}
a_{j} &=&\frac{k+1+j}{q}, \\
b_{j} &=&\frac{k}{q}+\frac{\beta +j}{p}.
\end{eqnarray*}

\begin{specialtable}[H] \centering%
\caption{The Mittag-Leffler functions derived for some values of
parameters $\alpha$ and $\beta$ by using (\ref{Mittag_Leffler_def}).}%
\begin{tabular*}{\hsize}{c@{\extracolsep{\fill}}cc}
\toprule
\boldmath$\alpha $ & \boldmath$\beta $ & \boldmath$\mathrm{E}_{\alpha ,\beta }\left( x\right) $ \\
\midrule
$\frac{1}{2}$ & $\frac{1}{2}$ & $\frac{1}{\sqrt{\pi }}+x\,e^{x^{2}}\left[
\mathrm{erf}\left( x\right) +1\right] $ \\ \midrule
$\frac{1}{2}$ & $1$ & $e^{x^{2}}\left[ \mathrm{erf}\left( x\right) +1\right]
$ \\ \midrule
$\frac{1}{2}$ & $\frac{3}{2}$ & $\frac{e^{x^{2}}\left[ \mathrm{erf}\left(
x\right) +1\right] -1}{x}$ \\ \midrule
$\frac{1}{2}$ & $2$ & $\frac{1}{x^{2}}\left[ e^{x^{2}}\left[ \mathrm{erf}%
\left( x\right) +1\right] -1-\frac{2x}{\sqrt{\pi }}\right] $ \\ \midrule
$\frac{1}{2}$ & $3$ & $-\frac{1}{3x^{4}}\left[ 3e^{x^{2}}\left[ \mathrm{erfc}%
\left( x\right) -2\right] +\frac{4x^{3}+6x}{\sqrt{\pi }}+3\left(
1+x^{2}\right) \right] $ \\ \midrule
$\frac{1}{2}$ & $4$ & $\frac{1}{30x^{4}}\left[ 30\,e^{x^{2}}\left[ \mathrm{%
erf}\left( x\right) +1\right] -\frac{4x\left( 4x^{4}+10x^{2}+15\right) }{%
\sqrt{\pi }}-15\left( x^{4}+2x^{2}+2\right) \right] $ \\ \midrule
$\frac{1}{2}$ & $\beta $ & $e^{x^{2}}x^{2\left( 1-\beta \right) }\left[ 2-%
\frac{\Gamma \left( \beta -1,x^{2}\right) }{\Gamma \left( \beta -1\right) }-%
\frac{\Gamma \left( \beta -\frac{1}{2},x^{2}\right) }{\Gamma \left( \beta -%
\frac{1}{2}\right) }\right] $ \\ \midrule
$1$ & $\frac{1}{2}$ & $\frac{1}{\sqrt{\pi }}+\sqrt{x}\,e^{x}\,\mathrm{erf}%
\left( \sqrt{x}\right) $ \\ \midrule
$1$ & $1$ & $e^{x}$ \\ \midrule
$1$ & $\frac{3}{2}$ & $\frac{e^{x}\,\mathrm{erf}\left( \sqrt{x}\right) }{%
\sqrt{x}}$ \\ \midrule
$1$ & $2$ & $\frac{e^{x}-1}{x}$ \\ \midrule
$\frac{3}{2}$ & $\frac{1}{2}$ & $\frac{1}{\sqrt{\pi }}\,_{1}F_{3}\left(
\left.
\begin{array}{c}
1 \\
\frac{1}{6},\frac{1}{2},\frac{5}{6}%
\end{array}%
\right\vert \frac{x^{2}}{27}\right) +\frac{x^{1/3}}{3}\left[
e^{x^{2/3}}-2\,e^{-x^{2/3}/2}\sin \left( \frac{\pi -3\sqrt{3}x^{2/3}}{6}%
\right) \right] $ \\ \midrule
$\frac{3}{2}$ & $1$ & $\frac{1}{3}\left[ \frac{4x}{\sqrt{\pi }}%
\,_{1}F_{3}\left( \left.
\begin{array}{c}
1 \\
\frac{5}{6},\frac{7}{6},\frac{3}{2}%
\end{array}%
\right\vert \frac{x^{2}}{27}\right) +e^{x^{2/3}}+2\,e^{-x^{2/3}/2}\cos
\left( \frac{\sqrt{3}}{2}x^{2/3}\right) \right] $ \\ \midrule
$\frac{3}{2}$ & $\frac{3}{2}$ & $\frac{2}{\sqrt{\pi }}\,_{1}F_{3}\left(
\left.
\begin{array}{c}
1 \\
\frac{1}{2},\frac{5}{6},\frac{7}{6}%
\end{array}%
\right\vert \frac{x^{2}}{27}\right) +\frac{x^{-1/3}}{3}\left[
e^{x^{2/3}}-2\,e^{-x^{2/3}/2}\sin \left( \frac{\pi +3\sqrt{3}x^{2/3}}{6}%
\right) \right] $ \\ \midrule
$\frac{3}{2}$ & $2$ & $\frac{8x}{15\sqrt{\pi }}\,_{1}F_{3}\left( \left.
\begin{array}{c}
1 \\
\frac{7}{6},\frac{3}{2},\frac{11}{6}%
\end{array}%
\right\vert \frac{x^{2}}{27}\right) +\frac{x^{-2/3}}{3}\left[
e^{x^{2/3}}-2\,e^{-x^{2/3}/2}\sin \left( \frac{\pi -3\sqrt{3}x^{2/3}}{6}%
\right) \right] $ \\ \midrule
$\frac{3}{2}$ & $\beta $ & $\frac{x}{\Gamma \left( \beta +\frac{1}{2}\right)
}\,_{1}F_{3}\left( \left.
\begin{array}{c}
1 \\
\frac{2\beta +3}{6},\frac{2\beta +5}{6},\frac{2\beta +7}{6}%
\end{array}%
\right\vert \frac{x^{2}}{27}\right) +\frac{1}{\Gamma \left( \beta \right) }%
\,_{1}F_{3}\left( \left.
\begin{array}{c}
1 \\
\frac{\beta +1}{3},\frac{\beta +2}{3},\frac{\beta }{3}%
\end{array}%
\right\vert \frac{x^{2}}{27}\right) $ \\ \midrule
$2$ & $\frac{1}{2}$ & $\frac{1}{\sqrt{\pi }}\,_{1}F_{2}\left( \left.
\begin{array}{c}
1 \\
\frac{1}{4},\frac{3}{4}%
\end{array}%
\right\vert \frac{x}{4}\right) $ \\ \midrule
$2$ & $1$ & $\cosh \left( \sqrt{x}\right) $ \\ \midrule
$2$ & $2$ & $\frac{\sinh \left( \sqrt{x}\right) }{\sqrt{x}}$ \\ \midrule
$2$ & $3$ & $\frac{\cosh \left( \sqrt{x}\right) -1}{x}$ \\ \midrule
$2$ & $4$ & $\frac{\sinh \left( \sqrt{x}\right) }{x^{3/2}}-\frac{1}{x}$ \\
\midrule
$2$ & $\beta $ & $\frac{1}{\Gamma \left( \beta \right) }\,_{1}F_{2}\left(
\left.
\begin{array}{c}
1 \\
\frac{\beta +1}{2},\frac{\beta }{2}%
\end{array}%
\right\vert \frac{x}{4}\right) $ \\ \bottomrule
\end{tabular*}%
\label{TableA}%
\end{specialtable}%

\begin{specialtable}[H] \centering%
\caption{The Mittag-Leffler functions derived for some values of
parameters $\alpha$ and $\beta$ by using (\ref{Mittag_Leffler_def}).}%
\begin{tabular*}{\hsize}{c@{\extracolsep{\fill}}cc}
\toprule
\boldmath$\alpha $ & \boldmath$\beta $ & \boldmath$\mathrm{E}_{\alpha ,\beta }\left( x\right) $ \\
\midrule
$3$ & $1$ & $\frac{1}{3}\left[ e^{x^{1/3}}+2\,e^{-x^{1/3}/2}\cos \left(
\frac{\sqrt{3}}{2}x^{1/3}\right) \right] $ \\ \midrule
$3$ & $2$ & $\frac{x^{-1/3}}{3}\left[ e^{x^{1/3}}-2\,e^{-x^{1/3}/2}\sin
\left( \frac{\pi -3\sqrt{3}x^{1/3}}{6}\right) \right] $ \\ \midrule
$3$ & $3$ & $\frac{x^{-2/3}}{3}\left[ e^{x^{1/3}}-2\,e^{-x^{1/3}/2}\sin
\left( \frac{\pi +3\sqrt{3}x^{1/3}}{6}\right) \right] $ \\ \midrule
$3$ & $\beta $ & $\frac{1}{\Gamma \left( \beta \right) }\,_{1}F_{3}\left(
\left.
\begin{array}{c}
1 \\
\frac{\beta +1}{3},\frac{\beta +2}{3},\frac{\beta }{3}%
\end{array}%
\right\vert \frac{x}{27}\right) $ \\ \midrule
$4$ & $1$ & $\frac{1}{2}\left[ \cos \left( x^{1/4}\right) +\cosh \left(
x^{1/4}\right) \right] $ \\ \midrule
$4$ & $2$ & $\frac{\sin \left( x^{1/4}\right) +\sinh \left( x^{1/4}\right) }{%
2x^{1/4}}$ \\ \midrule
$4$ & $3$ & $\frac{\cosh \left( x^{1/4}\right) -\cos \left( x^{1/4}\right) }{%
2\sqrt{x}}$ \\ \midrule
$4$ & $4$ & $\frac{\sinh \left( x^{1/4}\right) -\sin \left( x^{1/4}\right) }{%
2x^{3/4}}$ \\ \midrule
$4$ & $\beta $ & $\frac{1}{\Gamma \left( \beta \right) }\,_{1}F_{4}\left(
\left.
\begin{array}{c}
1 \\
\frac{\beta +1}{4},\frac{\beta +2}{4},\frac{\beta +3}{4},\frac{\beta }{4}%
\end{array}%
\right\vert \frac{x}{256}\right) $ \\ \midrule
$5$ & $1$ & $\,_{0}F_{4}\left( \left.
\begin{array}{c}
- \\
\frac{1}{5},\frac{2}{5},\frac{3}{5},\frac{4}{5}%
\end{array}%
\right\vert \frac{x}{3125}\right) $ \\ \midrule
$5$ & $2$ & $_{0}F_{4}\left( \left.
\begin{array}{c}
- \\
\frac{2}{5},\frac{3}{5},\frac{4}{5},\frac{6}{5}%
\end{array}%
\right\vert \frac{x}{3125}\right) $ \\ \midrule
$5$ & $3$ & $\frac{1}{2}\,_{0}F_{4}\left( \left.
\begin{array}{c}
- \\
\frac{3}{5},\frac{4}{5},\frac{6}{5},\frac{7}{5}%
\end{array}%
\right\vert \frac{x}{3125}\right) $ \\ \midrule
$5$ & $4$ & $\frac{1}{6}\,_{0}F_{4}\left( \left.
\begin{array}{c}
- \\
\frac{4}{5},\frac{6}{5},\frac{7}{5},\frac{8}{5}%
\end{array}%
\right\vert \frac{x}{3125}\right) $ \\ \midrule
$5$ & $5$ & $\frac{1}{24}\,_{0}F_{4}\left( \left.
\begin{array}{c}
- \\
\frac{6}{5},\frac{7}{5},\frac{8}{5},\frac{9}{5}%
\end{array}%
\right\vert \frac{x}{3125}\right) $ \\ \midrule
$5$ & $\beta $ & $\frac{1}{\Gamma \left( \beta \right) }\,_{1}F_{5}\left(
\left.
\begin{array}{c}
1 \\
\frac{\beta +1}{5},\frac{\beta +2}{5},\frac{\beta +3}{5},\frac{\beta +4}{5},%
\frac{\beta }{5}%
\end{array}%
\right\vert \frac{x}{3125}\right) $ \\ \bottomrule
\end{tabular*}%
\label{TableAa}%
\end{specialtable}%
\vspace{-6pt}
\begin{specialtable}[H] \centering%
\caption{The Laplace transforms Mittag-Leffler functions derived for some values of
parameters $\alpha$ and $\beta$ by using (\ref{LT_Mittag-Leffler}).}%
\begin{tabular*}{\hsize}{c@{\extracolsep{\fill}}cc}
\toprule
$\alpha $ & $\beta $ & $\quad \mathcal{L}\left[ \mathrm{E}_{\alpha ,\beta
}\left( t\right) \right] $ \\ \midrule
$1$ & $\frac{1}{2}$ & $\frac{\sqrt{s-1}+\csc ^{-1}\left( \sqrt{s}\right) }{%
\sqrt{\pi }\left( s-1\right) ^{3/2}}$ \\ \midrule
$\frac{1}{2}$ & $1$ & $\frac{1}{s-1}$ \\ \midrule
$1$ & $\frac{3}{2}$ & $\frac{2\csc ^{-1}\left( \sqrt{s}\right) }{\sqrt{\pi }%
\sqrt{s-1}}$ \\ \midrule
$1$ & $2$ & $\mathrm{\ln }\left( \frac{s}{s-1}\right) $ \\ \midrule
$1$ & $\beta $ & $\frac{1}{s\,\Gamma \left( \beta \right) }\,_{2}F_{1}\left(
\left.
\begin{array}{c}
1,1 \\
\beta%
\end{array}%
\right\vert \frac{1}{s}\right) $ \\ \midrule
$\frac{3}{2}$ & $\frac{1}{2}$ & $\frac{1}{\sqrt{\pi }s}\,_{2}F_{2}\left(
\left.
\begin{array}{c}
1,1 \\
\frac{1}{6},\frac{5}{6}%
\end{array}%
\right\vert \frac{4}{27s^{2}}\right) +\frac{1}{s^{2}}\,_{2}F_{2}\left(
\left.
\begin{array}{c}
1,\frac{3}{2} \\
\frac{2}{3},\frac{4}{3}%
\end{array}%
\right\vert \frac{4}{27s^{2}}\right) $ \\

\bottomrule
\end{tabular*}
\end{specialtable}

\begin{specialtable}[H]\ContinuedFloat
\small
\caption{{\em Cont.}}

\begin{tabular*}{\hsize}{c@{\extracolsep{\fill}}cc}
\toprule
$\alpha $ & $\beta $ & $\quad \mathcal{L}\left[ \mathrm{E}_{\alpha ,\beta
}\left( t\right) \right] $ \\ \midrule

$\frac{3}{2}$ & $1$ & $\frac{4}{3\sqrt{\pi }s^{2}}\,_{2}F_{2}\left( \left.
\begin{array}{c}
1,1 \\
\frac{5}{6},\frac{7}{6}%
\end{array}%
\right\vert \frac{4}{27s^{2}}\right) +\frac{1}{s}\,_{2}F_{2}\left( \left.
\begin{array}{c}
\frac{1}{2},1 \\
\frac{1}{3},\frac{2}{3}%
\end{array}%
\right\vert \frac{4}{27s^{2}}\right) $ \\ \midrule

$\frac{3}{2}$ & $\frac{3}{2}$ & $\frac{2}{\sqrt{\pi }s}\,_{2}F_{2}\left(
\left.
\begin{array}{c}
1,1 \\
\frac{5}{6},\frac{7}{6}%
\end{array}%
\right\vert \frac{4}{27s^{2}}\right) +\frac{1}{2s^{2}}\,_{2}F_{2}\left(
\left.
\begin{array}{c}
1,\frac{3}{2} \\
\frac{4}{3},\frac{5}{3}%
\end{array}%
\right\vert \frac{4}{27s^{2}}\right) $ \\ \midrule

$\frac{3}{2}$ & $2$ & $\frac{8}{15\sqrt{\pi }s^{2}}\,_{2}F_{2}\left( \left.
\begin{array}{c}
1,1 \\
\frac{7}{6},\frac{11}{6}%
\end{array}%
\right\vert \frac{4}{27s^{2}}\right) +\frac{1}{s}\,_{2}F_{2}\left( \left.
\begin{array}{c}
\frac{1}{2},1 \\
\frac{2}{3},\frac{4}{3}%
\end{array}%
\right\vert \frac{4}{27s^{2}}\right) $ \\ \midrule
$\frac{3}{2}$ & $\beta $ & $\frac{1}{\Gamma \left( \beta \right) s}%
\,_{3}F_{3}\left( \left.
\begin{array}{c}
\frac{1}{2},1,1 \\
\frac{\beta +1}{3},\frac{\beta +2}{3},\frac{\beta }{3}%
\end{array}%
\right\vert \frac{4}{27s^{2}}\right) +\frac{1}{\Gamma \left( \beta +\frac{3}{%
2}\right) s^{2}}\,_{3}F_{3}\left( \left.
\begin{array}{c}
1,1,\frac{3}{2} \\
\frac{2\beta +3}{6},\frac{2\beta +5}{6},\frac{2\beta +7}{6}%
\end{array}%
\right\vert \frac{4}{27s^{2}}\right) $ \\ \midrule
$2$ & $\frac{1}{2}$ & $\frac{1}{\sqrt{\pi }s}\,_{2}F_{2}\left( \left.
\begin{array}{c}
1,1 \\
\frac{1}{4},\frac{3}{4}%
\end{array}%
\right\vert \frac{1}{4s}\right) $ \\ \midrule
$2$ & $1$ & $\frac{1}{s}+\frac{\sqrt{\pi }}{2}s^{-3/2}e^{1/\left( 4s\right) }%
\mathrm{erf}\left( \frac{1}{2\sqrt{s}}\right) $ \\ \midrule
$2$ & $2$ & $\sqrt{\pi }s^{-1/2}e^{1/\left( 4s\right) }\mathrm{erf}\left(
\frac{1}{2\sqrt{s}}\right) $ \\ \midrule
$2$ & $3$ & $\frac{1}{2s}\,_{2}F_{2}\left( \left.
\begin{array}{c}
1,1 \\
\frac{3}{2},2%
\end{array}%
\right\vert \frac{1}{4s}\right) $ \\ \midrule
$2$ & $4$ & $\frac{1}{6s}\,_{2}F_{2}\left( \left.
\begin{array}{c}
1,1 \\
2,\frac{5}{2}%
\end{array}%
\right\vert \frac{1}{4s}\right) $ \\ \midrule
$2$ & $\beta $ & $\frac{1}{\Gamma \left( \beta \right) s}\,_{2}F_{2}\left(
\left.
\begin{array}{c}
1,1 \\
\frac{\beta +1}{2},\frac{\beta }{2}%
\end{array}%
\right\vert \frac{1}{4s}\right) $ \\ \bottomrule
\end{tabular*}%
\label{TableB}%
\end{specialtable}%

\vspace{-6pt}

\begin{specialtable}[H] \centering%
\caption{The Laplace transforms Mittag-Leffler functions derived for some values of
parameters $\alpha$ and $\beta$ by using (\ref{LT_Mittag-Leffler}).}%
\begin{tabular*}{\hsize}{c@{\extracolsep{\fill}}cc}
\toprule
\boldmath$\alpha $ & \boldmath$\beta $ & \boldmath$\quad \mathcal{L}\left[ \mathrm{E}_{\alpha ,\beta
}\left( t\right) \right] $ \\ \midrule
$3$ & $1$ & $\frac{1}{s}\,_{1}F_{2}\left( \left.
\begin{array}{c}
1 \\
\frac{1}{3},\frac{2}{3}%
\end{array}%
\right\vert \frac{1}{27s}\right) $ \\ \midrule
$3$ & $2$ & $\frac{1}{s}\,_{1}F_{2}\left( \left.
\begin{array}{c}
1 \\
\frac{2}{3},\frac{4}{3}%
\end{array}%
\right\vert \frac{1}{27s}\right) $ \\ \midrule
$3$ & $3$ & $\frac{1}{2s}\,_{1}F_{2}\left( \left.
\begin{array}{c}
1 \\
\frac{4}{3},\frac{5}{3}%
\end{array}%
\right\vert \frac{1}{27s}\right) $ \\ \midrule
$3$ & $\beta $ & $\frac{1}{\Gamma \left( \beta \right) s}\,_{2}F_{3}\left(
\left.
\begin{array}{c}
1,1 \\
\frac{\beta +1}{3},\frac{\beta +2}{3},\frac{\beta }{3}%
\end{array}%
\right\vert \frac{1}{27s}\right) $ \\ \midrule
$4$ & $1$ & $\frac{1}{s}\,_{1}F_{3}\left( \left.
\begin{array}{c}
1 \\
\frac{1}{4},\frac{1}{2},\frac{3}{4}%
\end{array}%
\right\vert \frac{1}{256s}\right) $ \\ \midrule
$4$ & $2$ & $\frac{1}{s}\,_{1}F_{3}\left( \left.
\begin{array}{c}
1 \\
\frac{1}{2},\frac{3}{4},\frac{5}{4}%
\end{array}%
\right\vert \frac{1}{256s}\right) $ \\ \midrule
$4$ & $3$ & $\frac{1}{2s}\,_{1}F_{3}\left( \left.
\begin{array}{c}
1 \\
\frac{3}{4},\frac{5}{4},\frac{3}{2}%
\end{array}%
\right\vert \frac{1}{256s}\right) $ \\ \midrule
$4$ & $4$ & $\frac{1}{6s}\,_{1}F_{3}\left( \left.
\begin{array}{c}
1 \\
\frac{5}{4},\frac{3}{2},\frac{7}{4}%
\end{array}%
\right\vert \frac{1}{256s}\right) $ \\

\bottomrule
\end{tabular*}
\end{specialtable}

\begin{specialtable}[H]\ContinuedFloat
\small
\caption{{\em Cont.}}
\begin{tabular*}{\hsize}{c@{\extracolsep{\fill}}cc}
\toprule
\boldmath$\alpha $ & \boldmath$\beta $ & \boldmath$\quad \mathcal{L}\left[ \mathrm{E}_{\alpha ,\beta
}\left( t\right) \right] $ \\ \midrule

$4$ & $\beta $ & $\frac{1}{\Gamma \left( \beta \right) s}\,_{2}F_{4}\left(
\left.
\begin{array}{c}
1,1 \\
\frac{\beta +1}{4},\frac{\beta +2}{4},\frac{\beta +3}{4},\frac{\beta }{4}%
\end{array}%
\right\vert \frac{1}{256s}\right) $ \\ \midrule
$5$ & $1$ & $\frac{1}{s}\,_{1}F_{4}\left( \left.
\begin{array}{c}
1 \\
\frac{1}{5},\frac{2}{5},\frac{3}{5},\frac{4}{5}%
\end{array}%
\right\vert \frac{1}{3125s}\right) $ \\ \midrule
$5$ & $2$ & $\frac{1}{s}\,_{1}F_{4}\left( \left.
\begin{array}{c}
1 \\
\frac{2}{5},\frac{3}{5},\frac{4}{5},\frac{6}{5}%
\end{array}%
\right\vert \frac{1}{3125s}\right) $ \\ \midrule
$5$ & $3$ & $\frac{1}{2s}\,_{1}F_{4}\left( \left.
\begin{array}{c}
1 \\
\frac{3}{5},\frac{4}{5},\frac{6}{5},\frac{7}{5}%
\end{array}%
\right\vert \frac{1}{3125s}\right) $ \\ \midrule

$5$ & $4$ & $\frac{1}{6s}\,_{1}F_{4}\left( \left.
\begin{array}{c}
1 \\
\frac{4}{5},\frac{6}{5},\frac{7}{5},\frac{8}{5}%
\end{array}%
\right\vert \frac{1}{3125s}\right) $ \\ \midrule
$5$ & $5$ & $\frac{1}{24s}\,_{1}F_{4}\left( \left.
\begin{array}{c}
1 \\
\frac{6}{5},\frac{7}{5},\frac{8}{5},\frac{9}{5}%
\end{array}%
\right\vert \frac{1}{3125s}\right) $ \\ \midrule
$5$ & $\beta $ & $\frac{1}{\Gamma \left( \beta \right) s}\,_{2}F_{5}\left(
\left.
\begin{array}{c}
1,1 \\
\frac{\beta +1}{5},\frac{\beta +2}{5},\frac{\beta +3}{5},\frac{\beta +4}{5},%
\frac{\beta }{5}%
\end{array}%
\right\vert \frac{1}{3125s}\right) $ \\ \bottomrule
\end{tabular*}%
\label{TableBa}%
\end{specialtable}%

\section{Representations of the Whittaker Functions and their Laplace
Transforms \label{Appendix: Whittaker}}

The Whittaker functions $\mathrm{M}_{\kappa ,\mu }\left( x\right) $ and $%
\mathrm{W}_{\kappa ,\mu }\left( x\right) $ defined in (\ref{Whittaker_def})
were derived by using the MATHEMATICA program, and they are presented in Tables 
\ref{TableC}, \ref{TableCb}, \ref{TableD} and \ref{TableDa}. The
corresponding Laplace transforms are in Tables \ref{TableF}, \ref{TableFa}, \ref%
{TableFb} and \ref{TableG}. Most of the reported results in these tables are
unknown in the mathematical reference literature.

\begin{specialtable}[H] \centering%
\caption{The Whittaker functions $\mathrm{M}_{\kappa, \mu}$ derived for some values of
parameters $\kappa$ and $\mu$ by using (\ref{Whittaker_def}).}%
\begin{tabular*}{\hsize}{c@{\extracolsep{\fill}}cc}
\toprule
\boldmath$\kappa $ & \boldmath$\mu $ & \boldmath$\mathrm{M}_{\kappa ,\mu }\left( x\right) $ \\
\midrule
$-\frac{5}{2}$ & $0$ & $\frac{\sqrt{x}}{2}e^{x/2}\left[ x\left( x+4\right) +2%
\right] $ \\ \midrule
$-\frac{3}{2}$ & $0$ & $\sqrt{x}e^{x/2}\left( x+1\right) $ \\ \midrule
$-\frac{3}{2}$ & $\frac{1}{2}$ & $\frac{x}{3}\left[ \left( 2x+3\right)
I_{0}\left( \frac{x}{2}\right) +\left( 2x+1\right) I_{1}\left( \frac{x}{2}%
\right) \right] $ \\ \midrule
$-\frac{3}{2}$ & $1$ & $e^{x/2}x^{3/2}$ \\ \midrule
$-\frac{3}{2}$ & $\frac{3}{2}$ & $\frac{4}{5}\left[ x\left( 2x-3\right)
I_{0}\left( \frac{x}{2}\right) +\left[ x\left( 2x-3\right) +4\right]
I_{1}\left( \frac{x}{2}\right) \right] $ \\ \midrule
$-\frac{3}{2}$ & $2$ & $4x^{-3/2}e^{-x/2}\left[ e^{x}\left(
x^{3}-3x^{2}+6x-6\right) +6\right] $ \\ \midrule
$-\frac{3}{2}$ & $\frac{5}{2}$ & $\frac{32}{7x}\left[ x\left(
2x^{2}-9x+24\right) I_{0}\left( \frac{x}{2}\right) +\left(
2x^{3}-11x^{2}+36x-96\right) I_{1}\left( \frac{x}{2}\right) \right] $ \\
\midrule
$-\frac{3}{2}$ & $3$ & $30x^{-5/2}e^{-x/2}\left[ e^{x}\left(
x^{4}-8x^{3}+36x^{2}-96x+120\right) -24\left( x+5\right) \right] $ \\ \midrule
$-\frac{1}{6}$ & $0$ & $e^{-x/2}\sqrt{x}L_{-2/3}\left( x\right) $ \\ \midrule
$-\frac{1}{4}$ & $0$ & $e^{-x/2}\sqrt{x}L_{-3/4}\left( x\right) $ \\ \midrule
$-\frac{1}{4}$ & $\frac{1}{4}$ & $\frac{\sqrt{\pi }}{2}e^{x/2}x^{1/4}\mathrm{%
erf}\left( \sqrt{x}\right) $ \\ \midrule
$-\frac{1}{3}$ & $0$ & $e^{-x/2}\sqrt{x}L_{-5/6}\left( x\right) $ \\ \midrule
$-\frac{1}{2}$ & $\frac{1}{2}$ & $x\left[ I_{0}\left( \frac{x}{2}\right)
+I_{1}\left( \frac{x}{2}\right) \right] $ \\ \midrule
$-\frac{1}{2}$ & $1$ & $x^{-1/2}e^{-x/2}\left[ 2e^{x}\left( x-1\right) +2%
\right] $ \\ \midrule
$-\frac{1}{2}$ & $\frac{3}{2}$ & $4\left[ xI_{0}\left( \frac{x}{2}\right)
+\left( x-4\right) I_{1}\left( \frac{x}{2}\right) \right] $ \\ \midrule
$-\frac{1}{2}$ & $2$ & $12x^{-3/2}e^{-x/2}\left[ e^{x}\left(
x^{2}-4x+6\right) -2\left( x+3\right) \right] $ \\

\bottomrule
\end{tabular*}
\end{specialtable}

\begin{specialtable}[H]\ContinuedFloat
\small
\caption{{\em Cont.}}
\begin{tabular*}{\hsize}{c@{\extracolsep{\fill}}cc}
\toprule
\boldmath$\kappa $ & \boldmath$\mu $ & \boldmath$\mathrm{M}_{\kappa ,\mu }\left( x\right) $ \\
\midrule

$0$ & $\frac{1}{8}$ & $x^{5/8}\,_{0}F_{1}\left( \left.
\begin{array}{c}
- \\
\frac{9}{8}%
\end{array}%
\right\vert \frac{x^{2}}{16}\right) $ \\ \midrule
$0$ & $\frac{1}{7}$ & $x^{9/14}\,_{0}F_{1}\left( \left.
\begin{array}{c}
- \\
\frac{8}{7}%
\end{array}%
\right\vert \frac{x^{2}}{16}\right) $ \\ \midrule
$0$ & $\frac{1}{6}$ & $x^{2/3}\,_{0}F_{1}\left( \left.
\begin{array}{c}
- \\
\frac{7}{6}%
\end{array}%
\right\vert \frac{x^{2}}{16}\right) $ \\ \midrule
$0$ & $\frac{1}{5}$ & $x^{7/10}\,_{0}F_{1}\left( \left.
\begin{array}{c}
- \\
\frac{6}{5}%
\end{array}%
\right\vert \frac{x^{2}}{16}\right) $ \\ \midrule
$0$ & $\frac{1}{4}$ & $x^{3/4}\,_{0}F_{1}\left( \left.
\begin{array}{c}
- \\
\frac{5}{4}%
\end{array}%
\right\vert \frac{x^{2}}{16}\right) $ \\ \midrule
$0$ & $\frac{1}{3}$ & $x^{5/6}\,_{0}F_{1}\left( \left.
\begin{array}{c}
- \\
\frac{4}{3}%
\end{array}%
\right\vert \frac{x^{2}}{16}\right) $ \\ \bottomrule
\end{tabular*}%
\label{TableC}%
\end{specialtable}%

\vspace{-6pt}
\begin{specialtable}[H] \centering%
\caption{The Whittaker functions $\mathrm{M}_{\kappa, \mu}$ derived for some values of
parameters $\kappa$ and $\mu$ by using (\ref{Whittaker_def}).}%
\begin{tabular*}{\hsize}{c@{\extracolsep{\fill}}cc}
\toprule
\boldmath$\kappa $ & \boldmath$\mu $ & \boldmath$\mathrm{M}_{\kappa ,\mu }\left( x\right) $ \\
\midrule
$0$ & $\frac{1}{2}$ & $2\sinh \left( \frac{x}{2}\right) $ \\ \midrule
$0$ & $1$ & $4\sqrt{x}I_{1}\left( \frac{x}{2}\right) $ \\ \midrule
$0$ & $\frac{3}{2}$ & $12\left[ \cosh \left( \frac{x}{2}\right) -\frac{2}{x}%
\sinh \left( \frac{x}{2}\right) \right] $ \\ \midrule
$0$ & $2$ & $32\sqrt{x}I_{2}\left( \frac{x}{2}\right) $ \\ \midrule
$0$ & $\frac{5}{2}$ & $\frac{120}{x^{2}}\left[ \left( x^{2}+12\right) \sinh
\left( \frac{x}{2}\right) -6x\cosh \left( \frac{x}{2}\right) \right] $ \\
\midrule
$\frac{1}{6}$ & $0$ & $e^{-x/2}\sqrt{x}L_{-1/3}\left( x\right) $ \\ \midrule
$\frac{1}{4}$ & $-\frac{5}{4}$ & $x^{-3/4}e^{-x/2}\left( \frac{2x}{3}%
+1\right) $ \\ \midrule
$\frac{1}{4}$ & $-\frac{3}{4}$ & $x^{-1/4}e^{x/2}$ \\ \midrule
$\frac{1}{4}$ & $-\frac{1}{4}$ & $x^{1/4}e^{-x/2}$ \\ \midrule
$\frac{1}{4}$ & $0$ & $e^{-x/2}\sqrt{x}L_{-1/4}\left( x\right) $ \\ \midrule
$\frac{1}{3}$ & $0$ & $e^{-x/2}\sqrt{x}L_{-1/6}\left( x\right) $ \\ \midrule
$\frac{1}{2}$ & $0$ & $e^{-x/2}\sqrt{x}$ \\ \midrule
$\frac{1}{2}$ & $\frac{1}{2}$ & $x\left[ I_{0}\left( \frac{x}{2}\right)
-I_{1}\left( \frac{x}{2}\right) \right] $ \\ \midrule
$\frac{1}{2}$ & $1$ & $2x^{-1/2}e^{-x/2}\left( e^{x}-x-1\right) $ \\ \midrule
$\frac{1}{2}$ & $\frac{3}{2}$ & $4\left[ -xI_{0}\left( \frac{x}{2}\right)
+\left( x+4\right) I_{1}\left( \frac{x}{2}\right) \right] $ \\ \midrule
$\frac{1}{2}$ & $\frac{5}{2}$ & $32\left[ \left( x+8\right) I_{0}\left(
\frac{x}{2}\right) -\left( x+4+\frac{32}{x}\right) I_{1}\left( \frac{x}{2}%
\right) \right] $ \\ \midrule

$1$ & $0$ & $\sqrt{x}\left[ -\left( x-1\right) I_{0}\left( \frac{x}{2}%
\right) +xI_{1}\left( \frac{x}{2}\right) \right] $ \\ \midrule
$1$ & $\frac{1}{2}$ & $x\,e^{-x/2}$ \\ \midrule
$1$ & $1$ & $\frac{4}{3}\sqrt{x}\left[ xI_{0}\left( \frac{x}{2}\right)
-\left( x+1\right) I_{1}\left( \frac{x}{2}\right) \right] $ \\ \midrule
$\frac{3}{2}$ & $0$ & $-\sqrt{x}e^{-x/2}\left( x-1\right) $ \\ \midrule
$\frac{3}{2}$ & $\frac{1}{2}$ & $-\frac{x}{3}\left[ \left( 2x-3\right)
I_{0}\left( \frac{x}{2}\right) +\left( 1-2x\right) I_{1}\left( \frac{x}{2}%
\right) \right] $ \\

\bottomrule
\end{tabular*}
\end{specialtable}

\begin{specialtable}[H]\ContinuedFloat
\small
\caption{{\em Cont.}}
\begin{tabular*}{\hsize}{c@{\extracolsep{\fill}}cc}
\toprule
\boldmath$\kappa $ & \boldmath$\mu $ & \boldmath$\mathrm{M}_{\kappa ,\mu }\left( x\right) $ \\
\midrule

$\frac{3}{2}$ & $1$ & $x^{3/2}\,e^{-x/2}$ \\ \midrule
$\frac{3}{2}$ & $\frac{3}{2}$ & $\frac{4}{5}\left[ x\left( 2x+1\right)
I_{0}\left( \frac{x}{2}\right) -\left( 2x^{2}+3x+4\right) I_{1}\left( \frac{x%
}{2}\right) \right] $ \\ \midrule
$2$ & $0$ & $\frac{1}{3}\sqrt{x}\left[ \left( 2x^{2}-6x+3\right) I_{0}\left(
\frac{x}{2}\right) -2x\left( x-2\right) I_{1}\left( \frac{x}{2}\right) %
\right] $ \\ \midrule
$2$ & $\frac{1}{2}$ & $-\frac{1}{2}e^{-x/2}x\left( x-2\right) $ \\ \midrule
$2$ & $1$ & $-\frac{4}{15}\sqrt{x}\left[ 2x\left( x-2\right) I_{0}\left(
\frac{x}{2}\right) +\left( -2x^{2}+2x+1\right) I_{1}\left( \frac{x}{2}%
\right) \right] $ \\ \midrule
$2$ & $\frac{3}{2}$ & $x^{2}\,e^{-x/2}$ \\ \midrule
$2$ & $2$ & $\frac{32}{35\sqrt{x}}\left[ x\left( 2x^{2}+2x+3\right)
I_{0}\left( \frac{x}{2}\right) -\left( x^{3}+2x^{2}+4x+6\right) I_{1}\left(
\frac{x}{2}\right) \right] $ \\ \midrule
$\frac{5}{2}$ & $0$ & $\frac{1}{2}e^{-x/2}\sqrt{x}\left( x^{2}-4x+2\right) $
\\ \midrule

$\frac{5}{2}$ & $\frac{1}{2}$ & $\frac{x}{15}\left[ \left(
4x^{2}-18x+15\right) I_{0}\left( \frac{x}{2}\right) +\left(
-4x^{2}+14x-3\right) I_{1}\left( \frac{x}{2}\right) \right] $ \\ \midrule
$\frac{5}{2}$ & $1$ & $-\frac{1}{3}e^{-x/2}x^{3/2}\left( x-3\right) $ \\
\midrule
$\frac{5}{2}$ & $2$ & $x^{5/2}\,e^{-x/2}$ \\ \midrule
$3$ & $\frac{1}{2}$ & $\frac{1}{6}e^{-x/2}x\left( x^{2}-6x+6\right) $ \\
\midrule
$3$ & $1$ & $\frac{4\sqrt{x}}{105}\left[ x\left( 4x^{2}-24x+27\right)
I_{0}\left( \frac{x}{2}\right) -\left( 4x^{3}+20x^{2}+9x+3\right)
I_{1}\left( \frac{x}{2}\right) \right] $ \\ \midrule
$3$ & $\frac{3}{2}$ & $-\frac{1}{4}e^{-x/2}x^{2}\left( x-4\right) $ \\ \midrule
$3$ & $\frac{5}{2}$ & $x^{3}\,e^{-x/2}$ \\ \midrule
$\frac{7}{2}$ & $0$ & $-\frac{1}{6}e^{-x/2}\sqrt{x}\left(
x^{3}-9x^{2}+18x-6\right) $ \\ \midrule
$4$ & $\frac{1}{2}$ & $-\frac{1}{24}e^{-x/2}x\left(
x^{3}-12x^{2}+36x-24\right) $ \\ \midrule
$4$ & $\frac{3}{2}$ & $\frac{1}{20}e^{-x/2}x^{2}\left( x^{2}-10x+20\right) $
\\ \bottomrule
\end{tabular*}%
\label{TableCb}%
\end{specialtable}%

\vspace{-6pt}

\begin{specialtable}[H] \centering%
\caption{The Laplace transforms of the Whittaker function $\mathrm{M}_{\kappa, \mu}$ derived for some values of
parameters $\kappa$ and $\mu$.}%
\begin{tabular*}{\hsize}{c@{\extracolsep{\fill}}cc}
\toprule
\boldmath$\kappa $ & \boldmath$\mu $ & \boldmath$\mathcal{L}\left[ \mathrm{M}_{\kappa ,\mu }\left(
t\right) \right] $ \\ \midrule
$-\frac{5}{2}$ & $0$ & $\sqrt{\frac{\pi }{2}}\frac{8s^{2}+16s+5}{\left(
2s-1\right) ^{7/2}}$ \\ \midrule
$-\frac{3}{2}$ & $0$ & $\frac{2\sqrt{2\pi }\left( s+1\right) }{\left(
2s-1\right) ^{5/2}}$ \\ \midrule
$-\frac{3}{2}$ & $\frac{1}{2}$ & $\left\{
\begin{array}{cc}
\frac{4\sqrt{4s^{2}-1}}{\left( 2s-1\right) ^{3}}, & s>\frac{1}{2} \\
0, & s<\frac{1}{2}%
\end{array}%
\right. $ \\ \midrule
$-\frac{3}{2}$ & $1$ & $\frac{2\sqrt{2\pi }}{\left( 2s-1\right) ^{5/2}}$ \\
\midrule
$-\frac{3}{2}$ & $\frac{3}{2}$ & $\left\{
\begin{array}{cc}
\frac{16\left\{ -4s\left[ 4s^{2}-2s\left( \sqrt{4s^{2}-1}+3\right) +3\left(
\sqrt{4s^{2}-1}+1\right) \right] +7\sqrt{4s^{2}-1}+2\right\} }{5\left(
2s-1\right) ^{3}}, & s>\frac{1}{2} \\
-\frac{32}{5}, & s<\frac{1}{2}%
\end{array}%
\right. $ \\ \midrule

$-\frac{1}{6}$ & $0$ & $\left\{
\begin{array}{cc}
\frac{\sqrt{2\pi }\left[ \left( 6s+3\right) \,_{2}F_{1}\left( \left.
\begin{array}{c}
-\frac{1}{2},\frac{2}{3} \\
1%
\end{array}%
\right\vert \frac{2}{2s+1}\right) +2\,_{2}F_{1}\left( \left.
\begin{array}{c}
\frac{1}{2},\frac{2}{3} \\
1%
\end{array}%
\right\vert \frac{2}{2s+1}\right) \right] }{5\left( 2s-1\right) ^{3}}, & s>%
\frac{1}{2} \\
-\frac{\pi \,_{2}F_{1}\left( \left.
\begin{array}{c}
\frac{2}{3},\frac{2}{3} \\
\frac{1}{6}%
\end{array}%
\right\vert s+\frac{1}{2}\right) }{\left( s+\frac{1}{2}\right) ^{5/6}\Gamma
\left( \frac{1}{6}\right) \Gamma \left( \frac{1}{3}\right) }, & s<\frac{1}{2}%
\end{array}%
\right. $ \\

\bottomrule
\end{tabular*}
\end{specialtable}

\begin{specialtable}[H]\ContinuedFloat
\small
\caption{{\em Cont.}}

\begin{tabular*}{\hsize}{c@{\extracolsep{\fill}}cc}
\toprule
\boldmath$\kappa $ & \boldmath$\mu $ & \boldmath$\mathcal{L}\left[ \mathrm{M}_{\kappa ,\mu }\left(
t\right) \right] $ \\ \midrule

$-\frac{1}{4}$ & $\frac{1}{4}$ & $\frac{2^{15/4}\Gamma \left( \frac{11}{4}%
\right) \left[ \left( 6s+3\right) \,_{2}F_{1}\left( \left.
\begin{array}{c}
-\frac{1}{4},\frac{1}{2} \\
\frac{3}{2}%
\end{array}%
\right\vert \frac{2}{1-2s}\right) -4\left( \frac{2}{2s+1}+1\right) ^{1/4}%
\right] }{21\left( 2s-1\right) ^{7/4}\left( 2s+1\right) }$ \\ \midrule
$-\frac{1}{3}$ & $0$ & $\left\{
\begin{array}{cc}
\frac{\sqrt{2\pi }\left[ \left( 6s+3\right) \,_{2}F_{1}\left( \left.
\begin{array}{c}
-\frac{1}{2},\frac{5}{6} \\
1%
\end{array}%
\right\vert \frac{2}{2s+1}\right) +4\,_{2}F_{1}\left( \left.
\begin{array}{c}
\frac{1}{2},\frac{5}{6} \\
1%
\end{array}%
\right\vert \frac{2}{2s+1}\right) \right] }{3\left( 2s-1\right) \left(
2s+1\right) ^{3/2}}, & s>\frac{1}{2} \\
-\frac{\pi \,_{2}F_{1}\left( \left.
\begin{array}{c}
\frac{5}{6},\frac{5}{6} \\
\frac{1}{3}%
\end{array}%
\right\vert s+\frac{1}{2}\right) }{\left( s+\frac{1}{2}\right) ^{2/3}\Gamma
\left( \frac{1}{6}\right) \Gamma \left( \frac{1}{3}\right) }, & s<\frac{1}{2}%
\end{array}%
\right. $ \\ \midrule
$-\frac{1}{2}$ & $0$ & $\frac{2\sqrt{2\pi }}{\left( 2s-1\right) ^{3/2}}$ \\
\midrule
$-\frac{1}{2}$ & $\frac{1}{2}$ & $\left\{
\begin{array}{cc}
\frac{4}{\left( 2s-1\right) \sqrt{4s^{2}-1}}, & s>\frac{1}{2} \\
0, & s<\frac{1}{2}%
\end{array}%
\right. $ \\ \midrule
$-\frac{1}{2}$ & $1$ & $2\sqrt{2\pi }\left[ \frac{1}{\sqrt{2s+1}}-\frac{1}{%
\sqrt{2s-1}}+\frac{1}{\left( 2s-1\right) ^{3/2}}\right] $ \\ \midrule
$0$ & $1$ & $\left\{
\begin{array}{cc}
\frac{2^{3/2}\left[ \left( 1-2s\right) \mathrm{K}\left( \frac{2}{2s+1}%
\right) +2s\,\mathrm{E}\left( \frac{2}{2s+1}\right) \right] }{\sqrt{\pi }%
\left( 2s-1\right) \sqrt{2s+1}}, & s>\frac{1}{2} \\
\frac{8\left[ \left( 1-2s\right) \mathrm{K}\left( s+\frac{1}{2}\right) +4s\,%
\mathrm{E}\left( s+\frac{1}{2}\right) \right] }{\sqrt{\pi }\left(
4s^{2}-1\right) }, & s<\frac{1}{2}%
\end{array}%
\right. $ \\ \midrule
$0$ & $2$ & $\left\{
\begin{array}{cc}
\frac{64\left[ 8s\left( 1-2s\right) \mathrm{K}\left( \frac{2}{2s+1}\right)
+\left( 16s^{2}-3\right) \,\mathrm{E}\left( \frac{2}{2s+1}\right) \right] }{%
\sqrt{\pi }\left( 2s-1\right) \sqrt{s+\frac{1}{2}}}, & s>\frac{1}{2} \\
\frac{64\left[ \left( -16s^{2}+2s+3\right) \mathrm{K}\left( s+\frac{1}{2}%
\right) +\left( 32s^{2}-6\right) \,\mathrm{E}\left( s+\frac{1}{2}\right) %
\right] }{\sqrt{\pi }\left( 4s^{2}-1\right) }, & s<\frac{1}{2}%
\end{array}%
\right. $ \\ \midrule
$\frac{1}{6}$ & $0$ & $\left\{
\begin{array}{cc}
\frac{\sqrt{2\pi }\left[ \left( 6s+3\right) \,_{2}F_{1}\left( \left.
\begin{array}{c}
-\frac{1}{2},\frac{1}{3} \\
1%
\end{array}%
\right\vert \frac{2}{2s+1}\right) -2\,_{2}F_{1}\left( \left.
\begin{array}{c}
\frac{1}{3},\frac{1}{2} \\
1%
\end{array}%
\right\vert \frac{2}{2s+1}\right) \right] }{3\left( 2s-1\right) \left(
2s+1\right) ^{3/2}}, & s>\frac{1}{2} \\
\frac{\sqrt{2\pi }\Gamma \left( \frac{7}{6}\right) \,_{2}F_{1}\left( \left.
\begin{array}{c}
\frac{1}{3},\frac{1}{3} \\
-\frac{1}{6}%
\end{array}%
\right\vert s+\frac{1}{2}\right) }{\left( 2s+1\right) ^{7/6}\Gamma \left(
\frac{5}{6}\right) \Gamma \left( \frac{1}{3}\right) }, & s<\frac{1}{2}%
\end{array}%
\right. $ \\ \bottomrule
\end{tabular*}%
\label{TableF}%
\end{specialtable}%
\vspace{-6pt}
\begin{specialtable}[H] \centering%
\caption{The Laplace transforms of the Whittaker functions $\mathrm{M}_{\kappa, \mu}$ derived for some values of
parameters $\kappa$ and $\mu$.}%
\begin{tabular*}{\hsize}{c@{\extracolsep{\fill}}cc}
\toprule
\boldmath$\kappa $ & \boldmath$\mu $ & \boldmath$\mathcal{L}\left[ \mathrm{M}_{\kappa ,\mu }\left(
t\right) \right] $ \\ \midrule
$\frac{1}{4}$ & $-\frac{5}{4}$ & $\frac{3\left( s+\frac{1}{2}\right) \Gamma
\left( \frac{1}{4}\right) +2\,\Gamma \left( \frac{5}{4}\right) }{3\left( s+%
\frac{1}{2}\right) ^{5/4}}$ \\ \midrule
$\frac{1}{4}$ & $-\frac{3}{4}$ & $\frac{\Gamma \left( \frac{3}{4}\right) }{%
\left( s-\frac{1}{2}\right) ^{3/4}}$ \\ \midrule
$\frac{1}{4}$ & $-\frac{1}{4}$ & $\frac{\Gamma \left( \frac{5}{4}\right) }{%
\left( s+\frac{1}{2}\right) ^{5/4}}$ \\ \midrule
$0$ & $\frac{1}{2}$ & $\frac{4}{4s^{2}-1}$ \\ \midrule
$\frac{1}{3}$ & $0$ & $\left\{
\begin{array}{cc}
\frac{\sqrt{2\pi }\left[ \left( 6s+3\right) \,_{2}F_{1}\left( \left.
\begin{array}{c}
-\frac{1}{2},\frac{1}{6} \\
1%
\end{array}%
\right\vert \frac{2}{2s+1}\right) -4\,_{2}F_{1}\left( \left.
\begin{array}{c}
\frac{1}{6},\frac{1}{2} \\
1%
\end{array}%
\right\vert \frac{2}{2s+1}\right) \right] }{3\left( 2s-1\right) \left(
2s+1\right) ^{3/2}}, & s>\frac{1}{2} \\
\frac{\Gamma \left( \frac{1}{6}\right) \Gamma \left( \frac{1}{3}\right)
\,_{2}F_{1}\left( \left.
\begin{array}{c}
\frac{1}{6},\frac{1}{6} \\
-\frac{1}{3}%
\end{array}%
\right\vert s+\frac{1}{2}\right) }{2^{2/3}\sqrt{3}\pi \left( 2s+1\right)
^{4/3}}, & s<\frac{1}{2}%
\end{array}%
\right. $ \\ \midrule

$\frac{1}{2}$ & $\frac{1}{2}$ & $\left\{
\begin{array}{cc}
\frac{4}{\left( 2s+1\right) \sqrt{4s^{2}-1}}, & s>\frac{1}{2} \\
0, & s<\frac{1}{2}%
\end{array}%
\right. $ \\

\bottomrule
\end{tabular*}
\end{specialtable}

\begin{specialtable}[H]\ContinuedFloat
\small
\caption{{\em Cont.}}

\begin{tabular*}{\hsize}{c@{\extracolsep{\fill}}cc}
\toprule
\boldmath$\kappa $ & \boldmath$\mu $ & \boldmath$\mathcal{L}\left[ \mathrm{M}_{\kappa ,\mu }\left(
t\right) \right] $ \\ \midrule

$\frac{1}{2}$ & $0$ & $\frac{\sqrt{2\pi }}{\left( 2s+1\right) ^{3/2}}$ \\
\midrule

$\frac{1}{2}$ & $1$ & $2\sqrt{2\pi }\left[ -\frac{1}{\sqrt{2s+1}}-\frac{1}{%
\left( 2s+1\right) ^{3/2}}+\frac{1}{\sqrt{2s-1}}\right] $ \\ \midrule
$1$ & $0$ & $\left\{
\begin{array}{cc}
\frac{2^{3/2}\left[ -\mathrm{K}\left( \frac{2}{2s+1}\right) +2\,\mathrm{E}%
\left( \frac{2}{2s+1}\right) \right] }{\sqrt{\pi }\left( 2s+1\right) ^{3/2}},
& s>\frac{1}{2} \\
\frac{2\left( 2s-3\right) \mathrm{K}\left( s+\frac{1}{2}\right) +8\,\mathrm{E%
}\left( s+\frac{1}{2}\right) }{\sqrt{\pi }\left( 2s+1\right) ^{2}}, & s<%
\frac{1}{2}%
\end{array}%
\right. $ \\ \midrule
$1$ & $\frac{1}{2}$ & $\frac{4}{\left( 2s+1\right) ^{2}}$ \\ \midrule
$1$ & $1$ & $\left\{
\begin{array}{cc}
-\frac{2^{7/2}\left[ -2\left( s+1\right) \mathrm{K}\left( \frac{2}{2s+1}%
\right) +\left( 2s+3\right) \,\mathrm{E}\left( \frac{2}{2s+1}\right) \right]
}{3\sqrt{\pi }\left( 2s+1\right) ^{3/2}}, & s>\frac{1}{2} \\
\frac{8\left( 2s+5\right) \mathrm{K}\left( s+\frac{1}{2}\right) -8\left(
4s+6\right) \,\mathrm{E}\left( s+\frac{1}{2}\right) }{3\sqrt{\pi }\left(
2s+1\right) ^{2}}, & s<\frac{1}{2}%
\end{array}%
\right. $ \\ \midrule
$\frac{3}{2}$ & $0$ & $\frac{2\sqrt{2\pi }\left( s-1\right) }{\left(
2s+1\right) ^{5/2}}$ \\ \midrule
$\frac{3}{2}$ & $\frac{1}{2}$ & $\left\{
\begin{array}{cc}
\frac{4\sqrt{4s^{2}-1}}{\left( 2s+1\right) ^{3}}, & s>\frac{1}{2} \\
0, & s<\frac{1}{2}%
\end{array}%
\right. $ \\ \midrule
$\frac{3}{2}$ & $1$ & $\frac{3\sqrt{2\pi }}{\left( 2s+1\right) ^{5/2}}$ \\
\midrule
$2$ & $0$ & $\left\{
\begin{array}{cc}
\frac{2^{3/2}\left[ 4\left( 1-2s\right) \mathrm{K}\left( \frac{2}{2s+1}%
\right) +\left( 14s-9\right) \,\mathrm{E}\left( \frac{2}{2s+1}\right) \right]
}{3\sqrt{\pi }\left( 2s+1\right) ^{5/2}}, & s>\frac{1}{2} \\
\frac{2\left( 2s-1\right) \left( 6s-13\right) \mathrm{K}\left( s+\frac{1}{2}%
\right) +4\left( 14s-9\right) \,\mathrm{E}\left( s+\frac{1}{2}\right) }{3%
\sqrt{\pi }\left( 2s+1\right) ^{3}}, & s<\frac{1}{2}%
\end{array}%
\right. $ \\ \midrule
$2$ & $\frac{1}{2}$ & $\frac{8s-4}{\left( 2s+1\right) ^{3}}$ \\ \midrule
$2$ & $1$ & $\left\{
\begin{array}{cc}
\frac{2^{7/2}\left[ \left[ -5+4s\left( s+2\right) \right] \mathrm{K}\left(
\frac{2}{2s+1}\right) -2\left[ \left( 2s+5\right) -6\right] \,\mathrm{E}%
\left( \frac{2}{2s+1}\right) \right] }{15\sqrt{\pi }\left( 2s+1\right) ^{5/2}%
}, & s>\frac{1}{2} \\
\frac{8\left( 2s-1\right) \left( 2s+17\right) \mathrm{K}\left( s+\frac{1}{2}%
\right) -32\left[ s\left( 2s+5\right) -6\right] \,\mathrm{E}\left( s+\frac{1%
}{2}\right) }{15\sqrt{\pi }\left( 2s+1\right) ^{3}}, & s<\frac{1}{2}%
\end{array}%
\right. $ \\ \midrule
$2$ & $\frac{3}{2}$ & $\frac{16}{\left( 2s+1\right) ^{3}}$ \\ \bottomrule
\end{tabular*}%
\label{TableFa}%
\end{specialtable}%
\vspace{-6pt}
\begin{specialtable}[H] \centering%
\caption{The Laplace transforms of the Whittaker functions $\mathrm{M}_{\kappa, \mu}$ derived for some values of
parameters $\kappa$ and $\mu$.}%
\begin{tabular*}{\hsize}{c@{\extracolsep{\fill}}cc}
\toprule
\boldmath$\kappa $ & \boldmath$\mu $ & \boldmath$\mathcal{L}\left[ \mathrm{M}_{\kappa ,\mu }\left(
t\right) \right] $ \\ \midrule
$\frac{5}{2}$ & $1$ & $\frac{2\sqrt{2\pi }\left( 3s-1\right) }{\left(
2s+1\right) ^{7/2}}$ \\ \midrule
$\frac{5}{2}$ & $2$ & $\frac{15\sqrt{2\pi }}{\left( 2s+1\right) ^{7/2}}$ \\
\midrule
$3$ & $\frac{1}{2}$ & $\frac{4\left( 1-2s\right) ^{2}}{\left( 2s+1\right)
^{4}}$ \\ \midrule
$3$ & $\frac{3}{2}$ & $\frac{8\left( 4s-1\right) }{\left( 2s+1\right) ^{4}}$
\\ \midrule
$3$ & $\frac{5}{2}$ & $\frac{96}{\left( 2s+1\right) ^{4}}$ \\ \midrule
$\frac{7}{2}$ & $0$ & $\frac{\sqrt{2\pi }\left( 8s^{3}-24s^{2}+15s-3\right)
}{\left( 2s+1\right) ^{9/2}}$ \\ \midrule
$4$ & $\frac{1}{2}$ & $\frac{4\left( 2s-1\right) ^{3}}{\left( 2s+1\right)
^{5}}$ \\ \midrule
$4$ & $\frac{3}{2}$ & $\frac{32\left( 10s^{2}-5s+1\right) }{5\left(
2s+1\right) ^{5}}$ \\ \bottomrule
\end{tabular*}%
\label{TableFb}%
\end{specialtable}%

\begin{specialtable}[H] \centering%
\caption{The Whittaker functions $\mathrm{W}_{\kappa, \mu}$ derived for some values of
parameters $\kappa$ and $\mu$ by using~(\ref{Whittaker_def}).}%
\begin{tabular*}{\hsize}{c@{\extracolsep{\fill}}cc}
\toprule
\boldmath$\kappa $ & \boldmath$\mu $ &\boldmath $\mathrm{W}_{\kappa ,\mu }\left( x\right) $ \\
\midrule
$-\frac{5}{2}$ & $0$ & $\frac{\sqrt{x}}{4}e^{-x/2}\left[ e^{x}\left(
x^{2}+4x+2\right) \Gamma \left( 0,x\right) -x-3\right] $ \\ \midrule
$-\frac{3}{2}$ & $0$ & $\sqrt{x}e^{-x/2}\left[ e^{x}\left( x+1\right) \Gamma
\left( 0,x\right) -1\right] $ \\ \midrule
$-\frac{3}{2}$ & $1$ & $x^{3/2}e^{x/2}\Gamma \left( -2,x\right) $ \\ \midrule
$-\frac{1}{4}$ & $\frac{1}{4}$ & $x^{1/4}e^{x/2}\Gamma \left( \frac{1}{2}%
,x\right) $ \\ \midrule
$-\frac{1}{2}$ & $0$ & $x^{1/2}e^{x/2}\Gamma \left( 0,x\right) $ \\ \midrule
$-\frac{1}{2}$ & $1$ & $x^{-1/2}e^{-x/2}$ \\ \midrule
$-\frac{1}{2}$ & $2$ & $x^{-3/2}e^{-x/2}\left( x+3\right) $ \\ \midrule
$-\frac{1}{2}$ & $3$ & $x^{-5/2}e^{-x/2}\left( x^{2}+8x+20\right) $ \\ \midrule
$-\frac{3}{4}$ & $\frac{3}{4}$ & $\frac{e^{-x/2}}{2x^{-1/4}}\left[ 2\sqrt{x}-%
\sqrt{\pi }e^{x}\left( 2x-1\right) \mathrm{erfc}\left( \sqrt{x}\right) %
\right] $ \\ \midrule
$0$ & $\beta $ & $\sqrt{\frac{x}{\pi }}K_{\beta }\left( \frac{x}{2}\right) $
\\ \midrule
$0$ & $\frac{1}{2}$ & $e^{-x/2}$ \\ \midrule
$0$ & $\frac{3}{2}$ & $e^{-x/2}\left( 1+\frac{2}{x}\right) $ \\ \midrule
$0$ & $\frac{5}{2}$ & $e^{-x/2}\left( 1+\frac{6}{x}+\frac{12}{x^{2}}\right) $
\\ \midrule
$\frac{1}{4}$ & $-\frac{5}{4}$ & $x^{-3/4}e^{-x/2}\left( x+\frac{3}{2}%
\right) $ \\ \midrule
$\frac{1}{4}$ & $-\frac{3}{4}$ & $x^{-1/4}e^{x/2}\Gamma \left( \frac{3}{2}%
,x\right) $ \\ \midrule
$\frac{1}{4}$ & $-\frac{1}{4}$ & $x^{1/4}e^{-x/2}$ \\ \midrule
$\frac{1}{2}$ & $0$ & $x^{1/2}e^{-x/2}$ \\ \midrule
$\frac{1}{2}$ & $1$ & $x^{-1/2}e^{x/2}\Gamma \left( 2,x\right) $ \\ \midrule
$\frac{3}{4}$ & $\frac{1}{4}$ & $x^{3/4}e^{-x/2}$ \\ \midrule
$\frac{3}{4}$ & $\frac{3}{4}$ & $\frac{1}{2}x^{-1/4}e^{-x/2}\left(
2x+1\right) $ \\ \bottomrule
\end{tabular*}%
\label{TableD}%
\end{specialtable}%
\vspace{-6pt}
\begin{specialtable}[H] \centering%
\caption{The Whittaker functions $\mathrm{W}_{\kappa, \mu}$ derived for some values of
parameters $\kappa$ and $\mu$ by using~(\ref{Whittaker_def}).}%
\begin{tabular*}{\hsize}{c@{\extracolsep{\fill}}cc}
\toprule
\boldmath$\kappa $ & \boldmath$\mu $ & \boldmath$\mathrm{W}_{\kappa ,\mu }\left( x\right) $ \\
\midrule
$\frac{3}{4}$ & $\frac{5}{4}$ & $x^{-3/4}e^{x/2}\Gamma \left( \frac{5}{2}%
,x\right) $ \\ \midrule
$1$ & $\frac{1}{2}$ & $x\,e^{-x/2}$ \\ \midrule
$\frac{3}{2}$ & $0$ & $\sqrt{x}e^{-x/2}\left( x-1\right) $ \\ \midrule
$\frac{3}{2}$ & $1$ & $x^{3/2}e^{-x/2}$ \\ \midrule
$\frac{3}{2}$ & $2$ & $x^{-3/2}e^{x/2}\Gamma \left( 4,x\right) $ \\ \midrule
$2$ & $\frac{1}{2}$ & $x\left( x-2\right) e^{-x/2}$ \\ \midrule
$2$ & $\frac{3}{2}$ & $x^{2}e^{-x/2}$ \\ \midrule

$2$ & $\frac{5}{2}$ & $x^{-2}e^{x/2}\Gamma \left( 5,x\right) $ \\ \midrule
$\frac{5}{2}$ & $0$ & $\sqrt{x}e^{-x/2}\left( x^{2}-4x+2\right) $ \\ \midrule
$\frac{5}{2}$ & $1$ & $x^{3/2}e^{-x/2}\left( x-3\right) $ \\ \midrule
$\frac{5}{2}$ & $2$ & $x^{5/2}e^{-x/2}$ \\

\bottomrule
\end{tabular*}
\end{specialtable}

\begin{specialtable}[H]\ContinuedFloat
\small
\caption{{\em Cont.}}

\begin{tabular*}{\hsize}{c@{\extracolsep{\fill}}cc}
\toprule
\boldmath$\kappa $ & \boldmath$\mu $ & \boldmath$\mathrm{W}_{\kappa ,\mu }\left( x\right) $ \\
\midrule

$\frac{5}{2}$ & $3$ & $x^{-5/2}e^{x/2}\Gamma \left( 6,x\right) $ \\ \midrule
$3$ & $\frac{1}{2}$ & $e^{-x/2}x\left( x^{2}-6x+6\right) $ \\ \midrule

$3$ & $\frac{3}{2}$ & $e^{-x/2}x^{2}\left( x-4\right) $ \\ \midrule
$3$ & $\frac{5}{2}$ & $x^{3}e^{-x/2}$ \\ \midrule
$\frac{7}{2}$ & $0$ & $e^{-x/2}\sqrt{x}\left( x^{3}-9x^{2}+18x-6\right) $ \\
\midrule
$4$ & $\frac{1}{2}$ & $e^{-x/2}x\left( x^{3}-12x^{2}+36x-24\right) $ \\
\midrule
$4$ & $\frac{3}{2}$ & $e^{-x/2}x^{2}\left( x^{2}-10x+20\right) $ \\ \bottomrule
\end{tabular*}%
\label{TableDa}%
\end{specialtable}%
\vspace{-6pt}
\begin{specialtable}[H] \centering%
\caption{The Laplace transforms of the Whittaker function $\mathrm{W}_{\kappa, \mu}$ derived for some values of
parameters $\kappa$ and $\mu$.}%
\begin{tabular*}{\hsize}{c@{\extracolsep{\fill}}cc}
\toprule
\boldmath$\kappa $ & \boldmath$\mu $ & \boldmath$\mathcal{L}\left[ \mathrm{W}_{\kappa ,\mu }\left(
t\right) \right] $ \\ \midrule
$-\frac{1}{2}$ & $0$ & $\sqrt{2\pi }\left[ \frac{\mathrm{\ln }\left( \sqrt{%
4s^{2}-1}+2s\right) }{\left( 2s-1\right) ^{3/2}}+\frac{2}{\left( 1-2s\right)
\sqrt{2s+1}}\right] $ \\ \midrule
$-\frac{1}{2}$ & $1$ & $\sqrt{\frac{2\pi }{2s+1}}$ \\ \midrule
$\frac{1}{4}$ & $-\frac{5}{4}$ & $\frac{2^{1/4}\Gamma \left( \frac{1}{4}%
\right) \left( 3s+1\right) }{\left( 2s+1\right) ^{5/4}}$ \\ \midrule
$\frac{1}{4}$ & $-\frac{1}{4}$ & $\frac{\Gamma \left( \frac{5}{4}\right) }{%
\left( s+\frac{1}{2}\right) ^{5/4}}$ \\ \midrule
$0$ & $\frac{1}{2}$ & $\frac{2}{2s+1}$ \\ \midrule
$0$ & $1$ & $\frac{8s\,\mathrm{E}\left( \frac{1}{2}-s\right) -2\left(
2s+1\right) \mathrm{K}\left( \frac{1}{2}-s\right) }{4s^{2}-1}$ \\ \midrule
$\frac{1}{2}$ & $0$ & $\frac{\sqrt{2\pi }}{\left( 2s+1\right) ^{3/2}}$ \\
\midrule
$\frac{3}{4}$ & $\frac{1}{4}$ & $\frac{\Gamma \left( \frac{7}{4}\right) }{%
\left( s+\frac{1}{2}\right) ^{7/4}}$ \\ \midrule
$\frac{3}{4}$ & $\frac{3}{4}$ & $\frac{2^{3/4}\Gamma \left( \frac{3}{4}%
\right) \left( s+2\right) }{\left( 2s+1\right) ^{7/4}}$ \\ \midrule
$\frac{3}{4}$ & $\frac{5}{4}$ & $\,4\,\Gamma \left( \frac{11}{4}\right)
\,_{2}F_{1}\left( \left.
\begin{array}{c}
\frac{1}{4},\frac{11}{4} \\
\frac{5}{4}%
\end{array}%
\right\vert \frac{1}{2}-s\right) $ \\ \midrule
$1$ & $\frac{1}{2}$ & $\frac{4}{\left( 2s+1\right) ^{2}}$ \\ \midrule
$\frac{3}{2}$ & $0$ & $\frac{2\sqrt{2\pi }\left( 1-s\right) }{\left(
2s+1\right) ^{5/2}}$ \\ \midrule
$\frac{3}{2}$ & $1$ & $\frac{3\sqrt{2\pi }}{\left( 2s+1\right) ^{5/2}}$ \\
\midrule
$2$ & $\frac{1}{2}$ & $\frac{8-16s}{\left( 2s+1\right) ^{3}}$ \\ \midrule
$2$ & $\frac{3}{2}$ & $\frac{16}{\left( 2s+1\right) ^{3}}$ \\ \midrule
$\frac{5}{2}$ & $0$ & $\frac{\sqrt{2\pi }\left( 8s^{2}-16s+5\right) }{\left(
2s+1\right) ^{7/2}}$ \\ \midrule
$\frac{5}{2}$ & $1$ & $\frac{6\sqrt{2\pi }\left( 1-3s\right) }{\left(
2s+1\right) ^{7/2}}$ \\ \midrule
$\frac{5}{2}$ & $2$ & $\frac{15\sqrt{2\pi }}{\left( 2s+1\right) ^{7/2}}$ \\

\bottomrule
\end{tabular*}
\end{specialtable}

\begin{specialtable}[H]\ContinuedFloat
\small
\caption{{\em Cont.}}

\begin{tabular*}{\hsize}{c@{\extracolsep{\fill}}cc}
\toprule
\boldmath$\kappa $ & \boldmath$\mu $ & \boldmath$\mathcal{L}\left[ \mathrm{W}_{\kappa ,\mu }\left(
t\right) \right] $ \\ \midrule

$3$ & $\frac{1}{2}$ & $\frac{24\left( 1-2s\right) ^{2}}{\left( 2s+1\right)
^{4}}$ \\ \midrule
$3$ & $\frac{3}{2}$ & $\frac{32\left( 4s-1\right) }{\left( 2s+1\right) ^{4}}$
\\ \midrule
$3$ & $\frac{5}{2}$ & $\frac{96}{\left( 2s+1\right) ^{4}}$ \\ \midrule
$\frac{7}{2}$ & $0$ & $-\frac{6\sqrt{2\pi }\left(
8s^{3}-24s^{2}+15s-3\right) }{\left( 2s+1\right) ^{9/2}}$ \\ \midrule
$4$ & $\frac{1}{2}$ & $\frac{96\left( 1-2s\right) ^{3}}{\left( 2s+1\right)
^{5}}$ \\ \midrule
$4$ & $\frac{3}{2}$ & $\frac{128\left( 10s^{2}-5s+1\right) }{\left(
2s+1\right) ^{5}}$ \\ \midrule
\end{tabular*}%
\label{TableG}%
\end{specialtable}%

\section{Representations of the Wright functions \label{Appendix: Wright}}

The Wright functions $\mathrm{W}_{\alpha ,\beta }\left( x\right) $, defined
in (\ref{Wright_def}), and presented in Tables \ref{TableH} and \ref{TableHa}%
, as well as the Mainardi functions $\mathrm{F}_{\alpha }\left( x\right) $
and $\mathrm{M}_{\alpha }\left( x\right) $, defined in (\ref{Mainardi_def}),
and presented in \mbox{Tables \ref{TableJ}} and \ref{TableK}, were derived by using
the MATHEMATICA program. Only a small part of these Wright functions is known in
the  mathematical reference literature.

In the case of positive rational $\alpha =p/q$ with $p$ and $q$ positive
coprimes, applying (\ref{Sum_split})\ and (\ref{Multiplication_Gamma}), it
is possible to express the Wright function by%
\begin{eqnarray}
&&\mathrm{W}_{p/q,\beta }\left( x\right)  \label{Wright_p/q_reduction} \\
&=&\sum_{k=0}^{q-1}\frac{x^{k}}{k!\,\Gamma \left( \frac{p}{q}k+\beta \right)
}\,_{0}F_{p+q-1}\left( \left.
\begin{array}{c}
- \\
b_{0},\ldots ,b_{p-1},c_{0}^{\ast },\ldots ,c_{q-2}^{\ast }%
\end{array}%
\right\vert \frac{x^{q}}{p^{p}q^{q}}\right) ,  \nonumber
\end{eqnarray}%
where%
\begin{equation}
\begin{array}{l}
\displaystyle%
b_{j}=\frac{k}{q}+\frac{\beta +j}{p}, \\
\displaystyle%
c_{j}=\frac{k+1+j}{q},%
\end{array}
\label{Wright_p/q_coefficients}
\end{equation}%
and the set of numbers $\left\{ c_{j}^{\ast }\right\} =\left\{ c_{j}\right\}
\backslash \left\{ 1\right\} $.

For the Mainardi functions, we have the following reduction formulas for
positive rational $\alpha =p/q$ with $p$ and $q$ positive coprimes:%
\begin{eqnarray}
&&\mathrm{F}_{p/q}\left( x\right)  \label{F_p/q} \\
&=&-\frac{1}{\pi }\sum_{k=1}^{q}\frac{\left( -x\right) ^{h}}{k!}\Gamma
\left( \frac{p}{q}k+1\right) \sin \left( \pi \frac{p}{q}k\right)  \nonumber
\\
&&\,_{p}F_{q-1}\left( \left.
\begin{array}{c}
a_{0},\ldots ,a_{p-1} \\
b_{0}^{\ast },\ldots ,b_{q-2}^{\ast }%
\end{array}%
\right\vert \frac{\left( -1\right) ^{p+q}x^{q}p^{p}}{q^{q}}\right) ,
\nonumber
\end{eqnarray}%
and
\begin{equation}
\mathrm{M}_{p/q}\left( x\right) =\frac{q}{px}\mathrm{F}_{p/q}\left( x\right)
,  \label{M_p/q}
\end{equation}%
where
\begin{eqnarray*}
a_{j} &=&\frac{k}{q}+\frac{j+1}{p}, \\
b_{j} &=&\frac{k+1+j}{q},
\end{eqnarray*}%
and the set of numbers $\left\{ b_{j}^{\ast }\right\} =\left\{ b_{j}\right\}
\backslash \left\{ 1\right\} $.

\begin{specialtable}[H] \centering%
\caption{The Wright functions $\mathrm{W}_{\alpha, \beta}$ derived for some values of
parameters $\alpha$ and $\beta$ by using (\ref{Wright_def}).}%
\begin{tabular*}{\hsize}{c@{\extracolsep{\fill}}cc}
\toprule
\boldmath$\alpha $ & \boldmath$\beta $ & \boldmath$\mathrm{W}_{\alpha ,\beta }\left( x\right) $ \\
\midrule
$-1$ & $\frac{1}{2}$ & $\frac{1}{2\sqrt{\pi }\left( x+1\right) ^{3/2}}$ \\
\midrule
$-1$ & $\frac{3}{2}$ & $\frac{1}{\sqrt{\pi }\left( x+1\right) ^{1/2}}$ \\
\midrule
$-1$ & $\beta $ & $\frac{\left( x+1\right) ^{\beta -1}}{\Gamma \left( \beta
\right) }$ \\ \midrule
$-\frac{1}{2}$ & $\beta $ & $\frac{1}{\Gamma \left( \beta \right) }%
\,_{1}F_{1}\left( \left.
\begin{array}{c}
1-\beta \\
\frac{1}{2}%
\end{array}%
\right\vert -\frac{x^{2}}{4}\right) +\frac{x}{\Gamma \left( \beta -\frac{1}{2%
}\right) }\,_{1}F_{1}\left( \left.
\begin{array}{c}
\frac{3}{2}-\beta \\
\frac{3}{2}%
\end{array}%
\right\vert -\frac{x^{2}}{4}\right) $ \\ \midrule
$-\frac{1}{2}$ & $-1$ & $\frac{1}{8\sqrt{\pi }}\left[ x\left( 6-x^{2}\right)
e^{-x^{2}/4}\right] $ \\ \midrule
$-\frac{1}{2}$ & $-\frac{1}{2}$ & $\frac{1}{4\sqrt{\pi }}\left[ \left(
x^{2}-2\right) e^{-x^{2}/4}\right] $ \\ \midrule
$-\frac{1}{2}$ & $0$ & $-\frac{x\,e^{-x^{2}/4}}{2\sqrt{\pi }}$ \\ \midrule
$-\frac{1}{2}$ & $\frac{1}{2}$ & $\frac{\,e^{-x^{2}/4}}{\sqrt{\pi }}$ \\
\midrule
$-\frac{1}{2}$ & $1$ & $\mathrm{erf}\left( \frac{x}{2}\right) +1$ \\ \midrule
$-\frac{1}{2}$ & $\frac{3}{2}$ & $x\left[ \mathrm{erf}\left( \frac{x}{2}%
\right) +1\right] +\frac{2\,}{\sqrt{\pi }}e^{-x^{2}/4}$ \\ \midrule
$0$ & $-\frac{3}{2}$ & $\frac{3\,e^{x}}{4\sqrt{\pi }}$ \\ \midrule
$0$ & $-\frac{1}{2}$ & $\frac{e^{x}}{2\sqrt{\pi }}$ \\ \midrule
$0$ & $1$ & $e^{x}$ \\ \midrule
$0$ & $\beta $ & $\frac{e^{x}}{\Gamma \left( \beta \right) }$ \\ \midrule
$\frac{1}{3}$ & $\beta $ & $%
\begin{array}{l}
\frac{1}{\Gamma \left( \beta \right) }\,_{0}F_{3}\left( \left.
\begin{array}{c}
- \\
\frac{1}{3},\frac{2}{3},\beta%
\end{array}%
\right\vert \frac{x^{3}}{27}\right) +\frac{x}{\Gamma \left( \beta +\frac{1}{3%
}\right) }\,_{0}F_{3}\left( \left.
\begin{array}{c}
- \\
\frac{2}{3},\frac{4}{3},\beta +\frac{1}{3}%
\end{array}%
\right\vert \frac{x^{3}}{27}\right) \\
+\frac{x^{2}}{2\,\Gamma \left( \beta +\frac{2}{3}\right) }\,_{0}F_{3}\left(
\left.
\begin{array}{c}
- \\
\frac{4}{3},\frac{5}{3},\beta +\frac{2}{3}%
\end{array}%
\right\vert \frac{x^{3}}{27}\right)%
\end{array}%
$ \\ \midrule
$\frac{1}{2}$ & $\beta $ & $\frac{1}{\Gamma \left( \beta \right) }%
\,_{0}F_{2}\left( \left.
\begin{array}{c}
- \\
\frac{1}{2},\beta%
\end{array}%
\right\vert \frac{x^{2}}{4}\right) +\frac{x}{\Gamma \left( \beta +\frac{1}{2}%
\right) }\,_{0}F_{3}\left( \left.
\begin{array}{c}
- \\
\frac{3}{2},\beta +\frac{1}{2}%
\end{array}%
\right\vert \frac{x^{2}}{4}\right) $ \\ \midrule
$1$ & $\beta $ & $x^{\left( 1-\beta \right) /2}I_{\beta -1}\left( 2\sqrt{x}%
\right) $ \\ \midrule
$1$ & $-\frac{3}{2}$ & $\frac{\left( 4x+3\right) \cosh \left( 2\sqrt{x}%
\right) -6\sqrt{x}\sinh \left( 2\sqrt{x}\right) }{4\sqrt{\pi }}$ \\ \midrule
$1$ & $-\frac{1}{2}$ & $\frac{2\sqrt{x}\sinh \left( 2\sqrt{x}\right) -\cosh
\left( 2\sqrt{x}\right) }{2\sqrt{\pi }}$ \\ \midrule
$1$ & $0$ & $\sqrt{x}I_{1}\left( 2\sqrt{x}\right) $ \\ \midrule
$1$ & $\frac{1}{2}$ & $\frac{\cosh \left( 2\sqrt{x}\right) }{\sqrt{\pi }}$
\\ \midrule
$1$ & $1$ & $I_{0}\left( 2\sqrt{x}\right) $ \\ \midrule
$1$ & $\frac{3}{2}$ & $\frac{\sinh \left( 2\sqrt{x}\right) }{\sqrt{\pi x}}$
\\ \midrule
$1$ & $\frac{5}{2}$ & $\frac{2\sqrt{x}\cosh \left( 2\sqrt{x}\right) -\sinh
\left( 2\sqrt{x}\right) }{2\sqrt{\pi }x^{3/2}}$ \\ \bottomrule
\end{tabular*}%
\label{TableH}%
\end{specialtable}%

\begin{specialtable}[H] \centering%
\caption{The Wright functions $\mathrm{W}_{\alpha, \beta}$ derived for some values of
parameters $\alpha$ and $\beta$ by using (\ref{Wright_def}).}%
\begin{tabular*}{\hsize}{c@{\extracolsep{\fill}}cc}
\toprule
\boldmath$\alpha $ & \boldmath$\beta $ & \boldmath$\mathrm{W}_{\alpha ,\beta }\left( x\right) $ \\
\midrule
$\frac{3}{2}$ & $\beta $ & $\frac{1}{\Gamma \left( \beta \right) }%
\,_{0}F_{4}\left( \left.
\begin{array}{c}
- \\
\frac{1}{2},\frac{\beta +1}{3},\frac{\beta +2}{3},\frac{\beta }{3}%
\end{array}%
\right\vert \frac{x^{2}}{108}\right) +\frac{x}{\Gamma \left( \beta +\frac{3}{%
2}\right) }\,_{0}F_{4}\left( \left.
\begin{array}{c}
- \\
\frac{3}{2},\frac{2\beta +3}{6},\frac{2\beta +5}{6},\frac{2\beta +7}{6}%
\end{array}%
\right\vert \frac{x^{2}}{108}\right) $ \\ \midrule
$2$ & $\beta $ & $\frac{1}{\Gamma \left( \beta \right) }\,_{0}F_{2}\left(
\left.
\begin{array}{c}
- \\
\frac{\beta +1}{2},\frac{\beta }{2}%
\end{array}%
\right\vert \frac{x}{4}\right) $ \\ \midrule
$3$ & $\beta $ & $\frac{1}{\Gamma \left( \beta \right) }\,_{0}F_{3}\left(
\left.
\begin{array}{c}
- \\
\frac{\beta +1}{3},\frac{\beta +2}{3},\frac{\beta }{3}%
\end{array}%
\right\vert \frac{x}{27}\right) $ \\ \midrule
$4$ & $\beta $ & $\frac{1}{\Gamma \left( \beta \right) }\,_{0}F_{4}\left(
\left.
\begin{array}{c}
- \\
\frac{\beta +1}{4},\frac{\beta +2}{4},\frac{\beta +3}{4},\frac{\beta }{4}%
\end{array}%
\right\vert \frac{x}{256}\right) $ \\ \midrule
$5$ & $\beta $ & $\frac{1}{\Gamma \left( \beta \right) }\,_{0}F_{5}\left(
\left.
\begin{array}{c}
- \\
\frac{\beta +1}{5},\frac{\beta +2}{5},\frac{\beta +3}{5},\frac{\beta +4}{5},%
\frac{\beta }{5}%
\end{array}%
\right\vert \frac{x}{3125}\right) $ \\ \bottomrule
\end{tabular*}%
\label{TableHa}%
\end{specialtable}%
\vspace{-6pt}
\begin{specialtable}[H] \centering%
\caption{The Mainardi function $\mathrm{F}_{\alpha}$ derived for some values of
parameter $\alpha$ by using (\ref{F_p/q}).}%
\begin{tabular*}{\hsize}{c@{\extracolsep{\fill}}c}
\toprule
\boldmath$\alpha $ & \boldmath$\mathrm{F}_{\alpha }\left( x\right) $ \\ \midrule
$\frac{3}{4}$ & $%
\begin{array}{l}
\frac{x\,\Gamma \left( \frac{7}{4}\right) }{\sqrt{2}\pi }\,_{2}F_{2}\left(
\left.
\begin{array}{c}
\frac{7}{12},\frac{11}{12} \\
\frac{1}{2},\frac{3}{4}%
\end{array}%
\right\vert -\frac{27x^{4}}{256}\right) +\frac{3x^{2}}{8\sqrt{\pi }}%
\,_{2}F_{2}\left( \left.
\begin{array}{c}
\frac{5}{6},\frac{7}{6} \\
\frac{3}{4},\frac{5}{4}%
\end{array}%
\right\vert -\frac{27x^{4}}{256}\right)  \\
+\frac{x^{3}\,\Gamma \left( \frac{13}{4}\right) }{6\sqrt{2}\pi }%
\,_{2}F_{2}\left( \left.
\begin{array}{c}
\frac{13}{12},\frac{17}{12} \\
\frac{5}{4},\frac{3}{2}%
\end{array}%
\right\vert -\frac{27x^{4}}{256}\right)
\end{array}%
$ \\ \midrule
$\frac{2}{3}$ & $\frac{\sqrt{3}x}{4\pi }\left[ 2\,\Gamma \left( \frac{5}{3}%
\right) \,_{1}F_{1}\left( \left.
\begin{array}{c}
\frac{5}{6} \\
\frac{2}{3}%
\end{array}%
\right\vert -\frac{4x^{3}}{27}\right) +x\,\Gamma \left( \frac{7}{3}\right)
\,_{1}F_{1}\left( \left.
\begin{array}{c}
\frac{7}{6} \\
\frac{4}{3}%
\end{array}%
\right\vert -\frac{4x^{3}}{27}\right) \right] $ \\ \midrule
$\frac{1}{2}$ & $\frac{x\,e^{-x^{2}/4}}{2\sqrt{\pi }}$ \\ \midrule
$\frac{1}{3}$ & $3^{-1/3}x\,\mathrm{Ai}\left( 3^{-1/3}x\right) $ \\ \midrule
$\frac{1}{4}$ & $%
\begin{array}{l}
\frac{x\,\Gamma \left( \frac{5}{4}\right) }{\sqrt{2}\pi }\,_{0}F_{2}\left(
\left.
\begin{array}{c}
- \\
\frac{1}{2},\frac{3}{4}%
\end{array}%
\right\vert -\frac{x^{4}}{256}\right) -\frac{x^{2}}{4\sqrt{\pi }}%
\,_{0}F_{2}\left( \left.
\begin{array}{c}
- \\
\frac{3}{4},\frac{5}{4}%
\end{array}%
\right\vert -\frac{x^{4}}{256}\right)  \\
+\frac{x^{3}\,\Gamma \left( \frac{7}{4}\right) }{6\sqrt{2}\pi }%
\,_{0}F_{2}\left( \left.
\begin{array}{c}
- \\
\frac{5}{4},\frac{3}{2}%
\end{array}%
\right\vert -\frac{x^{4}}{256}\right)
\end{array}%
$ \\ \bottomrule
\end{tabular*}%
\label{TableJ}%
\end{specialtable}%

\vspace{-6pt}

\begin{specialtable}[H] \centering%
\caption{The Mainardi function $\mathrm{M}_{\alpha}$ derived for some values of
parameter $\alpha$ by using (\ref{M_p/q}).}%
\begin{tabular*}{\hsize}{c@{\extracolsep{\fill}}c}
\toprule
\boldmath$\alpha $ & \boldmath$\mathrm{M}_{\alpha }\left( x\right) $ \\ \midrule
$\frac{3}{4}$ & $%
\begin{array}{l}
\frac{1}{\Gamma \left( \frac{1}{4}\right) }\,_{2}F_{2}\left( \left.
\begin{array}{c}
\frac{7}{12},\frac{11}{12} \\
\frac{1}{2},\frac{3}{4}%
\end{array}%
\right\vert -\frac{27x^{4}}{256}\right) +\frac{x}{2\sqrt{\pi }}%
\,_{2}F_{2}\left( \left.
\begin{array}{c}
\frac{5}{6},\frac{7}{6} \\
\frac{3}{4},\frac{5}{4}%
\end{array}%
\right\vert -\frac{27x^{4}}{256}\right)  \\
+\frac{x^{2}}{2\,\Gamma \left( -\frac{5}{4}\right) }\,_{2}F_{2}\left( \left.
\begin{array}{c}
\frac{13}{12},\frac{17}{12} \\
\frac{5}{4},\frac{3}{2}%
\end{array}%
\right\vert -\frac{27x^{4}}{256}\right)
\end{array}%
$ \\ \midrule
$\frac{2}{3}$ & $3^{-2/3}e^{-2x^{3}/27}\left[ 3^{1/3}\mathrm{Ai}\left(
3^{-4/3}x^{2}\right) \,-3\mathrm{Ai}^{\prime }\left( 3^{-4/3}x^{2}\right) %
\right] $ \\ \midrule
$\frac{1}{2}$ & $\frac{\,e^{-x^{2}/4}}{\sqrt{\pi }}$ \\ \midrule
$\frac{1}{3}$ & $3^{2/3}\,\mathrm{Ai}\left( 3^{-1/3}x\right) $ \\ \midrule
$\frac{1}{4}$ & $%
\begin{array}{l}
\frac{\,2\sqrt{2}\Gamma \left( \frac{5}{4}\right) }{\pi }\,_{0}F_{2}\left(
\left.
\begin{array}{c}
- \\
\frac{1}{2},\frac{3}{4}%
\end{array}%
\right\vert -\frac{x^{4}}{256}\right) -\frac{x}{\sqrt{\pi }}%
\,_{0}F_{2}\left( \left.
\begin{array}{c}
- \\
\frac{3}{4},\frac{5}{4}%
\end{array}%
\right\vert -\frac{x^{4}}{256}\right)  \\
+\frac{\sqrt{2}x^{2}\,\Gamma \left( \frac{7}{4}\right) }{3\pi }%
\,_{0}F_{2}\left( \left.
\begin{array}{c}
- \\
\frac{5}{4},\frac{3}{2}%
\end{array}%
\right\vert -\frac{x^{4}}{256}\right)
\end{array}%
$ \\ \bottomrule
\end{tabular*}%
\label{TableK}%
\end{specialtable}%
\end{paracol}
\reftitle{References}

\end{document}